\documentclass[11pt]{amsart}

\usepackage{amsmath}
\usepackage{amssymb}
\usepackage{bbm}
\usepackage{pdfsync}
\usepackage{esint} 
\usepackage[pagebackref]{hyperref}
\usepackage[utf8]{inputenc}
\usepackage[T1]{fontenc}
\usepackage{mathtools}

\usepackage{mathrsfs} 

\renewcommand*{\backref}[1]{}
\renewcommand*{\backrefalt}[4]{%
	\ifcase #1 (Not cited.)%
	\or        (Cited on page~#2)%
	\else      (Cited on pages~#2)%
	\fi}

\usepackage[capitalise]{cleveref} 
\crefname{equation}{}{} 
\crefname{claim}{Claim}{Claims} 
\crefname{rem}{Remark}{Remarks} 
\crefname{coro}{Corollary}{Corollaries} 
\crefname{propo}{Proposition}{Propositions} 

\crefname{page}{p.}{pp.}

\usepackage{crossreftools}
\usepackage{enumitem}

\newlist{subtask}{enumerate}{1}     
\setlist[subtask,1]{label={\bf Step~\arabic*:},ref={\arabic*}} 
\crefname{subtaski}{Step}{Steps} 

\usepackage{aligned-overset}

\hypersetup{ 
	linktoc=page, 
	colorlinks=true,
	linkcolor={red!90!black}, 
	citecolor={blue!90!black}, 
	urlcolor=magenta,
	filecolor=teal,
}

\let\oldtocsection=\tocsection
\let\oldtocsubsection=\tocsubsection
\let\oldtocsubsubsection=\tocsubsubsection
\renewcommand{\tocsection}[2]{\hspace{0em}\oldtocsection{#1}{#2}}
\renewcommand{\tocsubsection}[2]{\hspace{1em}\oldtocsubsection{#1}{#2}}
\renewcommand{\tocsubsubsection}[2]{\hspace{2em}\oldtocsubsubsection{#1}{#2}}

\makeatletter
\newcommand\@dotsep{4.5}
\def\@tocline#1#2#3#4#5#6#7{\relax
	\ifnum #1>\c@tocdepth 
	\else
	\par \addpenalty\@secpenalty\addvspace{#2}%
	\begingroup \hyphenpenalty\@M
	\@ifempty{#4}{%
		\@tempdima\csname r@tocindent\number#1\endcsname\relax
	}{%
		\@tempdima#4\relax
	}%
	\parindent\z@ \leftskip#3\relax
	\advance\leftskip\@tempdima\relax
	\rightskip\@pnumwidth plus1em \parfillskip-\@pnumwidth
	#5\leavevmode\hskip-\@tempdima #6\relax
	\leaders\hbox{$\m@th
		\mkern \@dotsep mu\hbox{.}\mkern \@dotsep mu$}\hfill
	\hbox to\@pnumwidth{\@tocpagenum{#7}}\par
	\nobreak
	\endgroup
	\fi}
\makeatother 

\newcommand{\C}{{\mathbb C}}       
\newcommand{\M}{{\mathcal M}}       %
\newcommand{\R}{{\mathbb R}}       
\newcommand{\N}{{\mathbb N}}

\newcommand{\Z}{{\mathbb Z}}       
\newcommand{\DD}{{\mathcal D}}

\newcommand{\HH}{{\mathcal H}}

\newcommand{\QQ}{{\mathcal Q}}

\newcommand{\NN}{{\mathcal N}}
\newcommand{\RR}{{\mathcal R}}

\newcommand{\SSSS}{{\mathscr S}}

\newcommand{\SSSSmod}{\SSSS_{\rm mod}} 
\newcommand{\Smod}{S_{\rm mod}} 

\newcommand{\diam}{{\rm diam}}
\newcommand{\dist}{{\rm dist}}

\newcommand{\rf}[1]{{(\ref{#1})}}

\newcommand{\supp}{\operatorname{supp}}

\newcommand{\vv}{{\vspace{2mm}}}

\newcommand{\vmo}{{\operatorname{VMO}}}
\newcommand{\bmo}{{\operatorname{BMO}}}

\usepackage{dsfont}
\newcommand{\characteristic}{\mathbf{1}} 

\newcommand{\tree}{{\mathsf{Tree}}} 

\newcommand{\BA}{{\mathsf{BA}}} 
\newcommand{\SA}{{\mathsf{SA}}} 

\newcommand{\Stop}{{\mathsf{Stop}}}

\newcommand{\divv}{\operatorname{div}}

\newcommand{\loc}{\operatorname{loc}} 

\newcommand{\rom}[1]{%
	\textup{\uppercase\expandafter{\romannumeral#1}}%
}

\def\Xint#1{\mathchoice
	{\XXint\displaystyle\textstyle{#1}}%
	{\XXint\textstyle\scriptstyle{#1}}%
	{\XXint\scriptstyle\scriptscriptstyle{#1}}%
	{\XXint\scriptscriptstyle\scriptscriptstyle{#1}}%
	\!\int}
\def\XXint#1#2#3{{\setbox0=\hbox{$#1{#2#3}{\int}$ }
		\vcenter{\hbox{$#2#3$ }}\kern-.58\wd0}}

\def\avint{\;\Xint-}

\makeatletter
\newcommand{\doublewidetilde}[1]{{%
		\mathpalette\double@widetilde{#1}%
}}
\newcommand{\double@widetilde}[2]{%
	\sbox\z@{$\m@th#1\widetilde{#2}$}%
	\ht\z@=.8\ht\z@
	\widetilde{\box\z@}%
}
\makeatother

\usepackage{todonotes} 

\newtheorem{theorem}{Theorem}[section]
\newtheorem{lemma}[theorem]{Lemma}

\newtheorem{coro}[theorem]{Corollary}
\newtheorem{propo}[theorem]{Proposition}

\newtheorem*{claim*}{Claim}

\newtheorem*{theorem*}{Theorem}

\theoremstyle{definition}
\newtheorem{definition}[theorem]{Definition}
\newtheorem*{notation}{Notation}

\newtheorem{rem}[theorem]{Remark} 

\numberwithin{equation}{section}

\makeatletter
\newbox\ackbox
\newenvironment{acknowledgements}{%
  \global\setbox\ackbox=\vtop\bgroup
    \normalfont\Small
    \list{}{\labelwidth\z@
      \leftmargin3pc \rightmargin\leftmargin
      \listparindent\normalparindent \itemindent\z@
      \parsep\z@ \@plus\p@
      
    }%
    \item[\hskip\labelsep\scshape Acknowledgements.]%
}{%
    \endlist
  \egroup
  \ifvoid\ackbox\else
    \skip@20\p@ \advance\skip@-\lastskip
    \advance\skip@-\baselineskip \vskip\skip@
    \box\ackbox
    \prevdepth\z@ 
  \fi
}
\makeatother

\begin{document}
\title[$L^p$-solvability in locally flat unbounded domains]{$L^p$-solvability of boundary value problems for the Laplacian in locally flat unbounded domains}

\author[Ignasi Guillén-Mola]{Ignasi Guillén-Mola}
\address{Ignasi Guillén-Mola, Departament de Matem\`atiques, Universitat Aut\`onoma de Barcelona.
}
\email{ignasi.guillen@uab.cat}

\date{\today}



\begin{abstract}
    We establish the solvability of the $L^p$-Dirichlet and $L^{p^\prime}$-Neumann problems for the Laplacian for $p\in (\frac{n}{n-1}-\varepsilon,\frac{2n}{n-1}]$ for some $\varepsilon>0$ in $2$-sided chord-arc domains with unbounded boundary that is sufficiently flat at large scales and outward unit normal vector whose oscillation fails to be small only at finitely many dyadic boundary balls.
\end{abstract}

\maketitle

{
	\tableofcontents
}


\section{Introduction and main results}

Let $\Omega\subset \R^{n+1}$, $n\geq 2$, be an ADR domain (see \cref{types of domains}), and denote its surface measure by
$$
\sigma = \sigma_\Omega \coloneqq \HH^n|_{\partial\Omega}.
$$
For $1<p<\infty$, we say that the Dirichlet problem is solvable in $L^p$ (denoted as $(D_p)$ is solvable) in $\Omega$ if there is $C_{D_p}\geq 1$ such that for any given $f\in L^p(\sigma)$, there exists $u:\Omega\to \R$ satisfying
\begin{equation}\label{dirichlet bvp}
\begin{cases*}
    \Delta u = 0 \text{ in } \Omega,\\
        \NN u \in L^p(\sigma),\\
    u|_{\partial\Omega}^{\rm nt} = f, \text{ }\sigma\text{-a.e.}, 
\end{cases*}
\end{equation}
and
\begin{equation}\label{lp norm nontang of dirichlet solution}
\|\NN u\|_{L^p(\sigma)} \leq C_{D_p} \|f\|_{L^p(\sigma)}.
\end{equation}
Here $\NN$ is the nontangential operator (see \rf{nontangential maximal function}) and $u|_{\partial\Omega}^{\rm nt}$ is the nontangential limit (see \rf{eq:nontan limit}). We say that the Neumann problem is solvable in $L^p$ (denoted as $(N_p)$ is solvable) in $\Omega$ if there is $C_{N_p}\geq 1$ such that for any given $f\in L^p(\sigma)$, there exists $u:\Omega\to\R$ satisfying
\begin{equation}\label{neumann bvp}
\begin{cases*}
    \Delta u = 0 \text{ in } \Omega,\\
        \NN(\nabla u) \in L^p(\sigma),\\
    \partial_{\nu_\Omega} u = f, \text{ }\sigma\text{-a.e.},
\end{cases*}
\end{equation}
and
\begin{equation}\label{lp norm nontang of neumann solution}
\|\NN (\nabla u)\|_{L^p (\sigma)} \leq C_{N_p} \|f\|_{L^p (\sigma)}.
\end{equation}
Here $\partial_{\nu_\Omega}$ is the interior nontangential derivative (see \rf{def:interior nontang derivative}).

We note that the definition of solvability for the Dirichlet and Neumann problems may vary slightly across the literature. We will not distinguish between these variations in the subsequent articles mentioned in the introduction. Furthermore, while several of the results below were originally established for more general divergence-form operators in their respective papers, we restrict our historical overview to the harmonic case.

In ADR domains, it is well known that $(D_p)$ is solvable for some $1<p<\infty$ if and only if the harmonic measure $\omega$ is locally in weak-$A_\infty(\sigma)$. Specifically, there exist constants $C,s>0$ such that for every ball $B=B(x,r)$ with $x\in\partial\Omega$ and $r<\diam(\partial\Omega)/4$, there holds $\omega^p (E)\leq C (\sigma(E)/\sigma(B))^s \omega^p (2B)$ for any $p\in \Omega\setminus 4B$ and any Borel set $E\subset B$. We refer the reader to \cite{Hofmann-Le-2018,Hofmann-2019-survey}. For ADR domains satisfying the corkscrew condition (see \cref{def:corks balls}), $(D_p)$ is solvable if and only if the harmonic measure satisfies a weak reverse Hölder inequality with exponent $p^\prime\coloneqq p/(p-1)$ (the Hölder conjugate of $p$), see, for instance, \cite[Theorem 9.2]{mourgoglou-tolsa-2024-Duke-regularityprobleminroughdomains}, \cite[Proposition 2.14]{MourgoglouPoggiTolsa2025}, and the references therein. Consequently, if $(D_p)$ is solvable for some $1<p<\infty$, then Gehring's lemma guarantees the existence of $\varepsilon>0$ such that $(D_q)$ is solvable for all $p-\varepsilon\leq q<\infty$.

The solvability of the Dirichlet and Neumann problems is a long-standing and active area of research. In 1963, Lavrent'ev \cite{Lavrentiev1963} showed that in bounded planar simply connected chord-arc domains (see \cref{def:CAD}), the harmonic measure is locally in weak-$A_\infty(\sigma)$, implying the solvability of $(D_p)$ for some $1<p<\infty$. However, Jerison \cite{Jerison1983} later proved in 1983 that for every $1<p<\infty$, there exists a planar chord-arc domain where $(D_p)$ fails to be solvable. A significant breakthrough came in 1977 with Dahlberg \cite{Dahlberg1977}, who proved that in any bounded Lipschitz domain $\Omega$, its harmonic measure satisfies a reverse Hölder inequality with exponent $2$. Consequently, there exists $\varepsilon_\Omega>$ such that $(D_p)$ is solvable for all $2-\varepsilon_\Omega\leq p<\infty$. This result is sharp: for every $\varepsilon>0$ there is a Lipschitz domain $\Omega_\varepsilon$ where $(D_{2-\varepsilon})$ is not solvable, see \cite[pp.~153-54]{Kenig1986-BeijingLectureNotes}. In 1978, Fabes, Jodeit and Rivière \cite{Fabes1978} employed double and single layer potentials (see \cref{sec:double layer potential,sec:single layer potential}) to establish the solvability of both $(D_p)$ and $(N_p)$ for all $1<p<\infty$ in bounded $C^1$ domains. Notably, Dahlberg \cite{Dahlberg1979} had already proven the solvability of $(D_p)$ for all $1<p<\infty$ in bounded $C^1$ domains without using layer potentials, though this was published later in 1979.\footnote{Fabes, Jodeit and Rivière cite a technical report of \cite{Dahlberg1979} in \cite{Fabes1978}.} Subsequent work by Jerison and Kenig \cite{Jerison1980,Jerison-Kenig-1981-Dirichlet} simplified the proofs of Dahlberg's results through the application of the so-called Rellich identity \cite{Rellich1940}. Furthermore, in \cite{Jerison-Kenig-1981-Neumann}, they resolved $(N_2)$ in bounded Lipschitz domains.

Briefly speaking, the layer potential approach relies on the invertibility of the operators $\frac{1}{2}Id + K$ and $-\frac{1}{2}Id + K^*$ in $L^p(\sigma)$ (here $K$ is the double layer potential and $K^*$ is its adjoint) and gives an ``explicit'' solution to the Dirichlet and Neumann problems, respectively. In $C^{1,\alpha}$ domains, $\alpha>0$, it is not difficult to see that the double layer potential is compact in $L^p(\sigma)$ and (by Fredholm theory) that those operators are invertible, for all $1<p<\infty$. For a proof of this, see the lecture notes \cite{Dahlberg-Kenig-HAPE-Lecturenotes1985,Kenig1986-BeijingLectureNotes}. While this argument does not directly extend to $C^1$ domains, Fabes, Jodeit and Rivière \cite{Fabes1978} succeeded in showing that in bounded $C^1$ domains, the double layer potential $K$ is compact in $L^p(\sigma)$, and both $\frac{1}{2}Id + K$ and $-\frac{1}{2}Id + K^*$ are invertible in $L^p(\sigma)$ for all $1<p<\infty$. For Lipschitz domains, however, the compactness of the double layer potential generally fails, as shown by Fabes, Jodeit and Lewis in \cite{Fabes-Jodeit-Lewis-1977}. Nevertheless, using the Rellich identity mentioned above, Verchota \cite{Verchota1984} established in 1984 the invertibility of $\frac{1}{2}Id + K$ and $-\frac{1}{2}Id + K^*$ in $L^2(\sigma)$, thereby recovering the solvability of $(D_2)$ and $(N_2)$ for bounded Lipschitz domains, as originally shown in \cite{Dahlberg1977,Jerison-Kenig-1981-Neumann} respectively.

In 1987, Dahlberg and Kenig \cite{Dahlberg-Kenig-1987} established that for every bounded Lipschitz domain $\Omega$, there exists $\varepsilon_\Omega>0$ such that $(N_p)$ is solvable for all $1<p\leq 2+\varepsilon_\Omega$, and in fact, they showed that its solution, as well as the solution of $(D_{p^\prime})$ in Dahlberg's result \cite{Dahlberg1977}, can be obtained using the method of layer potentials. As in the Dirichlet problem, this range of solvability is sharp: for every $\varepsilon>0$, there exists a Lipschitz domain $\Omega_\varepsilon$ for which $(N_{2+\varepsilon})$ fails to be solvable, see \cite[pp.~153-54]{Kenig1986-BeijingLectureNotes}.

Dahlberg's result \cite{Dahlberg1977} was extended to chord-arc domains (see \cref{def:CAD}) independently by David and Jerison \cite{David-Jerison-1990-LipschitzApprox}, and Semmes \cite{Semmes-1990-AnalysisVSGeometry} in 1990. They proved that for any bounded chord-arc domain, there exists $1<p<\infty$ such that the harmonic measure satisfies a reverse Hölder inequality with exponent $p$. A complete geometric characterization came in 2020 when Azzam, Hofmann, Martell, Mourgoglou, and Tolsa \cite{Azzam-Hofmann-Martell-Mourgoglou-Tolsa-INVENTIONES-2020} identified the class of ADR domains with interior corkscrews where $(D_p)$ is solvable for some $1<p<\infty$. These are precisely domains having interior big pieces of chord-arc domains (IBPCAD), as defined in \cite[Definition 2.12]{Azzam-Hofmann-Martell-Mourgoglou-Tolsa-INVENTIONES-2020}.

Let us roughly introduce the regularity boundary value problem for an ADR domain $\Omega$. When well-defined, we say that the regularity problem is solvable in $L^p$ (denoted as $(R_p)$ is solvable) in $\Omega$ if for any given $f\in C^{0,1}(\partial\Omega)$, there is a harmonic function $u:\Omega\to \R$ such that $u|_{\partial\Omega}^{\rm nt} = f$ holds $\sigma$-a.e.\ on $\partial\Omega$ and $\|\NN (\nabla u)\|_{L^p(\sigma)}\lesssim \|\nabla f\|_{L^p(\sigma)}$, where the implicit constant does not depend on $f$. The exact definition may vary slightly across the literature; see \cite[Definition 5.35]{HofmannSparrius-2025} and \cite[Definition 1.1]{MourgoglouPoggiTolsa2025} for technical formulations.

We observe that the regularity problem is closely connected to the Dirichlet problem: for bounded domains with the corkscrew condition and having UR boundary (see \cref{def:ur set}), Mourgoglou, Poggi and Tolsa proved in \cite[Theorem 1.13]{MourgoglouPoggiTolsa2025} that $(R_p)\iff(D_{p^\prime})$ for all $p\in (1,\infty)$, see also \cite[Theorems 1.2 and 1.6]{mourgoglou-tolsa-2024-Duke-regularityprobleminroughdomains}. One might expect similar results for $2$-sided chord-arc domains with unbounded boundaries. However, to the best of our knowledge, such results have not been established in the literature. It is quite likely that the arguments in \cite{mourgoglou-tolsa-2024-Duke-regularityprobleminroughdomains,MourgoglouPoggiTolsa2025} could be extended to domains with unbounded boundaries, but this extension has not yet been documented.

The question of whether $(N_p)$ or $(R_p)$ is solvable for some $1<p<\infty$ in chord-arc domains was first posed by Kenig \cite[Problem 3.2.2]{Kenig1994} in 1994 and later reintroduced by Toro \cite[Question 2.5]{Toro-2010ICM} at the ICM 2010. While the regularity part of this question was recently resolved in 2023 by Mourgoglou, Poggi, and Tolsa in \cite[Corollary 1.16]{MourgoglouPoggiTolsa2025}, the Neumann part remains an open problem.

Although the solvability of the Neumann problem remains open for chord-arc domains, significant progress has been made in understanding its extrapolation properties under the assumption that $(D_{p^\prime})$ or/and $(R_p)$ is solvable for some $1<p<\infty$. As noted earlier, in bounded chord-arc domains (which have UR boundary by \cref{thm:2-sided corkscrew + ADR implies UR boundary}), the equivalence $(R_p)\iff(D_{p^\prime})$ holds for any $1<p<\infty$. In 1993, Kenig and Pipher \cite[Theorem 6.3]{Kenig-Pipher-1993} showed the extrapolation $(D_{p^\prime})+(N_p)\implies (N_q)$ for all $1<q<p+\varepsilon$, for some $\varepsilon>0$, for bounded Lipschitz domains. In 2024, Feneuil and Li \cite[Corollary 1.22]{Feneuil-Li-2024-lppoissonneumannproblemrelation-arxiv} extended this result to bounded chord-arc domains. In the same year, Mourgoglou and Tolsa showed in \cite[Theorem 1.1]{Mourgoglou-Tolsa-2024-solvabilityneumannproblemelliptic-arxiv} that $(N_p)$ is solvable for a fixed $p\in (1,2)$, whenever $\Omega$ is a bounded chord-arc domain such that $(R_q)$ is solvable for some $q>p$, $\partial\Omega$ supports a weak $p$-Poincaré inequality, and $\Omega$ has very big pieces of chord-arc superdomains for which $(N_q)$ is solvable. Most recently in 2025, Hofmann and Sparrius proved in \cite[Theorem 5.54]{HofmannSparrius-2025} that $(N_p)+(R_p)+(D_{p^\prime}) \implies (N_q)$ for all $1<q<p$ in $2$-sided chord-arc domains with unbounded boundary.

Given the strong connection to the Neumann problem, we briefly address the extrapolation properties of the regularity problem. The equivalence $(R_p)\iff (D_{p^\prime})$ mentioned above, combined with the well-known extrapolation of solvability for the Dirichlet problem, immediately yields extrapolation of solvability of the regularity problem. A more subtle endpoint case was resolved in 2025 by Gallegos, Mourgoglou, and Tolsa, who proved solvability extrapolation for $(R_1)$ (even for $(R_{1-\varepsilon})$) in ADR domains satisfying the interior corkscrew condition, see \cite[Theorems 1.3 and 1.6]{Gallegos-Mourgoglou-Tolsa-ExtrapolationRegularity2025}.

Let us now return to our main discussion of the solvability of Dirichlet and Neumann problems. In 2010, Hofmann, M. Mitrea and Taylor studied the solvability of the Dirichlet and Neumann problems (among others) in bounded $\delta$-regular SKT (Semmes-Kenig-Toro) domains. Roughly speaking, a domain is $\delta$-regular SKT if it is Reifenberg flat (with small enough constant) and its geometric measure theoretic outward unit normal vector $\nu$ (see \cref{rem:def of geometric unit vector}) satisfies $\dist_{\bmo(\sigma)}(\nu,\vmo(\sigma))\leq\delta$. The precise definition of $\delta$-regular SKT domains can be found in \cite[Definition 4.9]{Hofmann2010}. In \cite[Section 5]{Hofmann2010}, the authors proved that for any $1<p<\infty$, both $(D_p)$ and $(N_{p^\prime})$ are solvable in $\delta$-regular SKT domains when $\delta$ is small enough. Consequently, $(D_p)$ and $(N_{p^\prime})$ are solvable for all $1<p<\infty$ in regular SKT domains, that is, domains that are $\delta$-regular SKT for all $\delta>0$, see \cite[Definition 4.8]{Hofmann2010}.

The approach in \cite{Hofmann2010} employs layer potentials. For any $1<p<\infty$ and a bounded regular SKT domain, the authors show in \cite[Theorem 4.36]{Hofmann2010} that the double layer potential $K$ is compact. More precisely, for bounded $\delta$-regular SKT domains with sufficiently small $\delta=\delta(p)>0$, \cite[Theorem 4.36]{Hofmann2010} establishes that $K$ is close enough in $L^p(\sigma)$ norm to the set of compact operators (see \cref{def:compact operator}). This property nevertheless allows the application of Fredholm theory to prove the invertibility of both $\frac{1}{2}Id + K$ and $-\frac{1}{2}Id+K^*$ from the injectivity of $\frac{1}{2}Id + K^*$ in $L^2(\sigma)$, see \cite[Proposition 5.11]{Hofmann2010}.

Recently in 2022, Marín, Martell, D. Mitrea, I. Mitrea and M. Mitrea \cite[Chapter 6]{Marin-Martell-Mitrea-x3} showed the $L^p$ solvability of the Dirichlet and Neumann problems (among others) for $2$-sided chord-arc domains with unbounded boundary, under the $\bmo$ smallness condition $\|\nu\|_*<\delta$ where $\delta>0$ is sufficiently small depending on $1<p<\infty$. They proved that the $L^p(\sigma)$ norm of $K$ tends to zero as $\delta\to 0$. This enabled them to see that both $\frac{1}{2}Id + K$ and $-\frac{1}{2}Id + K^*$ are invertible via Neumann series, thereby solving the Dirichlet and Neumann problems for any $1<p<\infty$. More specifically, the smallness of $\delta>0$ depends on the dimension, the chord-arc parameters of the domain and $p$. We emphasize that the results in \cite{Marin-Martell-Mitrea-x3} hold for weighted $L^p$ spaces and systems of divergence-form operators with constant coefficient matrices.

Throughout this work, we will work in $\R^{n+1}$ with $n\geq 2$, although similar results may hold in the planar case $n=1$. We focus on domains defined as follows:

\begin{definition}[$\delta$-$(s,S;R)$ domain]\label{def:final assumptions of domain}
    A domain $\Omega\subset \R^{n+1}$ is called $\delta$-$(s,S;R)$ domain if it is a $2$-sided chord-arc domain (see \cref{def:CAD}) with unbounded boundary and there exist $\delta>0$, scales $S>s>0$ and a radius $R>0$ such that the outward unit normal vector $\nu$ (see \cref{rem:def of geometric unit vector}) satisfies
    $$
    \avint_{B(x,r)} |\nu(z)-m_{B(x,r)}\nu|\, d\sigma(z) \leq\delta, \text{ provided}
    \begin{cases*}
        x\in\partial\Omega\setminus B_R (0) \text{ and } r\in (0,\infty),\text{ or}\\
        x\in B_R (0) \cap \partial\Omega \text{ and } r\not\in (s,S),
    \end{cases*}
    $$
    (here $m_{B(x,r)}\nu=\frac{1}{\sigma(B(x,r))}\int_{B(x,r)} \nu \, d \sigma$) and for all $x\in\partial\Omega$ and $r\geq S$, there holds
    $$
    \beta_{\infty,\partial\Omega} (B(x,r)) \coloneqq \inf_{n\text{-plane } L\ni x} \sup_{y\in\partial\Omega\cap B(x,r)} \frac{\dist(y,L)}{r}\leq \delta,
    $$
    where the infimum is taken over all $n$-planes $L\subset \R^{n+1}$ through $x$.
\end{definition}

That is, these are domains whose failure of sufficient flatness is limited to finitely many dyadic boundary balls. We note that the domains studied in \cite{Marin-Martell-Mitrea-x3} satisfy the first condition in \cref{def:final assumptions of domain} for every $x\in\partial\Omega$ and every scale $r\in (0,\infty)$, implying $\beta_{\infty,\partial\Omega} (B(x,r))\lesssim \delta^{\frac{1}{2n}}$ for all $x\in\partial\Omega$ and all $r\in (0,\infty)$, see \cite[Theorem 2.2]{Marin-Martell-Mitrea-x3}.

In this paper, we study the $L^p$-solvability of the Dirichlet and Neumann problems in $\delta$-$(s,S;R)$-domains, and we also provide uniqueness results for \rf{dirichlet bvp} and \rf{neumann bvp} respectively. In the subsequent results, $\DD$ denotes the interior double layer potential (see \rf{double layer potential}), $K$ the (boundary) double layer potential (see \rf{pv double layer potential}), $K^*$ the adjoint of $K$, and $\SSSSmod$ the modified interior single layer potential (see \rf{single layer potential}).

\begin{theorem}[$L^p$-Dirichlet problem]\label{thm:lp dirichlet problem}
    Let $\Omega\subset \R^{n+1}$ be a $\delta$-$(s,S;R)$ domain. There is $\varepsilon_D=\varepsilon_D (n,\text{CAD})\in (0,\frac{1}{n-1})$ such that for every $p_0\in (\frac{n}{n-1}-\varepsilon_D,\frac{2n}{n-1}]$ there exists $\delta_0=\delta_0(n,p_0,\text{CAD})>0$ such that if $\delta\leq \delta_0$, then $\frac{1}{2}Id + K$ is invertible in $L^{p_0} (\sigma)$ and given $f\in L^{p_0}(\sigma)$, the function
    \begin{equation}\label{eq:solution of lp dirichlet problem}
    u \coloneqq \DD\left(\left(\frac{1}{2}Id + K\right)^{-1} f\right),
    \end{equation}
    is the solution of $(D_{p_0})$. Furthermore, there exists $\varepsilon>0$ such that $(D_p)$ is solvable for all $p\in (p_0-\varepsilon,\infty)$, and the solution of $(D_p)$ is the unique solution of the Dirichlet problem \rf{dirichlet bvp}.
\end{theorem}

\begin{theorem}[$L^p$-Neumann problem]\label{thm:lp neumann problem}
    Let $\Omega\subset \R^{n+1}$ be a $\delta$-$(s,S;R)$ domain. There is $\varepsilon_N = \varepsilon_N (n,\text{CAD})>0$ such that for every $p\in [\frac{2n}{n+1},n+\varepsilon_N)$ there exists $\delta_0=\delta_0(n,p,\text{CAD})>0$ such that if $\delta\leq \delta_0$, then $-\frac{1}{2}Id + K^*$ is invertible in $L^p (\sigma)$, $(N_p)$ is solvable and given $f\in L^p(\sigma)$, the function
    \begin{equation}\label{eq:solution of lp neumann problem}
    u \coloneqq \SSSSmod\left(\left(-\frac{1}{2}Id + K^*\right)^{-1} f\right),
    \end{equation}
    is the unique (modulo constants) solution of $(N_p)$. Furthermore, it is the unique (modulo constants) solution of the Neumann problem \rf{neumann bvp}.
\end{theorem}

\Cref{thm:lp dirichlet problem,thm:lp neumann problem} mainly follow from the invertibility of $\frac{1}{2}Id + K$ and $-\frac{1}{2}Id + K^*$, respectively. We emphasize that the existence of $\varepsilon_D>0$ in \cref{thm:lp dirichlet problem} and $\varepsilon_N>0$ in \cref{thm:lp neumann problem} guarantees the solvability of both $(D_2)$ and $(N_2)$.

In contrast with \cite{Marin-Martell-Mitrea-x3}, the invertibility of $\frac{1}{2}Id+K$ in $L^{p^\prime}(\sigma)$ and $-\frac{1}{2}Id+K^*$ in $L^p(\sigma)$ on which \cref{thm:lp dirichlet problem,thm:lp neumann problem} rely holds only for a bounded range of exponents, $p\in[\frac{2n}{n+1},n+\varepsilon)$. This is the range for which we can prove the injectivity of $\pm\frac{1}{2}Id+K^*$ in $L^p(\sigma)$ (see \cref{sec:+-id+k* are injective}), which is governed by the uniqueness (modulo constants) of harmonic functions with $\NN(\nabla u)\in L^p(\sigma)$ (see \cref{sec:uniqueness of neumann problem}). Essentially, the lower endpoint is where the gradient of the solution becomes square-integrable, so that the energy-based regularity estimates apply, and the upper endpoint is where the solution still decays at infinity, forcing it to be constant. We note that, in the bounded case, it would suffice to prove this injectivity at the single endpoint $p=\frac{2n}{n+1}$, with the full range $[\frac{2n}{n+1},\infty)$ then following by Hölder's inequality.

Let us briefly discuss the extrapolation of solvability. For the Dirichlet problem, the range extends up to $p=\infty$ by Gehring's lemma, whereas for the Neumann problem no comparable extrapolation is available in our setting. While for $2$-sided chord-arc domains with unbounded boundary there holds $(N_p)+(R_p)+(D_{p^\prime}) \implies (N_q)$ for all $1<q<p$, our results for $\delta$-$(s,S;R)$ domains in \cref{thm:lp dirichlet problem,thm:lp neumann problem} establish $(D_{p^\prime})$ and $(N_p)$ solvability for a fixed $p\in [2n/(n+1),n+\varepsilon)$ ($\delta$ is sufficiently small), but not the solvability of $(R_p)$. Since we have not yet established the solvability of $(R_p)$, we cannot consequently derive that $(N_q)$ is solvable for all $1<q\leq p$.

Let us outline the proof strategy for \cref{thm:lp dirichlet problem,thm:lp neumann problem}. Unlike the approach in \cite{Marin-Martell-Mitrea-x3} discussed above, our setting with $\delta$-$(s,S;R)$-domains presents a key difference: the double layer potential $K$ does not a priori have small norm, and thus we cannot directly derive the invertibility of $\frac{1}{2}Id + K$ and $-\frac{1}{2}Id + K^*$ via Neumann series. As in \cite[Theorem 4.36]{Hofmann2010}, in the following result (one of the main points in this article), we show that $K$ and its adjoint $K^*$ are close to the set of compact operators (see \cref{def:compact operator}). This enables us to employ Fredholm theory to characterize their invertibility.

\begin{theorem}\label{thm:double layer potential is almost compact}
    Let $\Omega\subset \R^{n+1}$ be a $\delta$-$(s,S;R)$ domain, $1<p<\infty$, and let $p^\prime = p/(p-1)$ be its Hölder conjugate exponent. For all $\varepsilon > 0$ there exists $\delta_0 = \delta_0(\varepsilon,p,\text{CAD},n)>0$ such that if $\delta\leq \delta_0$, then there is a compact operator $T:L^p(\sigma)\to L^p(\sigma)$ such that $\|K-T\|_{L^p (\sigma)\to L^p(\sigma)} = \|K^*-T^*\|_{L^{p^\prime} (\sigma) \to L^{p^\prime}(\sigma)}<\varepsilon$. 
\end{theorem}

Here, $T^*:L^{p^\prime}(\sigma)\to L^{p^\prime}(\sigma)$ is the adjoint of $T$, which is compact by Schauder's theorem. To prove this theorem we truncate the double layer potential $K$ at small, intermediate (``close'' and ``far'' from $B_R(0)$) and large scales, see \rf{full truncation of double layer potential}. The operator on the ``close'' intermediate scales turns out to be compact (\cref{trucated and intermediate scales is compact}). The operators at small and ``far'' intermediate scales have small norm (\cref{lemma:double layer potentials in small scales has small norm,lemma:double layer potentials far from bad ball has small norm}), by using the already known behavior of the double layer potential at these scales (\cref{theorem:small scales Lp}) via Semmes' decomposition \cref{Semmes' decomposition}. One of the main difficulties in this article is establishing the small norm of the operator at large scales (\cref{large scale bounded by bmo norm}), which we address in \cref{sec:proof of large scale bounded by bmo}.

Let us explain why the large-scale estimate is the delicate point. At small and ``far'' intermediate scales, Semmes' decomposition shows that, away from an exceptional set of small measure, $\partial\Omega$ coincides with a Lipschitz graph with small constant, which makes the corresponding truncated operators have small norm. At large scales this is no longer available: $K$ has no small norm a priori and the flatness of $\partial\Omega$ may fail at finitely many bad balls. We overcome this by comparing $\partial\Omega$ with a Lipschitz graph that is now only close to it rather than coincident with it, built through a stopping-time argument in the spirit of David-Semmes (\cref{sec:construction of Lip graph}). This closeness is ensured by the large-scale smallness of $\beta_{\infty,\partial\Omega}$ and $\|\nu\|_*$.

Combining \cref{thm:double layer potential is almost compact} with the Fredholm alternative \cref{fredholm alternative}, in \cref{reduction:all equivalent}, for $\lambda \in \R\setminus \{0\}$ we establish several equivalent conditions in order to see that $\lambda Id + K$ (respectively $\lambda Id + K^*$) is invertible in $L^p(\sigma)$ (respectively $L^{p^\prime}(\sigma)$). Furthermore, the invertibility of $\lambda Id + K$ in $L^p(\sigma)$ and $\lambda Id + K^*$ in $L^{p^\prime}(\sigma)$ is shown to be equivalent to the injectivity of $\lambda Id + K^*$ in $L^{p^\prime}(\sigma)$.

As previously noted, the solvability of the Dirichlet and Neumann problems mainly relies on the invertibility of $\frac{1}{2}Id + K$ and $-\frac{1}{2}Id + K^*$. The equivalences discussed in the previous paragraph allow us to reduce this invertibility problem to showing the injectivity of $\pm\frac{1}{2}Id + K^*$ in $L^p(\sigma)$ for some range of $p$. Using the modified single layer potential (see \cref{sec:single layer potential}) and Moser estimates for harmonic functions with vanishing Dirichlet or Neumann boundary conditions, in \cref{sec:+-id+k* are injective} we prove that $\pm\frac{1}{2}Id + K^*$ are indeed injective in $L^p(\sigma)$ when $p\in [2n/(n+1),n+\varepsilon)$, for some $\varepsilon>0$.

In \cref{sec:proof of dirichlet and neumann problems}, we establish the $L^p$ solvability of the Dirichlet and Neumann problems in $\delta$-$(s,S;R)$ domains in \cref{thm:lp dirichlet problem,thm:lp neumann problem}. This mainly follows from the invertibility explained above, and also from some already-known properties of the single and double layer potentials. The uniqueness results for \rf{dirichlet bvp} and \rf{neumann bvp} are proved in \cref{sec:uniqueness of neumann problem,sec:uniqueness of dirichlet problem}.

A natural question is whether \cref{thm:lp dirichlet problem,thm:lp neumann problem,thm:double layer potential is almost compact} generalize to second-order divergence form operators $\divv(A\nabla\cdot)$ for uniformly elliptic matrices $A$ with Dini mean oscillation coefficients, in the spirit of recent work \cite{Molero2021} by Molero, Mourgoglou, Puliatti, and Tolsa.

\section{Preliminaries and definitions}

\subsection{Notation}\label{sec:notation}

\begin{itemize}
	\item We use $c,C\geq 1$ to denote constants that may depend only on the dimension and the constants appearing in the hypotheses of the results, and whose values may change at each occurrence.
	
	\item We write $a\lesssim b$ if there exists a constant $C\geq 1$ such that $a\leq Cb$, and $a\approx b$ if $C^{-1} b \leq a \leq C b$.
	
	\item If we want to stress the dependence of the constant on a parameter $\eta$, we write $a\lesssim_\eta b$ or $a\approx_\eta b$ meaning that $C=C(\eta)=C_\eta$.
	
	\item The ambient space is $\R^{n+1}$ with $n\geq 2$. We denote by $m$ the Lebesgue measure on $\R^{n+1}$.
	
	\item The diameter of a set $E\subset \R^{n+1}$ is denoted by $\diam\, E$. We allow $\diam\, E = \infty$ if $E$ is unbounded.
	
	\item We denote by $B_r (x)$ or $B(x,r)$ the open ball with center $x$ and radius $r>0$. We denote $B_r \coloneqq B_r (0)$.

        \item Given a domain $\Omega$, we denote the boundary ball centered at $x\in\partial\Omega$ with $r>0$ by $\Delta(x,r) = \Delta_r (x) \coloneqq B(x,r)\cap \partial\Omega$.
	
	\item Given a ball $B$, we denote by $r_B$ or $r(B)$ its radius, and by $c_B$ or $c(B)$ its center. Analogously, $r_\Delta$ or $r(\Delta)$ and $c_\Delta$ or $c(\Delta)$ for a boundary ball $\Delta$.

        \item Given a ball $B$ and $t>1$, $tB\coloneqq B(c_B, tr_B)$. Analogously, $t\Delta = \Delta(c_\Delta,tr_\Delta)$.

        \item We denote by $Q_r(x)$ or $Q(x,r)$ the open cube with center $x$ and side length $2r$, i.e., $Q_r(x) = Q(x,r) = \{y\in \R^{n+1} : |y_i-x_i|<r\text{ for all }1\leq i \leq n+1\}$.
    
        \item Given a cube $Q$, we denote by $\ell(Q)$ its side length, and by $c_Q$ or $c(Q)$ its center. That is, $Q=Q(c_Q, \ell(Q)/2)$.

        \item Given a cube $Q$ and $t>1$, $tQ \coloneqq Q(c_Q, t\ell(Q)/2)$, that is, $c_{tQ} = c_Q$ and $\ell(tQ)=t\ell(Q)$.
	
	\item We say that a function $f$ is Hölder continuous with exponent $\alpha \in (0,1]$ in a set $U$, or briefly $C^{0,\alpha} (U)$, if there exists a constant $C_\alpha >0$ (called the Hölder seminorm) such that $|f(x)-f(y)| \leq C_\alpha |x-y|^{\alpha}$ for all $x, y \in U$. For $\alpha\in (0,1)$ we write $C^\alpha$ instead of $C^{0,\alpha}$. When $\alpha = 1$ we say that $f$ is ``Lipschitz continuous'', and we call $C_L \coloneqq C_1$ its Lipschitz constant. We also say that $f$ is $\kappa$-Lipschitz if $C_L \leq \kappa$.
	
	\item We denote the characteristic function of a set $E$ by $\characteristic_E$.
	
	\item Denote $\DD(\R^{n+1})$ the standard dyadic grid. That is, $\DD(\R^{n+1}) = \bigcup_{k\in \Z} \DD_k (\R^{n+1})$ where $\DD_k (\R^{n+1})$ is the collection of all cubes of the form
	$$
	\{ x \in \R^{n+1} : m_i 2^{-k} \leq x_i < (m_i + 1)2^{-k} \text{ for } i=1, \ldots, n+1\},
	$$
	where $m_i \in \Z$.
	
	\item Given $t>0$ and a set $E\subset \R^{n+1}$, we write $U_t (E) \coloneqq \{x\in \R^{n+1} : \dist (x,E) < t \}$ for the $t$-neighborhood of $E$.

    \item For $0\leq s<\infty$, $\HH^s$ denotes the $s$-dimensional Hausdorff measure.

    \item $\mu|_{S}$ denotes the restriction of the measure $\mu$ to a set $S\subset \R^{n+1}$ defined as $\mu|_{S} (E) \coloneqq \mu(S\cap E)$ for $E\subset \R^{n+1}$.

    \item Given a measure $\mu$ and a set $E$, if $\mu(E)\not=0$ then we denote $m_E f \coloneqq \avint_E f\, d\mu \coloneqq \frac{1}{\mu(E)} \int_E f\, d\mu$ for $f\in L^1_{\loc} (\mu)$.

    \item Let $E\subset \R^{n+1}$ be an ADR closed set (see \cref{types of domains}) and let $\mu\coloneqq \HH^n|_{E}$. Given $f\in L^2_{\loc} (\mu)$, $x\in E$, $R>0$, we set
    $$
    \|f\|_* (B(x,R)) \coloneqq \sup_{B\subset B(x,R)} \left(\avint_B |f(z)-m_B f|^2 \, d\mu(z)\right)^{1/2},
    $$
    where the supremum is taken over all balls $B$ centered at $E$ included in $B(x,R)$, and $m_B f\coloneqq \avint_B f\, d\mu$.

    \item Given an ADR domain $\Omega$ satisfying the $2$-sided corkscrew condition (see \cref{types of domains}) with outward unit normal $\nu_\Omega$ (see \cref{rem:def of geometric unit vector}), the tangential gradient of a Lipschitz function $f$ in $\partial\Omega$ is
    $$
    \nabla_t f (y) \coloneqq \nabla \widetilde f(y) - \langle\nabla \widetilde f(y),\nu_\Omega(y)\rangle\nu_\Omega(y)\text{ for }\sigma\text{-a.e.\ } y\in \partial\Omega,
    $$
    where $\widetilde f:\R^{n+1}\to \R$ is any Lipschitz extension of $f$ to $\R^{n+1}$.
\end{itemize}

\subsection{ADR, UR, Uniform, NTA and CAD}\label{types of domains}

We say that a Radon measure $\mu$ in $\R^{n+1}$ is ($n$-dimensional) Ahlfors-David regular if there exists $C\geq 1$ (called the ADR constant) such that
\begin{equation}\label{def:adr measure}
C^{-1}r^n\leq \mu(B(x,r))\leq Cr^n\text{ for all }x\in \supp\mu \text{ and }0<r<\diam(\supp\mu),
\end{equation}
where $\diam(\supp\mu)$ may be infinite. A closed set $E\subset \R^{n+1}$ is said to be Ahlfors-David regular if $\HH^n|_{E}$ is Ahlfors-David regular. A domain $\Omega\subset \R^{n+1}$ is said to be an Ahlfors-David regular domain if $\partial\Omega$ is Ahlfors-David regular.

\begin{notation}
    From now on, the term Ahlfors-David regular may be shortened to AD regular or ADR. Moreover, given an ADR domain $\Omega$ we will denote its surface measure by
    $$
    \sigma = \sigma_{\Omega}\coloneqq \HH^n|_{\partial\Omega}.
    $$
\end{notation}

\begin{definition}[UR set]\label{def:ur set}
    A set $E\subset \R^{n+1}$ is called ($n$-dimensional) uniformly rectifiable, UR for short, if it is ADR and there exist $\varepsilon, M\in (0,\infty)$ (called the UR constants of $E$) such that for every $x\in E$ and $r\in (0,\diam\, E)$, there is a Lipschitz map $\varphi = \varphi_{x,r} : \{y\in \R^n : |y|<r\} \to \R^{n+1}$ with Lipschitz constant $\leq M$, such that
    $$
    \HH^n (E\cap B(x,r)\cap \varphi (\{y\in \R^n : |y|<r\})) \geq \varepsilon r^n.
    $$
\end{definition}

This is a quantitative version of rectifiability introduced by David and Semmes in \cite{David-Semmes-1991,David-Semmes-1993-AMS}. It is well known that any UR set is rectifiable, for a detailed proof see \cite[p.~2629]{Hofmann2010}.

\begin{definition}[Corkscrew ball conditions]\label{def:corks balls}
    We say that a domain $\Omega\subset \R^{n+1}$ satisfies
    \begin{itemize}
        \item the (interior) corkscrew condition if there is a constant $M>1$ such that for every $x\in \partial\Omega$ and $r\in (0,\diam(\partial\Omega))$ there exists a point $A_r (x)\in \Omega$ such that $B(A_r(x),M^{-1}r)\subset B(x,r)\cap \Omega$. $A_r(x)$ is called the corkscrew point of the point $x$ at radius $r$.
        \item the exterior corkscrew condition if $\R^{n+1}\setminus\overline\Omega$ satisfies the corkscrew condition.
        \item the $2$-sided corkscrew condition if it satisfies the interior and exterior corkscrew condition.
    \end{itemize}
\end{definition}

\begin{rem}[The geometric measure theoretic outward unit normal vector $\nu$]\label{rem:def of geometric unit vector}
    If $\Omega\subset\R^{n+1}$ is an ADR domain satisfying the $2$-sided corkscrew condition, then for $\sigma$-a.e.\ $x\in\partial\Omega$ there exists a unique unit vector $\nu(x)$ (called the geometric measure theoretic outward unit vector) satisfying for $\Omega^+ \coloneqq \Omega$ and $\Omega^- \coloneqq \R^{n+1}\setminus \Omega$, both
    $$
    \lim_{r\to 0} \frac{m(\Omega^\pm \cap \{y\in B(x,r) : \pm \langle \nu (x),y-x\rangle \geq 0\})}{r^{n+1}}=0.
    $$
    We remark that the vector $\nu$ exists under more general conditions, see for instance \cite[Section 2.2]{Hofmann2010} and the references therein.
\end{rem}

\begin{definition}[Harnack chain condition]
    We say that a domain $\Omega\subset \R^{n+1}$ satisfies the Harnack chain condition if there is a constant $M>1$ such that for every $\varepsilon>0$ and $x_1,x_2\in \Omega$ with $\dist(x_i,\partial\Omega)\geq \varepsilon$ ($i=1,2$) and $|x_1-x_2|\leq 2^j\varepsilon$ for some integer $j\geq 1$, there exists a chain of open balls $\{B_k\}_{1\leq k \leq N}$ inside $\Omega$ with $N\leq Mj$ satisfying that $x_1\in B_1$, $x_2\in B_N$, $B_k\cap B_{k+1}\not = \emptyset$ (for $1\leq k \leq N-1$) and $M^{-1}r(B_k) \leq \dist (B_k,\partial\Omega)\leq Mr(B_k)$ (for $1\leq k \leq N$).
\end{definition}

\begin{definition}[Uniform domain]
    A domain $\Omega\subset \R^{n+1}$ is called uniform domain if it satisfies the Harnack chain and interior corkscrew conditions.
\end{definition}

\begin{definition}[NTA domain]
    A domain $\Omega\subset \R^{n+1}$ is called nontangentially accessible (NTA for short) domain if it is a uniform domain and it satisfies the exterior corkscrew condition.
\end{definition}

\begin{definition}[CAD]\label{def:CAD}
    A domain $\Omega\subset \R^{n+1}$ is called chord-arc domain ($1$-sided CAD or CAD for short) if it is an NTA domain and $\partial\Omega$ is ADR. We say that $\Omega$ is a $2$-sided chord-arc domain ($2$-sided CAD) if $\Omega$ and $\R^{n+1}\setminus\overline{\Omega}$ are CAD.
\end{definition}

\begin{notation}
    Given a ($2$-sided) CAD $\Omega$, we will write $C=C(\text{CAD})$ if the constant $C$ depends on the CAD constants of $\Omega$.
\end{notation}

The following is from \cite{David-Jerison-1990-LipschitzApprox,Semmes-1990-AnalysisVSGeometry}, see also \cite[Corollary 3.9]{Hofmann2010}.

\begin{theorem}\label{thm:2-sided corkscrew + ADR implies UR boundary}
    If $\Omega\subset \R^{n+1}$ is a domain satisfying the $2$-sided corkscrew condition and whose boundary is ADR, then $\partial\Omega$ is UR.
\end{theorem}

\begin{rem}\label{rem:CAD-equivalence}
    Most of the results in \cite{Hofmann2010} are presented for $2$-sided local John domains with ADR boundary. However, it was shown in \cite{Tapiola-Tolsa-2024} that this apparently weaker condition is, in fact, equivalent to being $2$-sided CAD. This allows us to apply the known results in the literature for $2$-sided local John domains with ADR boundary when working with $2$-sided CAD. For instance, the Semmes decomposition in \cite[Theorem 4.16]{Hofmann2010}, restated in \cref{Semmes' decomposition} below.
\end{rem}

\subsection{Maximal operators and boundary dyadic cubes in ADR domains}

Given a domain $\Omega\subset\R^{n+1}$, for each $x\in\partial\Omega$ we define the nontangential approach cone with aperture $\alpha >0$ as
$$
\Gamma (x) 
= \Gamma_\alpha^\Omega (x)
\coloneqq \{y\in\Omega : |y-x|<(1+\alpha)\dist(y,\partial\Omega)\}.
$$
For a function $u:\Omega\to \R$ we define the nontangential maximal function
\begin{equation}\label{nontangential maximal function}
\NN u (x) = \NN_\alpha u(x)
\coloneqq \sup_{y\in \Gamma_\alpha^\Omega(x)} |u(y)|, \quad x\in\partial\Omega,
\end{equation}
and the $\delta$-nontangential maximal function $\NN^\delta u$ of $u$ by
\begin{equation}\label{delta-nontangential maximal function}
\NN^\delta u (x) = \NN^\delta_\alpha u(x)
\coloneqq \sup_{y\in \Gamma_\alpha^\Omega(x)\cap \overline{B_{2\delta} (x)}} |u(y)|, \quad x\in\partial\Omega.
\end{equation}
For a fixed $\alpha>0$, we introduce the following definitions whenever well-defined. The nontangential limit is
\begin{equation}\label{eq:nontan limit}
u|_{\partial\Omega} (x) = u|_{\partial\Omega}^{\rm nt}(x) \coloneqq \lim_{\Gamma^\Omega(x)\ni z\to x} u(z) , \quad x\in\partial\Omega.
\end{equation}
If in addition the outward unit normal $\nu_\Omega$ exists (see \cref{rem:def of geometric unit vector}) and $u\in C^1(\Omega)$, the interior normal nontangential derivative is
\begin{equation}\label{def:interior nontang derivative}
\partial_{\nu_\Omega}^{\rm int} u(x) \coloneqq \lim_{\Gamma^\Omega(x)\ni z\to x} \langle \nu_\Omega (x), \nabla u(z)\rangle, \quad x\in\partial\Omega.
\end{equation}
For shortness, we will also write $\partial_{\nu_\Omega}$ or $\partial_\nu$ instead of $\partial_{\nu_\Omega}^{\rm int}$. If $u\in C^1 (\overline{\Omega}^c)$, the exterior normal nontangential derivative is
\begin{equation}\label{def:exterior nontang derivative}
\partial_{\nu_\Omega}^{\rm ext} u(x) \coloneqq \lim_{\Gamma^{\overline{\Omega}^c}(x)\ni z\to x} \langle \nu_\Omega (x), \nabla u(z)\rangle, \quad x\in\partial\Omega.
\end{equation}

Given an $\HH^n$-measurable set $E\subset \R^{n+1}$ with $\HH^n(E)>0$ and $1<q<\infty$, we define the uncentered Hardy-Littlewood maximal operator $\M_{q,E}$, for $f\in L^{q}_{\loc}(\HH^n|_{E})$, as
\begin{equation}\label{eq:uncentered hardy-littlewood maximal operator}
\M_{q,E} f(x)\coloneqq \sup_{\substack{r>0\\ y\in E\\ B(y,r)\ni x}} \left(\avint_{B(y,r)} |f(z)|^{q}\, d\HH^n|_{E}(z)\right)^{\frac{1}{q}}, \quad x\in E,
\end{equation}
It is well known that it is bounded from $L^p (\HH^n|_{E})$ to $L^p(\HH^n|_{E})$ if $q<p\leq\infty$, with norm $C_{p,q,n}>0$. If the set $E$ is clear from the context we write $\M_{q} = \M_{q,E}$. We will use the uncentered Hardy-Littlewood maximal operator with $E=\partial\Omega$, where $\Omega$ is an ADR domain.

For ADR domains $\Omega\subset \R^{n+1}$, the $L^p(\sigma)$ norm (with $1<p<\infty$) of the nontangential maximal function $\NN$ in \rf{nontangential maximal function} does ``not'' depend on the aperture in the sense that, for every $\alpha,\beta>0$ and any $u:\Omega\to \R$ there holds
\begin{equation}\label{eq:Lp of nontangential maximal function does not depend on aperture}
\|\NN_\alpha u\|_{L^p (\sigma)}
\approx_{\alpha,\beta}  \|\NN_\beta u\|_{L^p (\sigma)},
\end{equation}
see \cite[Proposition 2.2]{Hofmann2010}. For this fact, we will omit the aperture $\alpha>0$ in $\NN_\alpha$ and $\Gamma_\alpha$ from now on. We may fix $\alpha = 1$ for instance.

The following two lemmas provide the control of interior integrals by the nontangential maximal function. The first is for solid interior integrals (see \cite[(2.3.25) in Proposition 2.12]{Hofmann2010}) and the second\footnote{The statement as written here appears in \cite[Lemma 5.1]{Mourgoglou-Tolsa-2024-solvabilityneumannproblemelliptic-arxiv-v1}.} is for interior sets with $n$-growth (see \cite[Lemma 5.1]{Mourgoglou-Tolsa-2024-solvabilityneumannproblemelliptic-arxiv}).

\begin{lemma}\label{lemma:2.3.25 in prop 2.12 in HMT}
    Let $\Omega\subset\R^{n+1}$ be an ADR domain, and fix $\alpha>0$. Then there exists $C=C(n,\alpha,\text{ADR})>0$ such that for any measurable function $u:\Omega\to \R$ there holds
    $$
    \frac{1}{\delta} \int_{U_\delta (\partial\Omega)} |u(z)|\, dm(z) \leq C \|\NN_\alpha^\delta u\|_{L^1 (\sigma)}, \text{ }0<\delta\leq\diam(\Omega).
    $$
\end{lemma}

\begin{lemma}\label{lemma:interior integral adr to boundary}
    Let $\Omega\subset \R^{n+1}$ be an ADR domain, $B_0$ a ball centered at $\partial\Omega$, and $E\subset B_0\cap\Omega$ such that
    $$
    \HH^n (B(x,r)\cap E)\leq C_0 r^n \text{ for all }x\in E\text{ and }r>0.
    $$
    Then, for any Borel function $u:\Omega\to \R$ such that $u\in L^1_{\loc}(\HH^n|_{E})$,
    $$
    \int_E |u(z)|\, d\HH^n(z) \lesssim \int_{2B_0} \NN^{2r(B_0)}_\beta u (z)\, d\sigma(z),
    $$
    assuming the aperture $\beta>0$ to be large enough (depending only on $n$). The implicit constant above depends only on $n$, $C_0$, and the ADR constants of $\partial\Omega$.
\end{lemma}

For the construction of the Lipschitz graph in \cref{sec:construction of Lip graph} we will follow \cite[Section 8]{David-Semmes-1991}. So, given an ADR domain $\Omega$ with surface measure $\sigma$ of $\partial\Omega$, we consider a dyadic lattice $\DD_\sigma$ of ``cubes'' on $\partial\Omega$ constructed in this generality by Christ \cite{Christ1990}. Here we use the formulation of David and Semmes in \cite[Chapter 3 of Part I]{David-Semmes-1993-AMS} with codimension $1$.

\begin{lemma}[Boundary dyadic cubes]\label{lemma:def of boundary dyadic cubes}
    Given an ADR domain $\Omega$ with surface measure $\sigma$, for each $j\in \Z$ there exists a family $\DD_{\sigma,j}$ of Borel subsets of $\supp\sigma = \partial\Omega$, called the dyadic cubes of the $j$-th generation, with the following properties:
    \begin{enumerate}
        \item each $\DD_{\sigma,j}$ is a partition of $\partial\Omega$, i.e., $\partial\Omega = \bigcup_{Q\in \DD_{\sigma,j}} Q$ with $Q\cap Q^\prime =\emptyset$ whenever $Q,Q^\prime\in \DD_{\sigma,j}$ with $Q\not = Q^\prime$,
        \item if $Q\in \DD_{\sigma,i}$ and $Q^\prime \in \DD_{\sigma,j}$ for some $i\geq j$, then either $Q\cap Q^\prime =\emptyset$ or $Q\subset Q^\prime$,
        \item\label{boundary dyadic cube diam control} for all $j\in \Z$ and all $Q\in \DD_{\sigma,j}$, we have that $2^{-j}\leq \diam\, Q \leq C_\DD 2^{-j}$ and $C^{-1}2^{-jn}\leq \sigma(Q) \leq C2^{-jn}$, and
        \item\label{boundary dyadic cube thin boundary condition} for all $j\in \Z$, $Q\in \DD_{\sigma,j}$ and $0<\tau <1$, we have the so-called ``thin boundary condition'':
        $$
        \sigma(\{x\in Q : \dist(x,\partial\Omega\setminus Q)\leq \tau 2^{-j}\})
        + \sigma (\{x\in \partial\Omega\setminus Q : \dist(x,Q)\leq \tau 2^{-j}\}) \leq C_\DD \tau^{1/C_\DD} 2^{-jn}.
        $$
    \end{enumerate}
    We set $\DD_\sigma = \bigcup_{j\in \Z} \DD_{\sigma,j}$. The constants $C_\DD,C\geq 1$ in \rf{boundary dyadic cube diam control} and \rf{boundary dyadic cube thin boundary condition} above do not depend on $j$, $Q$ or $\tau$.
\end{lemma}

We follow the convention that the size of the cubes decreases as the generation index grows, that is, $\diam\, Q\approx 2^{-j}$ for $Q\in\DD_{\sigma,j}$, in accordance with the standard dyadic grid of $\R^{n+1}$ in \cref{sec:notation}.

\subsection{Compact operators}

Let us briefly recall the definition of compact operators and the Fredholm alternative.

\begin{definition}\label{def:compact operator}
    Given Banach spaces $(X,\|\cdot\|_X)$ and $(Y,\|\cdot\|_Y)$, a bounded linear operator $T:X\to Y$ is called compact if for every bounded sequence $\{x_k\}_{k\geq 1}\subset X$, the sequence $\{Tx_k\}_{k\geq 1}\subset Y$ has a convergent subsequence.
\end{definition}

\begin{theorem}[Fredholm alternative]\label{fredholm alternative}
    Let $(X,\|\cdot\|_X)$ be a Banach space, $T:X\to X$ be a compact operator, and $\lambda \in \C\setminus\{0\}$. Then exactly one of the following holds:
    \begin{itemize}
        \item the equation $Tv-\lambda v = 0$ has a non-zero solution $v\in X$, or
        \item for every $u\in X$, the equation $Tv-\lambda v = u$ has a unique solution $v\in X$. In this case, the solution $v$ depends continuously on $u$.
    \end{itemize}
\end{theorem}

\begin{rem}\label{rem:coro of fredholm alternative}
    Since the composition $T\circ B$ of a compact operator $T$ with a bounded operator $B$ is again compact, the same holds when replacing $\lambda v$ by $\mathcal I v$ for an invertible bounded operator $\mathcal I:X\to X$.
\end{rem}

\subsection{Calderón-Zygmund operators and the Riesz transform}

We say that $k:\{(x,y)\in \R^{n+1}\times\R^{n+1} : x\not=y\}\to \C$ is a Calderón-Zygmund kernel if there exist constants $C\geq 1$ and $0<\tau \leq 1$ such that for all $x, x^\prime,y\in \R^{n+1}$ with $x\not = y$, $x^\prime \not = y$, there holds
$$
\begin{aligned}
|k(x,y)|&\leq C\frac{1}{|x-y|^n}, \text{ and}\\
|k(x,y)-k(x^\prime,y)|+|k(y,x)-k(y,x^\prime)|&\leq C \frac{|x-x^\prime|^\tau}{|x-y|^{n+\tau}}, \text{ if } |x-x^\prime|\leq |x-y|/2.
\end{aligned}
$$

Given a Radon measure $\mu$ and a Calderón-Zygmund kernel $k$, we define
$$
T^k \mu (x) \coloneqq \int k(x,y) \, d\mu(y), \quad x\in \R^{n+1}\setminus \supp\mu,
$$
and as this may not converge for $x\in \supp\mu$, for $\varepsilon>0$ we define the truncated operator
$$
T^k_\varepsilon \mu(x) \coloneqq \int_{|y-x|>\varepsilon} k(x,y) \, d\mu(y),\quad x\in \R^{n+1}.
$$

Given a Radon measure $\mu$ and $f\in L^1_{\loc}(\mu)$, we define
$$
\begin{aligned}
    T^k_\mu f(x)&\coloneqq T^k(f\mu)(x), \text{ for }x\in \R^{n+1}\setminus \supp\mu,\\
    T^k_{\mu,\varepsilon} f(x)&\coloneqq T^k_\varepsilon (f\mu)(x), \text{ for }\varepsilon>0\text{ and } x\in \R^{n+1},
\end{aligned}
$$
and the maximal operator
\begin{equation}\label{eq:CZ maximal operator}
T^k_{\mu,*} f (x) \coloneqq \sup_{\varepsilon>0} |T^k_{\mu,\varepsilon} f(x)|, \quad x\in \supp\mu.
\end{equation} 
We say that $T^k_\mu$ is bounded in $L^p(\mu)$, $1<p<\infty$, if the truncated operators $T^k_{\mu,\varepsilon}$ are bounded in $L^p(\mu)$ uniformly in $\varepsilon>0$. In this case, we write $T^k_\mu:L^p(\mu)\to L^p(\mu)$ is bounded. We remark that, if $\mu$ has growth of degree $n$ (i.e., $\mu$ satisfies the upper bound in \rf{def:adr measure}), then the boundedness of $T^k_\mu$ in $L^p(\mu)$ is equivalent to the boundedness of the maximal operator $T^k_{\mu,*}$ in $L^p(\mu)$, by \cite[Theorem 2.16]{Tolsa2014}\footnote{A quick inspection of its proof reveals that the same holds if $L^2(\mu)$ is replaced by $L^p(\mu)$ for any $1<p<\infty$.} and Cotlar's inequality (take $s=1$ in \cite[(2.26)]{Tolsa2014} for instance).

\begin{definition}[Riesz transform]\label{def:Riesz transform}
    The ($n$-dimensional) Riesz kernel is the Calderón-Zygmund vector-valued kernel (with $\tau=1$)
    $$
    \widetilde k_\RR (x)\coloneqq \frac{x}{|x|^{n+1}} \text{ for }x\in \R^{n+1}\setminus\{0\}.
    $$
    The ($n$-dimensional) Riesz transform $\RR$ is defined as
    $$
    \RR \coloneqq T^k, \text{ with } k(x,y)=\widetilde k_\RR (x-y) \text{ for } x\not = y.
    $$
\end{definition}

By \cite[Proposition 4 bis]{David-1988}, if $\Omega\subset\R^{n+1}$ is an ADR domain with the $2$-sided corkscrew condition (in particular $\partial\Omega$ is UR by \cref{thm:2-sided corkscrew + ADR implies UR boundary}), then
$$
\|\RR_{\sigma,*} f\|_{L^p(\sigma)} \lesssim_{p,\text{UR}} \|f\|_{L^p(\sigma)} \text{ for all }f\in L^p(\sigma),
$$
see also \cite[Proposition 3.18]{Hofmann2010} for instance.

\section{The double layer potential}\label{sec:double layer potential}

Let $\Omega\subset \R^{n+1}$ be an ADR domain with the $2$-sided corkscrew condition (in particular $\partial\Omega$ is UR by \cref{thm:2-sided corkscrew + ADR implies UR boundary}), let $\nu \coloneqq \nu_\Omega$ be the geometric measure theoretic outward unit vector of $\Omega$ (see \cref{rem:def of geometric unit vector}), and let $f\in L^1 \left(\frac{d\sigma(x)}{1+|x|^n}\right)$. The interior double layer potential operator associated with $\Omega$ is
\begin{equation}\label{double layer potential}
\DD f (x) = \DD_\Omega f (x) \coloneqq \frac{1}{w_n} \int_{\partial\Omega} \frac{\langle \nu(y), y-x \rangle}{|x-y|^{n+1}} f(y) \, d\sigma(y), \quad x\in \R^{n+1}\setminus \partial\Omega.
\end{equation}
Here $w_n$ is the surface area of the unit sphere in $\R^{n+1}$. The double layer potential satisfies $\Delta (\DD f) = 0$ in $\R^{n+1}\setminus \partial\Omega$.

\begin{rem}
    Note that $L^p(\sigma)\subset L^1\left(\frac{d\sigma(x)}{1+|x|^n}\right)$ for all $p>1$, as $\partial\Omega$ is ADR.
\end{rem}

The boundary double layer potential, that is, the principal value version of the interior double layer potential, is defined as
\begin{equation}\label{pv double layer potential}
Kf(x) = K_\Omega f(x) \coloneqq \lim_{\varepsilon \to 0^+} K_\varepsilon f(x), \quad x\in \partial\Omega,
\end{equation}
where
\begin{equation}\label{eq:truncated epsilon dlp}
K_\varepsilon f(x)\coloneqq \frac{1}{w_n} \int_{\{y\in\partial\Omega : |y-x|>\varepsilon\}} \frac{\langle \nu(y), y-x \rangle}{|x-y|^{n+1}} f(y) \, d\sigma(y),\quad x\in \partial\Omega.
\end{equation}
The principal value in \rf{pv double layer potential} exists for $\sigma$-a.e.\ $x\in\partial\Omega$, see \cite[Proposition 5.6.7]{MitreaX3-GeometricHarmonicAnalysis-Vol1}.

Also, the maximal operator of the boundary double layer potential is defined as
$$
K_* f(x)\coloneqq \sup_{\varepsilon >0} |K_\varepsilon f(x)|, \quad x\in \partial\Omega.
$$

We also define
$$
K^*_\varepsilon f(x) \coloneqq \frac{1}{w_n} \int_{\{y\in\partial\Omega : |y-x|>\varepsilon\}} \frac{\langle \nu(x), x-y \rangle}{|x-y|^{n+1}} f(y) \, d\sigma(y),\quad x\in \partial\Omega,
$$
and the maximal operator
$$
K^*_* f(x) \coloneqq \sup_{\varepsilon>0} |K^*_\varepsilon f(x)|, \quad x\in\partial\Omega.
$$
A quick computation shows that the operator $K^*$ defined as
$$
K^* f(x) \coloneqq \lim_{\varepsilon\to 0^+} K^*_\varepsilon f(x), \quad x\in\partial\Omega,
$$
is the adjoint operator of $K$.

The interior double layer potential satisfies the jump relation
\begin{equation}\label{jump formula double layer interior}
(\DD f)|_{\partial\Omega}^{\rm nt}(x) 
=\left(\frac{1}{2} Id + K\right) f(x), \text{ for }\sigma\text{-a.e.\ } x\in \partial\Omega,
\end{equation}
see \cite[(3.31)]{Marin-Martell-Mitrea-x3}. If in addition $f\in L^p(\sigma)$ with $1<p<\infty$, then the boundary and interior double layer potentials satisfy
\begin{subequations}
\begin{align}
\|K_* f\|_{L^p (\sigma)} &\lesssim_{p,\text{UR}} \|f\|_{L^p (\sigma)},\label{maximal operator of K is bounded}\\
\|\NN (\DD f)\|_{L^p (\sigma)} & \lesssim_{p,\text{UR}} \|f\|_{L^p (\sigma)},\label{lp norm of nontangential of double layer potential}
\end{align}
\end{subequations}
see \cite[(3.3.5) and (3.3.6)]{Hofmann2010} respectively. Here UR denotes that the constant depends on the UR constants of $\partial\Omega$. The second estimate also depends on the aperture of $\NN$.

\subsection{Proof of \cref{thm:double layer potential is almost compact} and consequences}

In order to study the boundary double layer potential $K$ in \cref{thm:double layer potential is almost compact}, for $0<t<T<\infty$ we define its truncations by
$$
\begin{aligned}
	K_{s} f(x)&\coloneqq\lim_{\varepsilon\to 0^+}\frac{1}{w_n} \int_{\{y\in\partial\Omega : \varepsilon<|y-x|\leq t\}} \frac{\langle \nu(y), y-x \rangle}{|x-y|^{n+1}} f(y) \, d\sigma(y) , & x\in\partial\Omega,\\
	K_{i} f(x)&\coloneqq\frac{1}{w_n} \int_{\{y\in\partial\Omega : t< |y-x|\leq T\}} \frac{\langle \nu(y), y-x \rangle}{|x-y|^{n+1}} f(y) \, d\sigma(y) , & x\in\partial\Omega,\\
	K_{l} f(x)&\coloneqq\frac{1}{w_n} \int_{\{y\in\partial\Omega : |y-x|>T\}} \frac{\langle \nu(y), y-x \rangle}{|x-y|^{n+1}} f(y) \, d\sigma(y) , & x\in\partial\Omega,
\end{aligned}
$$
where $s$, $i$ and $l$ stand for small, intermediate and large scales respectively. If we want to stress the scales we will write $K_{s(t)}$, $K_{i(t,T)}$ and $K_{l(T)}$ respectively. We define $K^*_{s(t)}$, $K^*_{i(t,T)}$ and $K^*_{l(T)}$ analogously from the definition of $K^*$.

For $\widetilde R>0$ we also break the intermediate scales operator as $K_{i} = \characteristic_{B_{\widetilde R} (0)} K_{i} + \characteristic_{B_{\widetilde R} (0)^c} K_{i}$. So, for $f\in L^p (\sigma)$ and $x\in \partial\Omega$ we will decompose the double layer potential as
\begin{equation}\label{full truncation of double layer potential}
K f(x) = K_s f(x) + K_l f(x) + \characteristic_{B_{\widetilde R} (0)} (x) K_{i} f(x) + \characteristic_{B_{\widetilde R} (0)^c}(x) K_{i} f(x).
\end{equation}

Next, we present a series of results to summarize the relevant properties of the operators in the decomposition \rf{full truncation of double layer potential}.

\begin{lemma}\label{trucated and intermediate scales is compact}
    Let $\Omega\subset \R^{n+1}$ be an ADR domain with the $2$-sided corkscrew condition. For any $0<t<T<\infty$ and $\widetilde R>0$, the operator $\characteristic_{B_{\widetilde R}(0)} K_{i(t,T)} : L^p (\sigma)\to L^p (\sigma)$ is compact for all $p\in (1,\infty)$.
\end{lemma}

We write the proof in \cref{sec:proof of compact}. The compact operator in \cref{thm:double layer potential is almost compact} will be in fact $T=\characteristic_{B_{\widetilde R}(0)} K_{i(t,T)}$, for some choice of parameters.

The following two results control the $L^p(\sigma)$ norm of the boundary double layer potential at small scales and, respectively, at scales far from the ``bad'' balls in the definition of $\delta$-$(s,S;R)$ domains (see \cref{def:final assumptions of domain}).

\begin{lemma}\label{lemma:double layer potentials in small scales has small norm}
    Let $\Omega\subset \R^{n+1}$ be a $2$-sided CAD. Assume also that there is $\delta>0$ and $s>0$ such that
    $$
    \sup_{\substack{x\in\partial\Omega\\0<r\leq s}}\avint_{B(x,r)} |\nu-m_{B(x,r)}\nu|\, d\sigma \leq\delta.
    $$
    Given $p\in (1,\infty)$, there exists $t_0 = t_0 (s,\delta,p,\text{CAD},n)>0$ such that for every $0<t\leq t_0$ there holds
    $$
    \|K_{s(t)}\|_{L^p (\sigma)\to L^p (\sigma)} \lesssim \delta^{1/4},
    $$
    where the involved constant depends on $n$, the CAD constants of $\Omega$ and $p$.
\end{lemma}

\begin{lemma}\label{lemma:double layer potentials far from bad ball has small norm}
    Let $\Omega\subset \R^{n+1}$ be a $\delta$-$(s,S;R)$ domain. Given $p\in (1,\infty)$ and $t>0$, there exists $\widetilde R=\widetilde R(R,\delta,t,p,\text{CAD},n)$ such that there holds 
    $$
    \|\characteristic_{B_{\widetilde R}(0)^c}K_{s(t)}\|_{L^p (\sigma)\to L^p (\sigma)} \lesssim \delta^{1/4},
    $$
    where the involved constant depends on $n$, the CAD constants of $\Omega$ and $p$.
\end{lemma}

We prove the preceding two lemmas in \cref{sec:double layer potential in small scales}. Note that, to obtain the claimed bound, the first lemma gives a sufficiently small parameter $t>0$, while the second one provides the truncation parameter $\widetilde R$ given a scale parameter $t>0$, which is not assumed to be small in this case.

The next result is the main work in this article, and provides the small norm of the large scales double layer potential maximal operator, defined for $f\in L^p (\sigma)$ as
$$
K_{l,*} f(x) \coloneqq \sup_{\varepsilon\geq T} |K_\varepsilon f (x)|,\quad x\in \partial\Omega,
$$
assuming flatness conditions at large scales. The proof is deferred to \cref{sec:proof of large scale bounded by bmo}.

\begin{theorem}\label{large scale bounded by bmo norm}
    Let $\Omega\subset \R^{n+1}$ be a $2$-sided CAD with unbounded boundary. Assume also that there are $\delta_\beta,\delta_*>0$ and $S>0$ such that the following two conditions hold:
    \begin{itemize}
        \item Small Jones' $\beta_{\infty,\partial\Omega}$ coefficient at large scales: for $x\in \partial\Omega$ and $r\geq S$,
        \begin{equation}\label{betas in large scales}
            \beta_{\infty,\partial\Omega} (B(x,r))  \leq \delta_\beta.
        \end{equation}
        \item Small $\bmo$ norm of $\nu$ at large scales: for $x\in \partial\Omega$ and $r\geq S$,
        \begin{equation}\label{bmo norm in large scales}
            \avint_{{B(x,r)}} |\nu(z)-m_{B(x,r)}\nu| \, d\sigma(z) \leq \delta_\ast.
        \end{equation}
    \end{itemize}
    Denote $\delta=\max\{\delta_\beta,\delta_*\}$. Given $1<p<\infty$, there exists $\theta = \theta (n,p)>0$ (see \rf{choice of A}) and $T = T (p,\delta,S,\text{CAD},n) \gg 100S$ such that
    \begin{equation}\label{eq:Lp norm of large scales dlp}
    \|K_{l(T),*}\|_{L^p (\sigma)\to L^p (\sigma)} \lesssim \delta^\theta,
    \end{equation}
    where the involved constant depends on $n$, the CAD constants of $\Omega$ and $p$.
\end{theorem}

Although the truncated operators appearing in \rf{full truncation of double layer potential} and their properties in \cref{trucated and intermediate scales is compact,lemma:double layer potentials in small scales has small norm,lemma:double layer potentials far from bad ball has small norm,large scale bounded by bmo norm} have not yet been studied and proved, we now turn to the proof of \cref{thm:double layer potential is almost compact}.

\begin{proof}[Proof of \cref{thm:double layer potential is almost compact}]
    Given $\varepsilon>0$, by \cref{lemma:double layer potentials in small scales has small norm,lemma:double layer potentials far from bad ball has small norm,large scale bounded by bmo norm} there exists $\delta_0 = \delta_0(\varepsilon,p,\text{CAD},n)$ and $0<t\ll 1\ll T \ll \widetilde R$ such that if $\delta\leq \delta_0$, then for all $f\in L^p(\sigma)$ there holds
    $$
    \|K f - \characteristic_{B_{\widetilde R} (0)} K_{i(t,T)} f\|_{L^p(\sigma)}\overset{\text{\rf{full truncation of double layer potential}}}{=} \|K_{s(t)} f + K_{l(T)} f + \characteristic_{B_{\widetilde R} (0)^c} K_{i(t,T)}f\|_{L^p(\sigma)}<\varepsilon,
    $$
    where we used that $|K_{i(t,T)} f| \leq |K_{s(t)}f|+|K_{s(T)}f|$ $\sigma$-a.e.\ on $\partial\Omega$. By \cref{trucated and intermediate scales is compact}, the operator $T=\characteristic_{B_{\widetilde R} (0)} K_{i(t,T)}$ is compact (with abuse of notation using $T$ for both the compact operator and the scale). Finally, since $T$ is compact, its adjoint $T^*$ is also compact by Schauder's theorem and moreover
    $$
    \|K^* - T^*\|_{L^{p^\prime}(\sigma)\to L^{p^\prime} (\sigma)} = \|K-T\|_{L^p(\sigma)\to L^p (\sigma)} < \varepsilon,
    $$
    as claimed.
\end{proof}

As a consequence of \cref{thm:double layer potential is almost compact}, in the following result we obtain that injectivity implies invertibility for $\frac{1}{2}Id + K$ and $-\frac{1}{2}Id + K^*$, under the assumption of enough flatness.

\begin{coro}\label{reduction:all equivalent}
    Let $\Omega\subset \R^{n+1}$ be a $\delta$-$(s,S;R)$ domain (see \cref{def:final assumptions of domain}), $\lambda_0>0$, $1<p<\infty$ and $p^\prime = p/(p-1)$ its Hölder conjugate exponent. There exists $\delta_0=\delta_0 (\lambda_0,p,\text{CAD},n)$ such that if $\delta \leq \delta_0$ and $\lambda\in\C$ with $|\lambda|\geq \lambda_0$, then the following are equivalent:
    \begin{enumerate}
        \item\label{item:id + dlp invertible} $\lambda Id +K$ is invertible in $L^p(\sigma)$,
        \item\label{item:id + dlp injective} $\lambda Id +K$ is injective in $L^p(\sigma)$,
        \item\label{item:id + dlp surjective} $\lambda Id +K$ is surjective in $L^p(\sigma)$,
        \item\label{item:id + adjoint dlp invertible} $\lambda Id + K^*$ is invertible in $L^{p^\prime}(\sigma)$,
        \item\label{item:id + adjoint dlp injective} $\lambda Id + K^*$ is injective in $L^{p^\prime}(\sigma)$,
        \item\label{item:id + adjoint dlp surjective} $\lambda Id + K^*$ is surjective in $L^{p^\prime}(\sigma)$.
    \end{enumerate}
	\begin{proof}
            By \cref{thm:double layer potential is almost compact} there is $\delta_0 = \delta_0 (\lambda_0,p,\text{CAD},n)$ such that if $\delta\leq \delta_0$ then there is a compact operator $T$ such that $\|K-T\|_{L^p(\sigma)\to L^p (\sigma)}=\|K^*-T^*\|_{L^{p^\prime}(\sigma)\to L^{p^\prime} (\sigma)}<\lambda_0\leq |\lambda|$.
        
            By Neumann series, the operator $\mathcal I \coloneqq \lambda Id + K - T$ is invertible. By the Fredholm alternative \cref{fredholm alternative} and \cref{rem:coro of fredholm alternative}, if either $\lambda Id+K$ is injective or surjective in $L^p(\sigma)$ we get that it is bijective, and by the bounded inverse theorem we conclude that $\left(\lambda Id + K\right)^{-1}:L^p(\sigma)\to L^p(\sigma)$ is a linear bounded operator. This concludes the equivalence between items \rf{item:id + dlp invertible}, \rf{item:id + dlp injective} and \rf{item:id + dlp surjective}.
  
            The same argument holds mutatis mutandis with $\lambda Id + K^*$ in $L^{p^\prime}(\sigma)$, whence we get the equivalences between \rf{item:id + adjoint dlp invertible}, \rf{item:id + adjoint dlp injective} and \rf{item:id + adjoint dlp surjective}.
    
            Now, since an (arbitrary) operator $U$ is injective if its adjoint $U^*$ is surjective, in particular \rf{item:id + adjoint dlp surjective} implies \rf{item:id + dlp injective}, and \rf{item:id + dlp surjective} implies \rf{item:id + adjoint dlp injective}.
	\end{proof}
\end{coro}

\subsection{The compact operator: Proof of \cref{trucated and intermediate scales is compact}}\label{sec:proof of compact}

We conclude this section by seeing that $\characteristic_{B_{\widetilde R} (0)} K_{i}$ is compact.

\begin{proof}[Proof of \cref{trucated and intermediate scales is compact}]
    First note that the kernels of $\characteristic_{B_{\widetilde R} (0)} K_{i(t,T)}$ and $(\characteristic_{B_{\widetilde R} (0)} K_{i(t,T)})^*$ are
    $$
    \begin{aligned}
    &\characteristic_{B_{\widetilde R}(0)} (x) \characteristic_{\{y\in \partial\Omega : t< |y-x|\leq T\}} (y)\frac{\langle \nu(y), y-x\rangle}{|x-y|^{n+1}}, \text{ and}\\
    &\characteristic_{\{y\in B_{\widetilde R}(0)\cap\partial\Omega : t< |y-x|\leq T\}} (y) \frac{\langle \nu(x), x-y\rangle}{|x-y|^{n+1}},
    \end{aligned}
    $$
    respectively. In particular, both satisfy the so-called Hilbert-Schmidt condition, see the line before \cite[Theorem 0.45]{Folland1995} for instance. Hence, by \cite[Theorem 0.45]{Folland1995} we have that both $\characteristic_{B_{\widetilde R} (0)} K_{i(t,T)}$ and $(\characteristic_{B_{\widetilde R} (0)} K_{i(t,T)})^*$ are bounded and compact from $L^2 (\sigma)$ to $L^2 (\sigma)$. This concludes the proof for the case $p=2$.

    For $p\in (1,2)$, fix any $p_0\in (1,p)$. Since the boundary is UR, we have that both $\characteristic_{B_{\widetilde R} (0)} K_{i(t,T)}$ and $(\characteristic_{B_{\widetilde R} (0)} K_{i(t,T)})^*$ are bounded from $L^{p_0}(\sigma)$ to $L^{p_0}(\sigma)$. By the interpolation theorem in \cite{krasnoselskii-1960} (see also \cite[Theorem 3.10]{krasnoselskii-book-1976}) between compact operators and bounded operators, we conclude that both $\characteristic_{B_{\widetilde R} (0)} K_{i(t,T)}$ and $(\characteristic_{B_{\widetilde R} (0)} K_{i(t,T)})^*$ are compact from $L^p (\sigma)$ to $L^p(\sigma)$. This gives the case $p\in (1,2)$.

    The case $p\in (2,\infty)$ follows from the result obtained in the previous paragraph and Schauder's theorem, which states that a bounded linear operator between Banach spaces is compact if and only if its adjoint is compact. This concludes the proof of the lemma.
\end{proof}

\section{The double layer potential at small scales. Proofs of \cref{lemma:double layer potentials in small scales has small norm,lemma:double layer potentials far from bad ball has small norm}}\label{sec:double layer potential in small scales}

In this section we study the domain and the double layer potential for a fixed scale. We conclude the section by proving \cref{lemma:double layer potentials in small scales has small norm,lemma:double layer potentials far from bad ball has small norm}.

\begin{theorem}[{\cite[Theorem 1.3]{Tapiola-Tolsa-2024}}]\label{Poincare inequality}
	Let $\Omega\subset \R^{n+1}$ be a $2$-sided CAD. Then the following weak $1$-Poincaré inequality for Lipschitz functions on $\partial\Omega$ holds: there exist constants\footnote{The notation $C_P$ is to specify the constant appearing in the Poincaré inequality.} $C_P\geq 1$ and $\Lambda\geq 1$ such that for every Lipschitz function $f$ on $\partial\Omega$, every $x\in\partial\Omega$ and every $r>0$, for $\Delta \coloneqq \Delta(x,r)$ we have
	$$
	\avint_{\Delta} |f(z)-m_\Delta f| \, d\sigma(z) \leq C_P r \avint_{\Lambda \Delta} |\nabla_t f(z)|\, d\sigma(z),
	$$
	where $\nabla_t$ is the tangential gradient of $f$.
\end{theorem}

This is a no-tail version of \cite[Proposition 4.13]{Hofmann2010} for Lipschitz functions, which is enough for our applications.

By the same proof in \cite[Theorem 4.14]{Hofmann2010}, but using the refined Poincaré inequality above, we obtain the following estimate for the theoretical unit normal vector $\nu$. We provide the proof for completeness.

\begin{theorem}\label{pointwise estimate vector dot product}
Let $\Omega\subset \R^{n+1}$ be a $2$-sided CAD. Then there exist $C_P^\prime = C (C_P) > 0$ such that for every $\alpha \in (0,1)$, every $x\in\partial\Omega$ and every $r>0$, for $\Delta\coloneqq\Delta(x,r)$ there holds
$$
\sup_{y\in 2\Delta} r^{-1} |\langle x-y, m_\Delta \nu \rangle| \leq C_P^\prime \left(\avint_{\Lambda\Delta} |\nu(z)-m_\Delta\nu|^{\frac{n}{1-\alpha}}\, d\sigma(z)\right)^{\frac{1-\alpha}{n}},
$$
with $\Lambda\geq 1$ as in the weak $1$-Poincaré inequality in \cref{Poincare inequality}.
\begin{proof}
    Define the Lipschitz function $g_x (z) \coloneqq \langle x-z, m_\Delta \nu \rangle$ for $z\in\R^{n+1}$. As in \cite[(4.2.25)]{Hofmann2010}, the claim follows from the particular case $y^\prime = x$ of
    \begin{equation}\label{holder reg of normal vector}
        |g_x (y) - g_x (y^\prime)| \leq C_P^\prime r^{1-\alpha} |y-y^\prime|^\alpha \left(\avint_{\Lambda\Delta} |\nu(z)-m_{\Delta}\nu|^{\frac{n}{1-\alpha}}\, d\sigma(z)\right)^{\frac{1-\alpha}{n}},
    \end{equation}
    for all $y,y^\prime\in 2\Delta$ and each $\alpha \in (0,1)$.
	
    Let us see this. Fixed $\alpha \in (0,1)$, let $p=n/(1-\alpha)>1$, equivalently $\alpha=1-n/p$. For $\sigma$-a.e.\ $z\in\partial\Omega$ we have
    $$
    |\nabla_t g_x (z)| = |m_\Delta \nu - \langle m_\Delta \nu,\nu(z)\rangle\nu(z)|
    = |m_\Delta \nu - \nu(z) - \langle m_\Delta \nu - \nu(z),\nu(z)\rangle\nu(z)|
    \leq 2|\nu(z)-m_\Delta \nu|.
    $$
    Now, for any arbitrary boundary ball $\Delta_s \subset \Delta$ (centered at $\partial\Omega$) of radius $s$, by the $1$-Poincaré inequality in \cref{Poincare inequality} we have
    \begin{multline*}
        \frac{1}{s}\avint_{\Delta_s}|g_x (z)-m_{\Delta_s} g_x| \, d\sigma(z)
        \leq C_P \avint_{\Lambda\Delta_s} |\nabla_t g_x (z)| \, d\sigma(z)
        \leq 2C_P \avint_{\Lambda\Delta_s} |\nu(z)-m_\Delta \nu| \, d\sigma(z)\\
        \leq 2C_P \left(\avint_{\Lambda\Delta_s} |\nu(z)-m_\Delta \nu|^p \, d\sigma(z)\right)^{1/p}
        \lesssim C_P \frac{r^{\frac{n}{p}}}{s^{\frac{n}{p}}} \left(\avint_{\Lambda\Delta} |\nu(z)-m_\Delta \nu|^p \, d\sigma(z) \right)^{1/p}.
    \end{multline*}
    By the choice of $p$ in terms of $\alpha$, for all $\Delta_s \subset \Delta$ we get
    $$
    \frac{1}{s^\alpha}\avint_{\Delta_s} |g_x (z)-m_{\Delta_s} g_x| \, d\sigma(z)
    \leq C_P^\prime r^{1-\alpha} \left(\avint_{\Lambda\Delta} |\nu(z)-m_{\Delta} \nu|^{\frac{n}{1-\alpha}}\, d\sigma(z)\right)^{\frac{1-\alpha}{n}}.
    $$
    This implies the Hölder regularity in \rf{holder reg of normal vector} by Meyer's criterion in \cite{Meyer-1964}.
\end{proof}
\end{theorem}

From \cref{pointwise estimate vector dot product} we deduce the following estimate at large scales.

\begin{coro}
Under the assumptions of \cref{pointwise estimate vector dot product} and \rf{bmo norm in large scales}, for every $x\in\partial\Omega$ and every $r\geq S/\Lambda$,
\begin{equation}\label{coro:pointwise estimate vector dot product}
\sup_{y\in \Delta (x,2r)} r^{-1} |\langle x-y, m_{\Delta(x,r)}\nu\rangle|
\lesssim_{\Lambda, C_P} \delta_*^{\frac{1-\alpha}{n}}.
\end{equation}
\begin{proof}
By \cref{pointwise estimate vector dot product}, it suffices to bound its right-hand side, which is controlled by
$$
\begin{aligned}
\left(\avint_{\Lambda\Delta} |\nu(z)-m_\Delta \nu|^{\frac{n}{1-\alpha}}\, d\sigma(z) \right)^{\frac{1-\alpha}{n}}
\lesssim&_\Lambda
\left(\avint_{\Lambda\Delta} |\nu(z)-m_{\Lambda\Delta}\nu|^{\frac{n}{1-\alpha}}\, d\sigma(z)\right)^{\frac{1-\alpha}{n}}\\
&+ \avint_{\Lambda\Delta} |\nu(z)-m_{\Lambda\Delta}\nu|\, d\sigma(z).
\end{aligned}
$$
Since $|\nu(z)-m_{\Lambda\Delta}\nu|\leq 2$, we uniformly have
\begin{equation}\label{eq:reduce high power by bounded oscillation}
\left(\avint_{\Lambda\Delta} |\nu(z)-m_{\Lambda\Delta}\nu|^{\frac{n}{1-\alpha}}\, d\sigma(z)\right)^{\frac{1-\alpha}{n}}
\lesssim \left(\avint_{\Lambda\Delta} |\nu(z)-m_{\Lambda\Delta}\nu|\, d\sigma(z)\right)^{\frac{1-\alpha}{n}}.
\end{equation}
Since $\Lambda r\geq S$, by \rf{bmo norm in large scales} both terms in the display above are $\lesssim \delta_*^{\frac{1-\alpha}{n}}$, which proves the claim.
\end{proof}
\end{coro}

For completeness, we state the Semmes decomposition, as in \cite[Theorem 4.16]{Hofmann2010} (see \cref{rem:CAD-equivalence} for the applicability to $2$-sided CAD).

\begin{theorem}\label{Semmes' decomposition}
    Let $\Omega\subset \R^{n+1}$ be a $2$-sided CAD. Then there exist $C_* \geq 1$ and $C_1, C_2, C_3, C_4>0$ (depending on the CAD constants of $\Omega$ and $n$) with the property that if for every compact set $\mathcal{K} \subset \R^{n+1}$, there exists $R_{\mathcal{K}}>0$ for which
    $$
    \sup _{x \in \mathcal{K} \cap \partial\Omega}\|\nu\|_*\left(\Delta\left(x, R_{\mathcal{K}}\right)\right) \leq \delta\leq 1/C_*,
    $$
    then for every compact set $\mathcal{K} \subset \R^{n+1}$ and $T\geq 1$, setting
    $$
    \widetilde{\mathcal{K}}_T \coloneqq \{y\in \R^{n+1} : \dist(y,\mathcal K) \leq T\},
    $$
    and
    \begin{equation}\label{choice of R_* in semmes' decomposition}
    R_{*,T,\mathcal K} = 
    \min\{
    \delta R_{\widetilde{\mathcal{K}}_T}/C_*,
    \diam(\partial\Omega)/C_*,T\},
    \end{equation}
    for $x \in \mathcal{K} \cap \partial\Omega$ and $0<r \leq R_{*,T,\mathcal K}$ the following holds:
    \begin{enumerate}
    \item There exists a unit vector $\vec{n}_{x, r}$ and a Lipschitz function
	$$
	h: H(x, r)\coloneqq\langle\vec{n}_{x, r}\rangle^{\perp} \to \R \text{ with }\|\nabla h\|_{L^{\infty}} \leq C_3 \sqrt{\delta},
	$$
	and whose graph
	$$
	\mathcal{G}\coloneqq\left\{y=x+\zeta+t \vec{n}_{x, r}: \zeta \in H(x, r), t=h(\zeta)\right\}
	$$
	(in the coordinate system $y=(\zeta, t) \Longleftrightarrow y=x+\zeta+t \vec{n}_{x, r}$, $\zeta \in H(x, r)$, $t \in \R$) is a good approximation of $\partial\Omega$ in the cylinder
	$$
	\mathcal{C}(x, r)\coloneqq\left\{x+\zeta+t \vec{n}_{x, r}: \zeta \in H(x, r),|\zeta| \leq r,|t| \leq r\right\}
	$$
	in the sense
	$$
	\sigma\left(\mathcal{C}(x, r) \cap \left((\partial\Omega\setminus \mathcal G)\cup (\mathcal G \setminus\partial\Omega)\right)\right) \leq C_1 w_n r^n \exp \left(-C_2 / \sqrt{\delta}\right) .
	$$
	\item There exist two disjoint sets $G(x,r)$ (``good'') and $E(x, r)$ (``evil'') such that
	$$
	\begin{aligned}
		& \mathcal{C}(x, r) \cap \partial\Omega=G(x, r) \cup E(x, r) \text { with } G(x,r) \subset \mathcal{G},\\
		& \sigma(E(x, r)) \leq C_1 w_n r^n \exp \left(-C_2 / \sqrt{\delta}\right),
	\end{aligned}
	$$
	and moreover, if $\Pi: \R^{n+1} \to H(x, r)$ is defined by $\Pi(y)=\zeta$ if $y=x+\zeta+t \vec{n}_{x, r} \in$ $\R^{n+1}$ with $\zeta \in H(x, r)$ and $t \in \R$, then
	$$
	\left|y-\left(x+\Pi(y)+h(\Pi(y)) \vec{n}_{x, r}\right)\right| \leq C_3 \sqrt{\delta} \dist(\Pi(y), \Pi(G(x, r))) \text{ for all } y \in E(x, r),
	$$
	and
	$$
	\begin{aligned}
		& \mathcal{C}(x, r) \cap \partial\Omega \subset\left\{x+\zeta+t \vec{n}_{x, r}:|t| \leq C_3 \sqrt{\delta} r, \zeta \in H(x, r)\right\} \\
		& \Pi(\mathcal{C}(x, r) \cap \partial\Omega)=\{\zeta \in H(x, r):|\zeta|<r\}.
	\end{aligned}
	$$
	\item $\left(1-C_4 \sqrt{\delta}\right) w_n r^n \leq \sigma(\Delta(x, r)) \leq\left(1+C_4 \sqrt{\delta}\right) w_n r^n$.
    \end{enumerate}
\end{theorem}

A few comments are in order. The case $T\geq 1$ follows from the case $T=1$ by scaling. For the case $T=1$, a careful inspection of \cite[Proof of Theorem 4.16]{Hofmann2010} reveals the following:
\begin{enumerate}
    \item The constants $C_*\geq 1$ and $C_1,C_2,C_3,C_4>0$ in \cite[Theorem 4.16]{Hofmann2010} do not depend on the compact set $\mathcal K$.
    \item In \cite[p.~2703, l.~6]{Hofmann2010} the authors define
    $$
    R_* = \min\{\delta R_{\widetilde{\mathcal K}_1}/(8C),R_{\widetilde{\mathcal K}_1}/8,R_0/100,1\},
    $$
    where $R_0$ is the constant used in the statement of \cite[Theorem 4.14]{Hofmann2010} and $C>0$ is the geometrical constant appearing in \cite[(4.2.20)]{Hofmann2010}, i.e., $C=C(\text{CAD})$. It turns out that $R_0\approx_{\text{CAD}} \diam(\partial\Omega)$ since $\Omega$ is a $2$-sided CAD, by \cref{def:CAD} and the fact that being a $2$-sided local John domain (see \cite[Definition 3.12]{Hofmann2010} with $R=\diam(\partial\Omega)$) with ADR boundary is equivalent to being $2$-sided CAD, see \cite[Theorem 1.2 and Corollary 1.5]{Tapiola-Tolsa-2024}.
\end{enumerate}
All in all, reusing the same notation for the constants, there exists a constant $C_*\geq 1$ such that the conclusions of the theorem hold with the choice of $R_{*,T,\mathcal K}$ in \rf{choice of R_* in semmes' decomposition}.

Given a $2$-sided CAD $\Omega\subset \R^{n+1}$, for a boundary ball $\Delta_0 \coloneqq \Delta (\xi,r_0)$, i.e., $\xi\in \partial\Omega$ and $r_0>0$, let $k_0\in \Z$ the maximal index satisfying $r_0 \leq 2^{-k_0}$ and we define
\begin{equation}\label{def:neighborhood}
I_0 (\Delta_0) \coloneqq \bigcup_{Q\in\QQ_0(\Delta_0)} Q,\text{ where }\QQ_0(\Delta_0)\coloneqq \{Q\in\DD_{\sigma,k_0} : Q\cap 2\Delta_0\not=\emptyset\},
\end{equation}
recall $\DD_{\sigma,k_0}$ is the family of dyadic cubes of $\partial\Omega$ in \cref{lemma:def of boundary dyadic cubes}.

Here we state the localized $L^p(\sigma)$ norm of the double layer potential. This estimate is not stated as a separate result in \cite{Hofmann2010}, but it is obtained by specializing the argument in the proof of \cite[Theorem 4.36]{Hofmann2010} to our localized setting. More precisely, it is the localized bound \cite[(4.5.59)]{Hofmann2010} obtained in the course of that proof, which follows from the good lambda inequality \cite[(4.5.17)]{Hofmann2010}; the latter relies in turn on Semmes' decomposition \cite[Theorem 4.16]{Hofmann2010} (see \cref{Semmes' decomposition}) and on the $L^p$ bounds for the double layer potential on Lipschitz graphs \cite[Theorem 4.34]{Hofmann2010}.

\begin{theorem}\label{theorem:small scales Lp}
    Let $\Omega\subset \R^{n+1}$ be a $2$-sided CAD, $p\in (1,\infty)$, $\xi_0\in\partial\Omega$ and $r_0>0$. There exists $C_0=C_0(\text{ADR},p,n)\geq 1$ such that if the conclusions in the Semmes decomposition \cref{Semmes' decomposition} are valid with $\delta>0$ for all $x\in I_0(\Delta(\xi_0,r_0))$ and all $0<r\leq C_0 r_0$, and
    $$
    \sup_{x\in I_0(\Delta(\xi_0,r_0))} \|\nu\|_* (\Delta(x,C_0 r_0)) \leq\delta,
    $$
    then
    $$
    \int_{I_0(\Delta(\xi_0,r_0))} |K_* f|^p \, d\sigma \lesssim \delta^{p/4} \int_{\Delta(\xi_0,r_0)} |f|^p \, d\sigma \text{ for any }f\in L^p (\Delta(\xi_0,r_0)),
    $$
    where the involved constant depends on the CAD constants of $\Omega$, $p$ and $n$.
\end{theorem}

\begin{rem}
    When invoking \cref{theorem:small scales Lp}, we may suppose that $C_0 = C_0(\text{ADR},p,n)$ is sufficiently large so that $I_0(\Delta(x,r))\subset \Delta(x,C_0r)$ for all $x\in\partial\Omega$ and all $r>0$.
\end{rem}

Using \cref{theorem:small scales Lp}, we are now ready to prove \cref{lemma:double layer potentials in small scales has small norm,lemma:double layer potentials far from bad ball has small norm}.

\begin{proof}[Proof of \cref{lemma:double layer potentials in small scales has small norm}]
    By the John-Nirenberg inequality, let $C_1\geq 1$ be a constant satisfying
    $$
    \sup_{x\in\partial\Omega}\|\nu\|_{*}(\Delta(x,s/2)) \leq C_1\delta.
    $$
    Let $C_*$ and $C_0$ be the constants in \cref{Semmes' decomposition,theorem:small scales Lp} respectively. We can assume $\delta\leq 1/(C_1 C_*)$, since otherwise, for any $t>0$, we have
    $$
    \|K_{s(t)}\|_{L^p (\sigma)\to L^p (\sigma)} \leq 2\|K_*\|_{L^p(\sigma)\to L^p (\sigma)} \overset{\text{\rf{maximal operator of K is bounded}}}{\lesssim} 1
    \lesssim \delta^{1/4}.
    $$
    
    We fix
    $$
    t_0=\frac{1}{2C_0}\min\{s/2,\delta s/(2C_*),\diam(\partial\Omega)/C_*,1\},
    $$ 
    (recall $R_{*,1,\cdot}$ from \rf{choice of R_* in semmes' decomposition}) so that for all $0<t\leq t_0$ we have
    $$
    \sup_{x\in\partial\Omega} \|\nu\|_* (B(x,2C_0 t)) \leq C_1 \delta,
    $$
    and the conclusions in the Semmes decomposition \cref{Semmes' decomposition} hold (with $C_1\delta$) for all $x\in\partial\Omega$ and all $0<r\leq 2 C_0 t$. For a fixed $0<t\leq t_0$ and any $x\in\partial\Omega$, applying \cref{theorem:small scales Lp} we obtain
    \begin{equation}\label{eq:localization at small scale}
    \|\characteristic_{I_0 (\Delta(x,2t))} K_{s(t)} (f\characteristic_{\Delta(x,2t)}) \|_{L^p(\sigma)}
    \leq 2\|\characteristic_{I_0 (\Delta(x,2t))} K_* (f\characteristic_{\Delta(x,2t)}) \|_{L^p(\sigma)}
    \lesssim \delta^{1/4} \|f\characteristic_{\Delta(x,2t)}\|_{L^p(\sigma)},
    \end{equation}
    with $I_0(\cdot)$ as in \rf{def:neighborhood}.

    By the $5R$-covering theorem, let $\{\Delta_i\}_{i\in \N}$ be a subfamily of $\{\Delta(x,t)\}_{x\in\partial\Omega}$ such that $\partial\Omega\subset \bigcup_{i\in\N} \Delta_i$ and $\{\Delta_i/5\}_{i\in \N}$ is pairwise disjoint. Since $\{\Delta_i/5\}_{i\in \N}$ is pairwise disjoint and all balls have the same radius, we have that the family $\{2\Delta_i\}_{i\in \N}$ has finite overlap, with constant depending only on the dimension. For any $f\in L^p (\sigma)$ we have
    $$\|K_{s(t)} f\|_{L^p(\sigma)}^p
    \leq \sum_{i\in \N} \|\characteristic_{\Delta_i} K_{s(t)} f \|_{L^p(\sigma)}^p
    = \sum_{i\in \N} \|\characteristic_{\Delta_i} K_{s(t)} (f\characteristic_{2\Delta_i}) \|_{L^p(\sigma)}^p.
    $$
    For each $i\in \N$, let $I_i \coloneqq I_0 (2\Delta_i) \supset 4\Delta_i$. Applying \rf{eq:localization at small scale}, we obtain
    $$
    \|\characteristic_{\Delta_i} K_{s(t)} (f\characteristic_{2\Delta_i}) \|_{L^p(\sigma)}^p
    \leq \|\characteristic_{I_i} K_{s(t)} (f\characteristic_{2\Delta_i}) \|_{L^p(\sigma)}^p
    \lesssim \delta^{p/4} \|f\characteristic_{2\Delta_i}\|_{L^p(\sigma)}^p.
    $$
    By the finite overlap of the family $\{2\Delta_i\}_{i\in \N}$, we conclude
    $$\|K_{s(t)} f\|_{L^p(\sigma)}
    \lesssim \delta^{1/4} \|f\|_{L^p(\sigma)}
    $$
    as claimed.
\end{proof}

The proof of \cref{lemma:double layer potentials far from bad ball has small norm} follows a similar approach to the proof of \cref{lemma:double layer potentials in small scales has small norm}. However, in the previous proof, we chose the small truncation parameter to ensure that \cref{theorem:small scales Lp} is satisfied. In contrast, here, given a truncation parameter for the ``small'' scales, we must choose a sufficiently large radius $\widetilde R$ so that \cref{theorem:small scales Lp} holds in the complement of the ball $B_{\widetilde R} (0)$.

\begin{proof}[Proof of \cref{lemma:double layer potentials far from bad ball has small norm}]
    By the John-Nirenberg inequality, let $C_1\geq 1$ be a constant such that for any $x\in\partial\Omega$ and any $s>0$ there holds
    $$
    \|\nu\|_* (B(x,s/2)) \leq C_1 \sup_{\substack{B\subset B(x,s)\\c_B\in\partial\Omega}} \avint_{B} |\nu(z)-m_B\nu | \, d\sigma(z),
    $$
    where the supremum is taken over all balls $B$ centered on $\partial\Omega$ and contained in $B(x,s)$.
    
    Let $C_*$ and $C_0$ be the constants in \cref{Semmes' decomposition,theorem:small scales Lp} respectively. As in the proof of \cref{lemma:double layer potentials in small scales has small norm}, we may assume $\delta\leq 1/(C_1C_*)$ and the lemma will follow by finding $\widetilde R = \widetilde R (t,R)$ such that for all $\xi_0\in \partial\Omega\setminus B_{\widetilde R} (0)$ there holds
    $$
     \|\characteristic_{I_0(\Delta(\xi_0,2t))} K_* (f\characteristic_{\Delta(\xi_0,2t)}) \|_{L^p(\sigma)}
    \lesssim \delta^{1/4} \|f\characteristic_{\Delta(\xi_0,2t)}\|_{L^p(\sigma)},
    $$
    with $I_0(\cdot)$ as in \rf{def:neighborhood}. By \cref{theorem:small scales Lp}, to see this it suffices to find $\widetilde R$ such that for all $\xi_0\in\partial\Omega\setminus B_{\widetilde R}(0)$, we have
    $$
    \sup_{x\in I_0 (\Delta(\xi_0,2t))} \|\nu\|_* (\Delta(x,2C_0 t)) \leq C_1 \delta,
    $$
    and the conclusions of the Semmes decomposition \cref{Semmes' decomposition} hold (with $C_1\delta$) for all $x\in I_0 (\Delta(\xi_0,2t))$ and all $0<r\leq 2C_0 t$.

    We now determine $\widetilde R$ to ensure these two conditions hold. By the John-Nirenberg inequality, the $\delta$-$(s,S;R)$ domain $\Omega$ satisfies the assumption of the Semmes decomposition \cref{Semmes' decomposition}. Recall that for every $x\in \partial\Omega\setminus B_R (0)$ we have
    $$
    \avint_{B(x,r)} |\nu(z)-m_{B(x,r)}\nu|\, d\sigma(z) \leq\delta, \text{ for all } r\in (0,\infty).
    $$
    In particular, by the John-Nirenberg inequality, we have
    \begin{equation}\label{sup in far from ball of nu at T}
    \sup_{x\in\partial\Omega\setminus B_{R+T}(0)} \|\nu\|_* (\Delta(x,T))\leq C_1\delta,
    \end{equation}
    This implies that for any compact set $\mathcal K \subset \partial\Omega \setminus B_{R+T} (0)$, we can take $R_{\widetilde{\mathcal K}_T}=T$ with $C_1\delta\leq 1/C_*$ in \cref{Semmes' decomposition}. Consequently, the value in \rf{choice of R_* in semmes' decomposition} is $R_{*,T,\mathcal K} = \delta T/C_*$, since $\diam(\partial\Omega)=\infty$.
    
    We fix $T = \max\{2C_* C_0 t/\delta,2C_0 t\}$ and we take any
    $$
    \widetilde R > R + T + 2C_0 t.
    $$
    With this choice, for any $\xi_0\in\partial\Omega\setminus B_{\widetilde R}(0)$, we have
    $$
    I_0 (\Delta(\xi_0,2t))\subset \Delta (\xi_0,2C_0t)\subset\partial\Omega\setminus B_{R+T} (0),
    $$
    and therefore,
    $$
    \sup_{x\in I_0 (\Delta(\xi_0,2t))} \|\nu\|_* (\Delta(x,2C_0 t)) \leq \sup_{x\in \partial\Omega\setminus B_{R+T} (0)} \|\nu\|_* (\Delta(x,T))\overset{\text{\rf{sup in far from ball of nu at T}}}{\leq} C_1\delta,
    $$
    and the conclusions in the Semmes decomposition \cref{Semmes' decomposition} hold (with $C_1\delta$) for all $x\in \partial\Omega\setminus B_{R+T} (0)\supset I_0 (\Delta(\xi_0,2t))$ and all $0<r\leq 2C_0 t$. This concludes the proof.
\end{proof}

\section{The double layer potential at large scales. Proof of \cref{large scale bounded by bmo norm}}\label{sec:proof of large scale bounded by bmo}

This section is dedicated to the proof of \cref{large scale bounded by bmo norm}. During this section we assume $\Omega\subset \R^{n+1}$ as in the statement of \cref{large scale bounded by bmo norm}. Briefly, $\Omega$ is a $2$-sided CAD with unbounded boundary satisfying the small Jones' $\beta_{\infty,\partial\Omega}$ coefficient and small $\bmo$ norm of $\nu$ conditions \rf{betas in large scales} and \rf{bmo norm in large scales} at large scales.

We allow the constants to depend on $n$, $1<p<\infty$ and the CAD constants of $\Omega$, and this will not be specified in the computations from now on in this section.

\subsection{Relation between \texorpdfstring{$\nu$}{nu} and the unit normal vector from \texorpdfstring{$\beta_{\infty,\partial\Omega}$}{beta}}

Given by \rf{betas in large scales}, for each $x\in\partial\Omega$ and $r\geq S$, we fix an $n$-plane $L_{B(x,r)}\ni x$ such that
\begin{equation}\label{choice of plane L in betas large scales}
\sup_{y\in\partial\Omega\cap B(x,r)} \frac{\dist(y,L_{B(x,r)})}{r}\leq 2\delta_\beta,
\end{equation}
and we denote its orthogonal unit vectors as $\pm N_{B(x,r)}$. The orientation of $N_{B(x,r)}$ is fixed below in \rf{two normal vectors are almost parallel in distance sense}.

\begin{lemma}\label{lemma:unit normal vectors betas and bmo}
	Let $x\in \partial\Omega$ and $r\geq S$. For all $y\in B(x,r)\cap \partial\Omega$,
	\begin{enumerate}
	\item $|\langle x-y, N_{B(x,r)} \rangle| \leq 2\delta_\beta r$, and
	\item $|\langle x-y, m_{B(x,r)}\nu \rangle| \lesssim \delta_\ast^{\frac{1}{2n}} r$.
	\end{enumerate}
	\begin{proof}
            The second item is proved in \rf{coro:pointwise estimate vector dot product}. For the first item, write $B=B(x,r)$ and let $N_B$ and $L_B\ni x$ be as in \rf{choice of plane L in betas large scales}. Hence $\dist(y,L_B)\leq 2\delta_\beta r$ for all $y\in B(x,r)\cap \partial\Omega$. For each $y\in B(x,r)$ let $\alpha_y = \angle (N_B, y-x)$ denote the angle between $N_B$ and $y-x$. Therefore, $|{\cos \alpha_y}|=|{\sin (\pi/2-\alpha_y)}| = \dist(y,L_B)/|x-y|\leq 2\delta_\beta r/|x-y|$, and so $|\langle x-y, N_B \rangle| = |x-y| |{\cos \alpha_y}| \leq 2\delta_\beta r$.
	\end{proof}
\end{lemma}

Under the conditions of the lemma above, as in \cite[Lemma 5.8]{David-Semmes-1991}, if $A$ is big enough, depending only on the dimension and the ADR constant, then there exists $\{y_j\}_{j=0}^{n}\subset B(x,r)\cap\partial\Omega$ with $y_0=x$ and $\dist(y_j, L_{j-1}) \geq A^{-1} r$, where $L_k$ is the $k$-plane passing through the points $y_0, \dots, y_k$. In particular $|x-y_j|\approx_A r$ and by \cref{lemma:unit normal vectors betas and bmo} we have
$$
\begin{aligned}
&\left|\left\langle\frac{x-y_j}{|x-y_j|},N_{B(x,r)}\right\rangle\right|\leq 2A \delta_\beta,
\text{ and}\\
&\left|\left\langle\frac{x-y_j}{|x-y_j|},m_{B(x,r)}\nu \right\rangle\right|\lesssim A \delta_*^{\frac{1}{2n}}.
\end{aligned}
$$
Therefore, (recall $\delta = \max\{\delta_\beta,\delta_*\}$)
\begin{equation}\label{two normal vectors are almost parallel}
|\langle m_{B(x,r)} \nu, N_{B(x,r)} \rangle| \geq 1-C_A\max\{\delta_\beta, \delta_*^{\frac{1}{2n}}\} \geq 1-C_A\delta^{\frac{1}{2n}},
\end{equation}
meaning that $N_{B(x,r)}$ and $\nu_{B(x,r)}$ are almost parallel provided both $\delta$'s are small enough. By \rf{two normal vectors are almost parallel}, choosing appropriately the orientation of $N_{B(x,r)}$, we assume that
\begin{equation}\label{two normal vectors are almost parallel in distance sense}
    |m_{B(x,r)}\nu - N_{B(x,r)} |^2 \lesssim \delta^{\frac{1}{2n}}.
\end{equation}

\subsection{Notation and parameters for the proof}

Here we write the notation we use in the following sections.
\begin{itemize}
	\item $S>0$: the first (large) scale where we have the smallness condition on the Jones' $\beta_{\infty,\partial\Omega}$ coefficient (bound given by $\delta_\beta$) in \rf{betas in large scales} and the $\bmo$ norm of the unit normal vector $\nu$ (bound given by $\delta_*$) in \rf{bmo norm in large scales}.
	\item $T\gg 100S$: scale where we truncate the kernel on the large scales.
	\item $L$, $L_B$, $L_j$, etc: planes, depending on the situation.
	\item $N_B$ or similar: unit normal vector in the Jones' $\beta_{\infty,\partial\Omega}$ coefficient in \rf{choice of plane L in betas large scales}.
	\item $M$: big parameter in the proof of \cref{lemma:large scales localization - step 1}.
	\item $\alpha$: the stopping condition in the construction of the Lipschitz graph, to obtain the Lipschitz graph with norm $\lesssim \alpha$.
	\item We fix 
    $$
    \gamma\coloneqq \frac{1}{2}\min\left\{1,p-1,\frac{1}{2n-1}\right\}>0.
    $$
    In particular, $\gamma/(1+\gamma)<1/2n$. We will use $\tau^{\frac{1}{2n}}<\tau^{\frac{\gamma}{1+\gamma}}$ when $\tau\in (0,1)$ repeatedly.
    \item Let $\M_{1+\gamma} = \M_{1+\gamma,\partial\Omega}$ denote the uncentered Hardy-Littlewood maximal function defined in \rf{eq:uncentered hardy-littlewood maximal operator}. Recall it is bounded from $L^p(\sigma)$ to $L^p(\sigma)$ with norm $C_p$, since $1+\gamma<p\leq \infty$ and $\gamma>0$ is already fixed depending on $p$.
\end{itemize}

Given $S>0$ and $\delta=\max\{\delta_\beta,\delta_*\}$ we take
\begin{equation}\label{choice of A}
A=\delta^{-\theta}\text{ with } \theta = \frac{\gamma}{2(n+3)(1+\gamma)},
\end{equation}
$M=A^{1+\frac{1}{n+1}}$, $\alpha = \delta^{\frac{\gamma}{3(1+\gamma)}}$ and $T=SA^2\gg 100S$. With this choice, we claim that, as $\delta\to 0$, there holds:
\begin{enumerate}[label=(C{\arabic*}),ref=C{\arabic*}]
	\item\label{parameter condition 1} $A\to \infty$,
	\item\label{parameter condition 2} $A/M\to 0$,
	\item\label{parameter condition 3} $\delta^{\frac{\gamma}{1+\gamma}} A\leq \delta^{\frac{1}{4n}(1-\frac{1}{\sqrt{1+\gamma}})} A^2 \leq \delta^{\frac{\gamma}{1+\gamma}} M^{n+1} A \to 0$,
	\item\label{parameter condition 4} $\delta^{\frac{\gamma}{1+\gamma}} \leq \alpha^3 \ll \alpha \to 0$,
	\item\label{parameter condition 5} $\alpha A^2\to 0$, and
	\item\label{parameter condition 6} $A S/T \to 0$.
\end{enumerate}
Indeed, \rf{parameter condition 1}, \rf{parameter condition 2}, \rf{parameter condition 4} and \rf{parameter condition 6} are immediate from the choice of parameters. For \rf{parameter condition 3} we have
$$
\delta^{\frac{\gamma}{1+\gamma}} M^{n+1} A
= \delta^{\frac{\gamma}{1+\gamma}} A^{(1+\frac{1}{n+1})(n+1)+1}
= \delta^{\frac{\gamma}{1+\gamma}-\theta (n+3)} = \delta^{\frac{1}{2}\frac{\gamma}{1+\gamma}} \to 0.
$$
For \rf{parameter condition 5},
$$
\alpha A^2 = \delta^{\frac{\gamma}{3(1+\gamma)}-2\theta} = \delta^{(\frac{1}{3}-\frac{1}{n+3})\frac{\gamma}{1+\gamma}}\to 0.
$$

\subsection{Reduction to a good lambda inequality}

In order to estimate the $L^p(\sigma)$ norm of $K_{l,*}=K_{l(T),*}$ in \rf{eq:Lp norm of large scales dlp}, it suffices to prove the following good lambda inequality
\begin{equation}\label{good lambda 1}
\sigma (\{x\in \partial\Omega : K_{l,*} f(x) > 101\lambda, \M_{1+\gamma} f(x) \leq A\lambda\})\leq c_\delta \sigma (\{x\in \partial\Omega : K_{l,*} f (x) > \lambda\}),
\end{equation}
for the parameters fixed above in \rf{choice of A} with $\delta$ small enough, and $c_\delta\to 0$ as $\delta\to 0$. A standard and routine computation using Fubini's theorem (see \cite[Theorem 1.15]{Mattila1995} for instance), the $L^p$ bound of the maximal operator $\M_{1+\gamma}$ (see \rf{eq:uncentered hardy-littlewood maximal operator}) and the good lambda inequality \rf{good lambda 1} provides that there exists $\delta_p>0$ such that if $\delta\leq \delta_p$ then
\begin{equation}\label{eq:exact Lp norm of large scales dlp}
\| K_{l,*}f\|_{L^p(\sigma)} \lesssim_{p,n} \frac{1}A\|f\|_{L^p(\sigma)}.
\end{equation}
By the choice of $A=\delta^{-\theta}$ in \rf{choice of A}, we obtain \rf{eq:Lp norm of large scales dlp} with constants depending only on the dimension and $p$. In the other case, i.e., if $\delta>\delta_p$, then \rf{eq:Lp norm of large scales dlp} follows directly from \rf{maximal operator of K is bounded}, with constants depending also on the CAD constants of $\Omega$. 

\subsection{Localization of the good lambda inequality \rf{good lambda 1}}

First, let us see the classical Whitney decomposition.

\begin{lemma}[Whitney decomposition]\label{Whitney decomposition lemma}
    If $U\subset \R^{n+1}$ is open, $U\not=\R^{n+1}$, then $U$ can be covered as
    $$
    U = \bigcup_{i\in I} Q_i,
    $$
    where $Q_i$, $i\in I$, are (classical) dyadic cubes in $\DD(\R^{n+1})$ with disjoint interiors such that the following holds:
    \begin{enumerate}
        \item $5Q_i \subset U$ for each $i\in I$.
        \item $11Q_i\cap U^c\not=\emptyset$ for each $i\in I$.
        \item If $3Q_i\cap 3Q_j\not=\emptyset$, $i,j\in I$, then $1/(7\sqrt{n+1})\leq \ell(Q_i)/\ell(Q_j)\leq 7\sqrt{n+1}$. 
        \item The family $\{3Q_i\}_{i\in I}$ has finite overlap with constant depending only on the dimension.
    \end{enumerate}
    \begin{proof}
        The family $\mathcal F = \{Q_i\}_{i\in I}$ of maximal dyadic cubes $Q\in \DD(\R^{n+1})$ such that $5Q\subset U$ satisfies the above properties by standard arguments.
    \end{proof}
\end{lemma}

For $\lambda > 0$, let $W_\lambda$ be the family of dyadic cubes given by \cref{Whitney decomposition lemma} of the open set $\Omega_\lambda \coloneqq \R^{n+1}\setminus \{x\in \partial\Omega : K_{l,*} f(x)\leq \lambda\}$. We will refer to $W_\lambda$ as the Whitney decomposition of $\Omega_\lambda$. So, we can write
$$
V_\lambda\coloneqq\partial\Omega\cap \Omega_\lambda = \{x\in \partial\Omega : K_{l,*} f (x)>\lambda\} = \bigcup_{Q\in W_\lambda} Q\cap \partial\Omega.
$$
Note that $11Q\cap (\partial\Omega\setminus V_\lambda) \not =\emptyset$ for all $Q\in W_\lambda$ and $\{3Q\}_{Q\in W_\lambda}$ has finite overlap.

\begin{notation}
    We will sometimes use $Q$ instead of $Q\cap\partial\Omega$ when it is clear from the context.
\end{notation}

\begin{lemma}\label{lemma:large scales localization}
    For $Q\in W_\lambda$, there holds
    \begin{equation}\label{eq:large scales localization}
    \{x\in Q \cap \partial\Omega : K_{l,*} f(x) > 101\lambda, \M_{1+\gamma} f(x)\leq A\lambda\}
    \subset 
    \{x\in Q \cap \partial\Omega : K_{l,*} (f\characteristic_{3Q})(x) > 50\lambda\}.
    \end{equation}
\end{lemma}

We note that if $A$ were small, the proof of \rf{eq:large scales localization} would be standard. However, its dependence on $\delta$ makes the proof more involved, relying on the flatness assumptions \rf{betas in large scales} and \rf{bmo norm in large scales}.

\begin{proof}[Proof of \cref{lemma:large scales localization}]
Let $x\in Q\cap\partial\Omega$ satisfy $K_{l,*} f(x) > 101\lambda$ and $\M_{1+\gamma} f(x)\leq A\lambda$, and let $z\in 11Q\cap (\partial\Omega\setminus V_\lambda) \not=\emptyset$ so that $K_{l,*} f(z)\leq \lambda$. We have $x,z\in 11Q\subset B(z,11(n+1)^{1/2}\ell(Q))\eqqcolon B$, and so $r(B)\approx \ell(Q)$.

The lemma will follow from
\begin{subequations}
\begin{align}
\label{lemma:large scales localization - step 1}&K_{l,*} (f\characteristic_{B^c}) (x) \leq 41\lambda,\text{ and}\\
\label{lemma:large scales localization - step 2}&K_{l,*} (f\characteristic_{B\setminus 3Q})(x)\leq 10\lambda,
\end{align}
\end{subequations}
because
$$
K_{l,*} (f\characteristic_{3Q})(x) \geq
K_{l,*} f(x) 
- K_{l,*} (f \characteristic_{B^c}) (x) 
- K_{l,*} (f \characteristic_{B\setminus 3Q})(x)\geq 50\lambda.
$$

We first prove \rf{lemma:large scales localization - step 1}. For $M>0$ big enough (auxiliary parameter only used in this proof, already defined depending on $A$, see \rf{parameter condition 2} and \rf{parameter condition 3}) write
$$
    \begin{aligned}
        K_{l,*} (f\characteristic_{B^c}) (x)
        &=  K_{l,*} (f \characteristic_{B^c}) (z) 
        + K_{l,*} (f \characteristic_{B^c}) (x)-K_{l,*} (f \characteristic_{B^c}) (z)\\
        &\leq K_{l,*} (f \characteristic_{B^c}) (z)
        +\left| K_{l,*} (f\characteristic_{MQ\setminus B})(x)- K_{l,*} (f\characteristic_{MQ\setminus B})(z)\right|\\
        &\quad +\left| K_{l,*} (f\characteristic_{(MQ)^c})(x)- K_{l,*} (f\characteristic_{(MQ)^c})(z)\right|\\
        &\eqqcolon \, \boxed{i}+\boxed{ii}+\boxed{iii}.
    \end{aligned}
$$
Since $B$ is centered at $z$ and by the choice of $z$ we have
$$
\boxed{i}\leq K_{l,*} f(z)\leq \lambda.
$$

Let us bound the term $\boxed{ii}$. We simply have
$$
\boxed{ii}
\lesssim \frac{1}{r_B^{n+1}}\left(\int_{MQ\setminus B} |\langle \nu(y),y-x \rangle| |f(y)| \, d\sigma(y)
+ \int_{MQ\setminus B} |\langle \nu(y),y-z \rangle| |f(y)| \, d\sigma(y)\right).
$$
For $\xi=x$ or $\xi=z$, we have
$$
\int_{MQ\setminus B} |\langle \nu(y),y-\xi \rangle| |f(y)| \, d\sigma(y)
\leq \int_{B(\xi, \sqrt{n+1}M\ell(Q))} |\langle \nu(y),y-\xi \rangle| |f(y)| \, d\sigma(y)
\eqqcolon \boxed{ii_\xi}.
$$
Set $B_\xi \coloneqq B(\xi, \sqrt{n+1}M\ell(Q))$. By the triangle inequality (after adding $\pm m_{B_\xi} \nu$), Hölder's inequality, \rf{coro:pointwise estimate vector dot product} and $\M_{1+\gamma} f(x)\leq A\lambda$ we have
$$
\begin{aligned}
\boxed{ii_\xi} &\leq \int_{B_\xi} |\langle \nu(y)-m_{B_\xi} \nu,y-\xi \rangle| |f(y)| \, d\sigma(y)
+ \int_{B_\xi} |\langle m_{B_\xi}\nu,y-\xi \rangle| |f(y)| \, d\sigma(y)\\
&\lesssim M^{n+1}\ell(Q)^{n+1} \left(\avint_{B_\xi} |\nu(y)-m_{B_\xi}\nu|^{\frac{1+\gamma}{\gamma}}\, d\sigma(y) \right)^{\frac{\gamma}{1+\gamma}}
\left(\avint_{B_\xi} |f(y)|^{1+\gamma}\, d\sigma(y)\right)^{\frac{1}{1+\gamma}}\\
&\quad + \int_{B_\xi} |\langle m_{B_\xi}\nu,y-\xi \rangle| |f(y)| \, d\sigma(y)\\
\overset{\text{\rf{coro:pointwise estimate vector dot product}}}&{\lesssim} M^{n+1}\ell(Q)^{n+1} \delta_*^{\frac{\gamma}{1+\gamma}}
\M_{1+\gamma}f(x) 
+ M^{n+1}\ell(Q)^{n+1} \delta_*^{\frac{1}{2n}} \M_1 f(x)\\
&\lesssim M^{n+1}\ell(Q)^{n+1} \delta_*^{\frac{\gamma}{1+\gamma}}A\lambda.
\end{aligned}
$$
Note that even though we are studying the term $\boxed{ii_\xi}$, in any case we controlled $\avint_{B_\xi} |f|\, d\sigma$ and $(\avint_{B_\xi} |f|^{1+\gamma}\, d\sigma)^{1/(1+\gamma)}$ by $\M_1 f(x)$ and $\M_{1+\gamma} f(x)$ respectively, and so we finally used that $\M_{1+\gamma} f(x)\leq A\lambda$. Hence,
$$
    \boxed{ii}
    \lesssim \frac{1}{r_B^{n+1}}M^{n+1}\ell(Q)^{n+1} \delta_*^{\frac{\gamma}{1+\gamma}}A\lambda \approx M^{n+1} \delta_*^{\frac{\gamma}{1+\gamma}} A\lambda.
$$

Let us bound $\boxed{iii}$. For $\varepsilon\geq T$, we have
$$
|K_\varepsilon (f\characteristic_{(MQ)^c})(x)-K_\varepsilon (f\characteristic_{(MQ)^c})(z)|\leq \boxed{A}+\boxed{B}+\boxed{C},
$$
where
$$
\begin{aligned}
    \boxed{A}&\coloneqq \left|\int_{\substack{|y-x|>\varepsilon\\|y-z|>\varepsilon}} \left\langle \nu(y)f(y)\characteristic_{(MQ)^c}(y),\frac{y-x}{|x-y|^{n+1}}-\frac{y-z}{|z-y|^{n+1}}\right\rangle \, d\sigma(y)\right|,\\
    \boxed{B}&\coloneqq \left|\int_{\substack{|y-x|>\varepsilon\\|y-z|\leq \varepsilon}} \frac{\langle \nu(y),y-x\rangle}{|x-y|^{n+1}} f(y)\characteristic_{(MQ)^c}(y) \, d\sigma(y)\right|,\text{ and}\\
    \boxed{C}&\coloneqq \left|\int_{\substack{|y-z|>\varepsilon\\|y-x|\leq \varepsilon}} \frac{\langle \nu(y),y-z\rangle}{|z-y|^{n+1}} f(y)\characteristic_{(MQ)^c}(y) \, d\sigma(y)\right|.
\end{aligned}
$$
We first study the term $\boxed{A}$. Using the cancellation property of the Calderón-Zygmund kernel first (this is just the mean value theorem), and $B(x,M\ell(Q)/2)\subset MQ$ and $|x-z|\lesssim\ell(Q)$ second, we have
$$
\boxed{A}
\lesssim \int_{y\not\in MQ} |f(y)| \frac{|x-z|}{|x-y|^{n+1}} \, d\sigma(y)
\lesssim \ell(Q) \int_{y\not\in B(x,\frac{M}{2}\ell(Q))} \frac{|f(y)|}{|x-y|^{n+1}} \, d\sigma(y).
$$
Note that the latter integral does not depend on $z$. Breaking the latter integral into annuli and by the AD regularity of $\sigma$, using standard estimates we get
$$
\boxed{A}
\lesssim \frac{1}{M} \M_{1+\gamma} f(x)
\leq \frac{A}{M} \lambda.
$$
We now turn to the terms $\boxed{B}$ and $\boxed{C}$. First we note that we can assume that $\varepsilon\geq \frac{M-10}{2}\ell(Q)\geq \frac{M}{4}\ell(Q)$, otherwise $B(z,\varepsilon)\subset MQ$ and so $\boxed{B}=0$. For $\varepsilon >\frac{M}{4}\ell(Q)$ we have $B(z,\varepsilon)\subset B(x,2\varepsilon)$, and we get
$$
\begin{aligned}
\boxed{B}
\leq& \int_{B(x,2\varepsilon)\setminus\overline{B(x,\varepsilon)}} \frac{|\langle \nu(y),y-x\rangle|}{|x-y|^{n+1}} |f(y)| \, d\sigma(y)
\lesssim \frac{1}{\varepsilon} \avint_{B(x,2\varepsilon)} |\langle \nu(y),y-x\rangle| |f(y)| \, d\sigma(y)\\
\leq & \, \frac{1}{\varepsilon} \avint_{B(x,2\varepsilon)} |\langle \nu(y) - m_{B(x,2\varepsilon)}\nu,y-x\rangle| |f(y)| \, d\sigma(y)\\
& + \frac{1}{\varepsilon} \avint_{B(x,2\varepsilon)} |\langle m_{B(x,2\varepsilon)}\nu,y-x\rangle| |f(y)| \, d\sigma(y)\\
\lesssim & \left(\avint_{B(x,2\varepsilon)} |\nu(y) - m_{B(x,2\varepsilon)}\nu|^{\frac{1+\gamma}{\gamma}} \, d\sigma(y)\right)^{\frac{\gamma}{1+\gamma}}
\left(\avint_{B(x,2\varepsilon)} |f(y)|^{1+\gamma}\, d\sigma(y)\right)^{\frac{1}{1+\gamma}}\\
&+\left(\sup_{y\in \partial\Omega\cap B(x,2\varepsilon)} \frac{|\langle x-y, m_{B(x,2\varepsilon)}\nu\rangle|}{\varepsilon}\right)
\left(\avint_{B(x,2\varepsilon)} |f(y)|^{1+\gamma}\, d\sigma(y)\right)^{\frac{1}{1+\gamma}}\\
\overset{\mathclap{\text{\rf{coro:pointwise estimate vector dot product}}}}{\lesssim}&\ (\delta_*^{\frac{\gamma}{1+\gamma}} + \delta_*^{\frac{1}{2n}}) \M_{1+\gamma} f(x) \lesssim \delta_*^{\frac{\gamma}{1+\gamma}} \M_{1+\gamma} f(x)
\leq \delta_*^{\frac{\gamma}{1+\gamma}} A\lambda.
\end{aligned}
$$
By symmetry, but using that $B(z,\varepsilon)\subset B(x,2\varepsilon)$ (recall we can assume $\varepsilon>\frac{M}{4}\ell(Q)$) before applying the bound of the maximal operator $\M_{1+\gamma}$, we have the same bound for $\boxed{C}$, that is, $\boxed{C}\lesssim \delta_*^{\frac{\gamma}{1+\gamma}} \M_{1+\gamma} f(x)\leq \delta_*^{\frac{\gamma}{1+\gamma}} A\lambda$. Therefore, the term $\boxed{iii}$ is controlled by
$$
\boxed{iii}\leq \boxed{A}+\boxed{B}+\boxed{C}
\lesssim \left(\frac{A}{M} + \delta_*^{\frac{\gamma}{1+\gamma}}A\right)\lambda
$$

From the bound of $\boxed{i}$, $\boxed{ii}$ and $\boxed{iii}$, and conditions \rf{parameter condition 2} and \rf{parameter condition 3} we conclude
$$
K_{l,*} (f\characteristic_{B^c}) (x)\leq \lambda + C\left(M^{n+1}\delta_*^{\frac{\gamma}{1+\gamma}} A + \frac{A}{M} + \delta_*^{\frac{\gamma}{1+\gamma}}A\right)\lambda
\leq 41\lambda,
$$
as claimed in \rf{lemma:large scales localization - step 1}.

It remains to prove \rf{lemma:large scales localization - step 2}. For $\varepsilon\geq T$, we may assume $\{y:|y-x|>\varepsilon\}\cap B \not=\emptyset$, otherwise $K_{\varepsilon} (f \characteristic_{B\setminus 3Q})(x)=0$. Thus, $r(B)\gtrsim\varepsilon\geq T$, allowing us to apply \rf{coro:pointwise estimate vector dot product}. Using that $\ell(Q)\approx r(B)$ and that $|y-x|\gtrsim \ell(Q)$ if $y\not\in 3Q$, and arguing in the last inequality as in the bound for $\boxed{ii_\xi}$, we get
$$
    \begin{aligned}
    K_{\varepsilon} (f \characteristic_{B\setminus 3Q})(x)
    &=\left|\int_{|y-x|>\varepsilon} \frac{\langle\nu(y),y-x\rangle}{|x-y|^{n+1}} f(y) \characteristic_{B\setminus 3Q}(y) \, d\sigma(y)\right|\\
    &\lesssim \frac{1}{\ell(Q)} \left(\avint_B |\langle\nu(y)-m_B \nu,y-x\rangle| |f(y)| \, d\sigma(y) 
    +\avint_B |\langle m_B \nu, y-x \rangle| |f(y)| \, d\sigma(y)\right)\\
    &\lesssim \delta_*^{\frac{\gamma}{1+\gamma}} A\lambda,
    \end{aligned}
$$
uniformly in $\varepsilon\geq T$. Then \rf{lemma:large scales localization - step 2} follows by condition \rf{parameter condition 3}.
\end{proof}

By \rf{eq:large scales localization}, in order to prove \rf{good lambda 1} we claim that it suffices to see 
\begin{equation}\label{good lambda 2}
	\sigma(\{x\in Q: K_{l,*} (f\characteristic_{3Q})(x) > 50\lambda, \M_{1+\gamma} f (x) \leq A\lambda\}) \leq c_\delta \sigma(3Q),
\end{equation}
for all $Q\in W_\lambda$. Indeed,
\begin{multline*}
	\sigma(\{ x\in\partial\Omega :  K_{l,*} f(x) > 101\lambda, \M_{1+\gamma} f(x) \leq A\lambda\}) \\
	\begin{aligned}
	&\leq \sum_{Q\in W_\lambda} \sigma(\{x\in Q: K_{l,*} f(x) > 101\lambda, \M_{1+\gamma} f(x) \leq A\lambda\})\\
	\overset{\text{\rf{eq:large scales localization}}}&{\leq} \sum_{Q\in W_\lambda} \sigma(\{x\in Q: K_{l,*} (f\characteristic_{3Q})(x) > 50\lambda, \M_{1+\gamma} f(x) \leq A\lambda\})\\
	\overset{\text{\rf{good lambda 2}}}&{\leq} c_\delta \sum_{Q\in W_\lambda} \sigma (3Q) \lesssim c_\delta \sigma (\Omega_\lambda) = c_\delta \sigma(V_\lambda),
	\end{aligned}
\end{multline*}
where in the last step we used the finite overlap of the family $\{3Q\}_{Q\in W_\lambda}$.

\begin{rem}\label{rem:we can assume cube in good lambda is big}
    From now on we can assume without loss of generality that $\ell(Q)\geq \frac{T}{10(n+1)^{1/2}}$, otherwise, if $x\in Q$ then $\{y \in 3Q : |y-x|> T\}=\emptyset$ and hence
    $$
    K_{l,*} (f\characteristic_{3Q}) (x) = \sup_{\varepsilon\geq T} \left|\int_{y\in 3Q\setminus \overline{B(x,\varepsilon)}}\frac{\langle\nu(y),y-x\rangle}{|x-y|^{n+1}} f(y) \characteristic_{3Q}(y)\, d\sigma(y)\right|=0.
    $$
\end{rem}

\subsection{Construction of a Lipschitz graph}\label{sec:construction of Lip graph}

During this section, let $\Omega\subset\R^{n+1}$ be an ADR domain with unbounded boundary, and let $\DD_\sigma$ denote its dyadic lattice in \cref{lemma:def of boundary dyadic cubes}. For $t > 1$ and $Q\in \DD_\sigma$ we set
\begin{equation}\label{dilated boundary cube}
t Q \coloneqq \{x\in \partial\Omega : \dist (x,Q)\leq (t-1)\diam\, Q\}.
\end{equation}

\begin{notation}[Dependence of parameters]
    Fixed $m>1$, we choose $k_0=k_0(m)\gg m$ and $k=k(k_0)\gg k_0$. We consider $\varepsilon,\alpha\ll 1$ such that $\varepsilon\ll \alpha$ and $\varepsilon \ll 1/k_0$. Throughout this section, we allow the constant $C$ to depend on the ADR constant, but not on the parameters mentioned above.
\end{notation}

Following \cite[Section 8]{David-Semmes-1991}, in this section we construct a Lipschitz approximating graph at large scales (see \cref{prop:construction Lip graph}), which only uses the following flatness assumption at large scales: We assume that for each cube $Q\in \DD_\sigma$ with $\diam\, Q \geq S$, there exists an $n$-plane $L_Q$ such that
\begin{equation}\label{line approx cube}
\dist(x,L_Q) \leq \varepsilon \diam\, Q \text{ for all }x\in kQ.
\end{equation}
We will see in \cref{rem:candidate of approx plane of cube} in \cref{sec:proof of the localized good lambda ineq} below that this is guaranteed under the assumption \rf{betas in large scales}.

We recall that the $n$-plane in \rf{line approx cube} is almost unique in the following sense.

\begin{lemma}[See {\cite[Lemma 5.13]{David-Semmes-1991}}]\label{lemma:DS91 lemma 5.13}
    Let $Q\in\DD_\sigma$ and assume that $L_1$ and $L_2$ are two $n$-planes such that $\dist(x,L_i)\leq \varepsilon^\prime \diam\, Q$ for all $x\in Q$, $i=1,2$. Then $\angle(L_1,L_2)\leq C\varepsilon^\prime$.
\end{lemma}

Let us fix $R\in \DD_{\sigma,j_0}$ with $j_0\in \Z$ so that $2^{-j_0}\geq S$; in particular $\diam\, Q\geq S$ for all $Q\in \DD_{\sigma,j_0}$. Fixed $m>1$, we define the ``$m$-neighborhood cubes'' of $R$ as the collection of cubes in $\DD_{\sigma,j_0}$ touching $mR$. That is,
$$
U_m(R) \coloneqq \{Q \in \DD_{\sigma,j_0} : Q\cap mR \not = \emptyset\}.
$$
So, for every $Q\in U_m (R)$ there holds $\dist(Q,R)\leq (m-1)\diam\, R$. For any two cubes $Q_1,Q_2\in U_m (R)$, if $x\in Q_1$ then
$$
\dist(x,Q_2) \leq \diam\, Q_1 + \dist(Q_1,R) + \diam\, R + \dist(R,Q_2) \leq C_\DD 2m \diam \, Q_2,
$$
and therefore
$$
Q_1 \subset (C_\DD 2m+1)Q_2.
$$
Fix a constant $k\gg m$ (at least $k\geq C_\DD 2m+1$), so that the inclusion above implies
$$
Q_1\subset kQ_2 \text{ for any }Q_1,Q_2\in U_m (R).
$$
As a consequence of this, \rf{line approx cube} and \cref{lemma:DS91 lemma 5.13}, we have
\begin{equation}\label{eq:angle between cubes in UmR}
\angle(L_{Q_1},L_{Q_2})\lesssim \varepsilon\text{ for any }Q_1,Q_2\in U_m (R).
\end{equation}

\begin{definition}\label{stopping conditions}
	We say that $Q\in \bigcup_{j\geq j_0} \DD_{\sigma,j}$ is a stopping cube with respect to $R$, and we write $Q\in \Stop(R)$, if $Q\subset \bigcup_{Q^\prime \in U_m (R)} Q^\prime$ and $Q$ is maximal such that either
	\begin{enumerate}
		\item $\diam\, Q<S$, or
		\item $\diam\, Q\geq S$ and $\angle (L_R, L_Q) > \alpha$, recall $L_R$ and $L_Q$ as in \rf{line approx cube}.
	\end{enumerate}
\end{definition}

\begin{definition}
	We define the set $\tree(R)$ as the collection of cubes $P\in \bigcup_{j\geq j_0} \DD_{\sigma,j}$ such that
        \begin{enumerate}
            \item $P\subset \bigcup_{Q^\prime \in U_m (R)} Q^\prime$, and
            \item $P\not\subset Q$ for any $Q\in \Stop(R)$.
        \end{enumerate}
        In particular, every $Q\in \tree (R)$ satisfies $\diam\, Q\geq S$ and $\angle (L_R,L_Q)\leq \alpha$.
\end{definition}

\begin{rem}
    In contrast to the construction in \cite[Section 8]{David-Semmes-1991}, where only dyadic subcubes of $R$ are considered, here we take dyadic subcubes of $U_m(R)$ in the construction of $\tree(R)$. In our situation, since $\varepsilon\ll\alpha$ and $\diam\, Q\geq S$ for all $Q\in U_m(R)$, we note that \rf{eq:angle between cubes in UmR} ensures that $Q\not\in\Stop(R)$ for all $Q\in U_m (R)$. Consequently, given a cube in $\Stop(R)$, its parent is in $\tree(R)$.
\end{rem}

Consider
$$
d(x) = d_R (x) \coloneqq \inf_{Q\in \tree(R)} \{\dist(x,Q) + \diam\, Q\}, \quad x\in \R^{n+1}.
$$ 
Since $\# \tree(R)<\infty$, the infimum is in fact a minimum. Note that $d(\cdot)\geq S$ and is $1$-Lipschitz, as it is the infimum of $1$-Lipschitz functions. Clearly, if $x\in Q\in \tree(R)$ then $d(x)\leq \diam\, Q$.

In this section, we prove the following result, which is the analog of \cite[Proposition 8.2]{David-Semmes-1991}. In this lemma and throughout this section, $L_R^\perp$ denotes the orthogonal line to the $n$-plane $L_R$ (see \rf{line approx cube}), $\Pi = \Pi_R$ denotes the orthogonal projection onto $L_R$, and $\Pi^\perp = \Pi^\perp_R$ denotes the orthogonal projection onto $L_R^\perp$.

\begin{propo}\label{prop:construction Lip graph}
    For $R\in \DD_{\sigma,j_0}$ with $j_0\in\Z$ so that $2^{-j_0}\geq S$ (in particular $\diam\, R\geq S$), there exists $A:L_R\to L_R^\perp$ Lipschitz graph with norm $\leq C\alpha$ such that
    \begin{equation}\label{eq:construction Lip graph}
    \dist\left(x,\left(\Pi_R(x),A(\Pi_R(x))\right)\right)\leq C \varepsilon d_R(x) \text{ for all } x\in k_0 R.
    \end{equation}
\end{propo}

We define in $L_R$ the $1$-Lipschitz function
\begin{equation}\label{def of D}
D(p) = D_R(p)\coloneqq \inf_{x\in \Pi^{-1}(p)} d_R (x), \quad p\in L_R,
\end{equation}
which can be rewritten as
$$
D (p) = \inf_{x\in \Pi^{-1}(p)} \inf_{Q\in \tree(R)} \{\dist(x,Q) + \diam\, Q\}
=\inf_{Q\in\tree(R)} \{\dist(p,\Pi_R (Q))+\diam\, Q\}.
$$
As $d \geq S$ in $\R^{n+1}$, in particular $D \geq S$ in $L_R$. As before, since $\# \tree(R)<\infty$, the latter infimum is in fact a minimum.

Arguing as in the proof of \cite[Lemma 8.4]{David-Semmes-1991}, we obtain the following, as the proof relies only on the fact that every $Q\in \tree(R)$ satisfies $\angle (L_R, L_Q)\leq \alpha$.

\begin{lemma}[{\cite[Lemma 8.4]{David-Semmes-1991}}]\label{far points imply flat graph}
    If $x,y\in 10 k_0 R$ satisfy $|x-y|\geq 10^{-3}\min \{d(x),d(y)\}$, then
    $$
    |\Pi^\perp (x)-\Pi^\perp (y)|\leq 2\alpha |\Pi(x)-\Pi(y)|.
    $$
\end{lemma}

Now we decompose $L_R$ in classical dyadic cubes using the function $D_R$. First, we shall identify $L_R$ with $\R^n$, and in particular we equip $L_R$ with classical dyadic cubes $\DD_{L_R} \coloneqq \DD_{\R^n}$. For each $x\in L_R$ (recall $D_R \geq S >0$ in $L_R$), let $R_x \in \DD_{L_R}$ be the largest dyadic cube containing $x$ and satisfying
\begin{equation}\label{choice of dyadic decomposition of hyperplane}
\diam\, R_x \leq 20^{-1} \inf_{u\in R_x} D(u).
\end{equation}
We relabel them without repetition as $\{R_i\}_{i\in I}$. Thus, the family $\{R_i\}_{i\in I}$ is pairwise disjoint and covers $L_R$. (As in \cite[Section 8]{David-Semmes-1991}, we use the convention that dyadic cubes are closed but are called disjoint if their interiors are disjoint.)

Recall that the definition of $tQ$ is slightly different depending on whether $Q\in \DD_\sigma$ or $Q\in \DD_{L_R}=\DD_{\R^n}$, see \rf{dilated boundary cube} and \cref{sec:notation} respectively.

From the definition of the family $\{R_i\}_{i\in I}$ using \rf{choice of dyadic decomposition of hyperplane} and that $D$ is $1$-Lipschitz in $L_R$, by the same proof of \cite[Lemma 8.7]{David-Semmes-1991} we have:

\begin{lemma}[{\cite[Lemma 8.7]{David-Semmes-1991}}]
    For all $y\in 10R_i$, $i\in I$,
    \begin{equation}\label{comparability D and diamRi}
        10\diam\, R_i \leq D(y) \leq 60\diam\, R_i,
        \end{equation}
    and in particular, if $10R_i \cap 10R_j \not=\emptyset$, $i,j\in I$, then
    \begin{equation}\label{close dyadic cubes have comparable diameter}
    6^{-1} \diam\, R_i \leq \diam\, R_j \leq 6\diam\, R_i.
    \end{equation}
\end{lemma}

Let $U_0 \coloneqq L_R \cap B(\Pi(x_R),2k_0 \diam\, R)$, where $x_R$ is any fixed point in $R$, and $I_0 \coloneqq \{i\in I : R_i \cap U_0 \not=\emptyset\}$. We fix $k_0\gg m$ so that $\bigcup_{Q^\prime \in U_m(R)} \Pi(Q^\prime) \subset U_0$. For each $i\in I_0$, let $Q(i)\in \tree(R)$ be such that\footnote{The correspondence $I_0 \ni i \mapsto Q(i)$ may not be injective.}
\begin{subequations}
\begin{align}
    \label{choice 1 of Qi}&C^{-1}\frac{1}{k_0} \diam\, R_i \leq \diam\, Q(i) \leq C\diam\, R_i, \text{ and}\\
    \label{choice 2 of Qi}&\dist (\Pi (Q(i)),R_i) \leq C \diam\, R_i.
\end{align}
\end{subequations}
Such cubes $Q(i)$ exist, because if $p\in R_i$, then there exists $Q\in \tree (R)$ such that
$$
\dist(p,\Pi(Q)) + \diam\, Q \leq 2D(p) \overset{\text{\rf{comparability D and diamRi}}}{\leq} 120 \diam\, R_i,
$$
and we can take $Q(i)$ to be the maximal cube in $\tree(R)$ with $Q\subset Q(i)$ and $\diam\, Q(i)\leq 120\diam\, R_i$.

\begin{definition}
    For $i\in I_0$, let $A_i : L_R \to L_R^\perp$ denote the affine function whose graph is the $n$-plane $L_{Q(i)}$. By the stopping condition, the Lipschitz norm is $\leq 2\alpha$.
\end{definition}

We consider a partition of unity for $V=\bigcup_{i\in I_0} 2R_i$. That is, a family $\{\phi_i\}_{i\in I_0}$ satisfying for all $i\in I_0$ that $\phi_i\geq 0$, $\phi_i\in C^1_c(3R_i)$ (hence $\{\supp\phi_i\}_{i\in I_0}$ has finite overlap, by \cref{close dyadic cubes have comparable diameter}) and
$$
|\nabla\phi_i|\lesssim \frac{1}{\diam\, R_i}.
$$

Finally we define $A$ on $V$ by
$$
A(p) \coloneqq \sum_{i\in I_0} \phi_i (p) A_i (p) \text{ for }p\in V.
$$

By the same proof in \cite[Lemma 8.17]{David-Semmes-1991} we have:

\begin{lemma}[{\cite[Lemma 8.17]{David-Semmes-1991}}]\label{bound of affine maps in terms of diameter}
            If $10R_i\cap 10R_j \not = \emptyset$, then $\dist(Q(i),Q(j))\leq C\diam\, R_j$ and
            \begin{equation}\label{same point but two affine maps}
            |A_i (q)-A_j (q)|\leq C\varepsilon \diam\, R_j \text{ for all } q\in 100R_j.
            \end{equation}
        \end{lemma}

Using the previous lemma, the same proof in \cite[Section 8, p.~46]{David-Semmes-1991} applies to see that the restriction of $A$ to $2R_j$, $j\in I_0$, is $3\alpha$-Lipschitz.

\begin{lemma}[{\cite[(8.19)]{David-Semmes-1991}}]\label{lemma:3alpha for close points}
    If $p,q\in 2R_j$, $j\in I_0$, then
    $$
    |A(p)-A(q)|\leq 3\alpha|p-q|.
    $$
\end{lemma}

We now aim to show that $A$ is $C\alpha$-Lipschitz in $U_0$. We remark that this is the analog of \cite[(8.20)]{David-Semmes-1991}, with the difference that in our case, we have $D(\cdot) \geq S >0$ in $L_R$. For the sake of completeness, we provide the full details of the proof.

\begin{lemma}\label{Lip graph is Calpha}
    If $p,q\in U_0$, then $|A(p)-A(q)| \leq C\alpha |p-q|$.
    \begin{proof}
        Since $D\geq S>0$, all points in $U_0$ belong to some (unique) $R_k$, where $k\in I_0$. Let $i,j\in I_0$ be such that $p\in R_i$ and $q\in R_j$. If $p,q\in 2R_i$ or $p,q\in 2R_j$, then $|A(p)-A(q)|\leq 3\alpha |p-q|$ by \cref{lemma:3alpha for close points}, and we are done. Hence, from now on, we may assume that $q\not\in 2R_i$ and $p\not\in 2R_j$. In particular, this implies
        \begin{equation}\label{case points are far for each dyadic cube}
        |p-q|\geq \max\{\diam\, R_i,\diam\, R_j\}.
        \end{equation}
        
        Let $y\in Q(i)$ and $z\in Q(j)$ be such that $|y-z|=\sup_{a\in Q(i),b\in Q(j)} |a-b|$, which in particular implies
        \begin{equation}\label{choice of y and z consequence}
        |y-z|\geq \frac{1}{2}\min\{\diam\, Q(i),\diam\, Q(j)\}.
        \end{equation}
        We have
        $$
        \begin{aligned}
        |A(p)-A(q)| \leq&\, |A(p)-A_i(p)| + |A_i (p) - A_i (\Pi(y))| + |A_i(\Pi(y))-\Pi^\perp (y)|\\
        & +|\Pi^\perp (y)- \Pi^\perp (z)|\\
        & +|A(q)-A_j(q)| + |A_j (q) - A_j (\Pi(z))| + |A_j(\Pi(z))-\Pi^\perp (z)|.
        \end{aligned}
        $$
        Let us see that all the terms are bounded by $\lesssim \alpha$. We bound the first four terms, and the last three follow by symmetry.
        
        Term $|A(p)-A_i(p)|$: Using the partition of unity in the second equality, and by \rf{same point but two affine maps} in \cref{bound of affine maps in terms of diameter} (since $p\in R_i$, the fact that $\phi_j (p)\not = 0$ implies that $\supp \phi_j \cap R_i \not = \emptyset$, and in particular $3R_j\cap R_i\not = \emptyset$) and since the last sum runs over the index $j\in I_0$ around $R_i$ (in particular a finite number of candidates), we have
        \begin{equation}\label{compare Lip graph with one hyperplane}
        \begin{aligned}
        |A(p)-A_i (p)| &= \left|\left(\sum_{j\in I_0} \phi_j (p) A_j (p)\right) - A_i (p)\right|
        = \left|\sum_{j\in I_0} \phi_j (p)(A_j (p) - A_i (p))\right|\\
        \overset{\text{\rf{same point but two affine maps}}}&{\leq} C\varepsilon \diam\, R_i
        \overset{\text{\rf{case points are far for each dyadic cube}}}{\leq} C\varepsilon |p-q| \leq \alpha |p-q|.
        \end{aligned}
        \end{equation}
        
        Term $|A_i (p) - A_i (\Pi(y))|$: First, $|A_i (p)-A_i (\Pi(y))| \leq 2\alpha |p-\Pi(y)|$ since $A_i$ is $2\alpha$-Lipschitz, and second $|p-\Pi(y)|\leq \diam\, R_i + \dist(R_i,\Pi(Q(i)))+\diam\, Q(i)$, which by the choice of $Q(i)$ in \rf{choice 1 of Qi} and \rf{choice 2 of Qi} the last two terms are controlled by $\lesssim \diam\, R_i$. We conclude this term by using \rf{case points are far for each dyadic cube}.
        
        Term $|A_i(\Pi(y))-\Pi^\perp (y)|$: Since $\angle (L_{Q(i)},L_R)\leq \alpha$ is small and by \rf{line approx cube}, we have $|A_i (\Pi(y))-\Pi^\perp (y)|\leq 2\dist(y,L_{Q(i)}) \leq 2\varepsilon \diam\, Q(i)$. By the choice of $Q(i)$ in \rf{choice 1 of Qi} this is controlled by $\leq C\varepsilon \diam\, R_i$, and by \rf{case points are far for each dyadic cube} the last term is controlled by $\leq C\varepsilon |p-q|\leq \alpha |p-q|$.
        
        Term $|\Pi^\perp (y)- \Pi^\perp (z)|$: First, $y\in Q(i)$ implies that $d(y)\leq \diam\, Q(i)$, and $z\in Q(j)$ implies $d(z)\leq \diam\, Q(j)$. Thus, from \rf{choice of y and z consequence} and this we have $|y-z|\geq \frac{1}{2}\min\{\diam\, Q(i),\diam\, Q(j)\} \geq \frac{1}{2} \min\{d(y),d(z)\}$. Using \cref{far points imply flat graph} in the first inequality we get
        $$
        |\Pi^\perp (y)-\Pi^\perp (z)|\leq 2\alpha |\Pi (y)-\Pi(z)|
        \leq 2\alpha \left(|\Pi (y)-p|+|p-q|+|q-\Pi(z)|\right).
        $$
        As in the bound of the second term (that is, $|A_i (p) - A_i (\Pi(y))|$), the first term inside the brackets is bounded by $\leq C\diam\, R_i$. By symmetry, the third term is controlled by $\leq C\diam\, R_j$. This term follows by \rf{case points are far for each dyadic cube} as well.
    \end{proof}
\end{lemma}

Now we can use the Whitney extension theorem to extend $A$ from $U_0$ to a $C\alpha$-Lipschitz function on all $L_R$. This gives the existence of the $C\alpha$-Lipschitz graph in \cref{prop:construction Lip graph}. It remains to see now that this $C\alpha$-Lipschitz graph approximates $k_0 R$ in the \rf{eq:construction Lip graph} sense.

\begin{proof}[Proof of \rf{eq:construction Lip graph}]
    Let $x\in k_0 R$ and set $p=\Pi (x)$. Recall $D(p)\geq S>0$, and let $R_i$, $i\in I_0$, so that $p\in R_i$. We first break
    $$
    \dist(x,(\Pi(x),A(\Pi(x)))) \leq |\Pi^\perp (x) - A_i (\Pi(x))|+|A (\Pi(x)) - A_i (\Pi(x))|.
    $$
    
    Applying \cite[Lemma 8.21]{David-Semmes-1991} with $r=D(p)$ and $Q=Q(i)$, we get that $x\in C_{k_0} Q(i)$ for some $C_{k_0}$ depending on $k_0$, so taking $k\gg C_{k_0}$ we have $x\in kQ(i)$. Using that $\angle (L_R, L_{Q(i)})\leq \alpha$ is small, the flatness condition in \rf{line approx cube} and the choice of $Q(i)$ in \rf{choice 1 of Qi} respectively, we have
    $$
    |\Pi^\perp (x) - A_i (\Pi(x))|
    \leq 2 \dist (x,L_{Q(i)})
    \leq 2\varepsilon \diam\, Q(i)
    \lesssim \varepsilon \diam\, R_i.
    $$
    
    The term $|A (\Pi(x)) - A_i (\Pi(x))|$ is simply the first term in the proof of \cref{Lip graph is Calpha}, whence we obtain directly from \rf{compare Lip graph with one hyperplane} that
    $$
    |A (\Pi(x)) - A_i (\Pi(x))| \leq C \varepsilon\diam\, R_i.
    $$

    All in all, $\dist(x,(\Pi(x),A(\Pi(x))))\leq C\varepsilon \diam\, R_i$. From the definition of the family $\{R_i\}_{i\in I}$ in \rf{choice of dyadic decomposition of hyperplane}, we have $\diam\, R_i \leq 20^{-1} D(p)$. Just from the definition $D$ in \rf{def of D}, $D(p)\leq d(x)$. This concludes the proof of \rf{eq:construction Lip graph}.
\end{proof}

\subsection{Proof of the localized good lambda inequality \rf{good lambda 2}}\label{sec:proof of the localized good lambda ineq}

Let $\Omega$ be as in the statement of \cref{large scale bounded by bmo norm}. To apply the notation of the section of the Lipschitz graph construction, for $\lambda > 0$ we write $R_W\in W_\lambda$ for the cube from the Whitney covering of $\Omega_\lambda$ (the letters $Q,R$ are reserved for boundary dyadic cubes $\DD_\sigma$), and we rewrite the good lambda inequality \rf{good lambda 2} as
\begin{equation}\label{good lambda 3}
	\sigma(\{x\in R_W: K_{l,*} (f\characteristic_{3R_W})(x) > 50\lambda, \M_{1+\gamma} f (x) \leq A\lambda\}) \leq c_\delta \sigma(3R_W).
\end{equation}

We will use without mention that we assume that there exists $x_0\in R_W$ with $\M_{1+\gamma} f(x_0)\leq A\lambda$, otherwise the left-hand side in \rf{good lambda 3} is zero and we are done. Recall also from \cref{rem:we can assume cube in good lambda is big} that we have $\ell(R_W) \gtrsim T \gg S$.

The thin boundary condition \rf{boundary dyadic cube thin boundary condition} in \cref{lemma:def of boundary dyadic cubes} implies, for each $Q\in\DD_\sigma$, the existence of a point $c_Q \in Q$, called the center of $Q$, such that $\dist(c_Q, \partial\Omega\setminus Q)\gtrsim_{\text{ADR},n} \diam\, Q$, see \cite[Lemma 3.5 of Part I]{David-Semmes-1993-AMS}. We denote
\begin{equation}\label{ball inside outside the cube}
B(Q) \coloneqq B(c_Q, c_1 \diam\, Q), \quad B_Q \coloneqq B(c_Q, \diam\, Q) \supset Q,
\end{equation}
where $c_1 \in (0,1]$ is chosen so that $B(Q^\prime)\cap \partial\Omega \subset Q^\prime$ for all $Q^\prime\in\DD_\sigma$, satisfying also that $B(Q_1)\cap B(Q_2)=\emptyset$ for all disjoint $Q_1,Q_2\in \DD_\sigma$. Indeed, assuming without loss of generality that $\diam\, Q_2\leq \diam \, Q_1$, if $B(Q_1)\cap B(Q_2)\not=\emptyset$, then
$$
\diam\, Q_1 \lesssim_{\text{ADR},n} \dist(c_{Q_1},\partial\Omega\setminus Q_1) \leq |c_{Q_1}-c_{Q_2}|\leq c_1 (\diam\, Q_1 + \diam\, Q_2)\leq 2c_1 \diam\, Q_1,
$$
which is not possible if $c_1$ is chosen to be sufficiently small, depending only on the ADR character of $\partial\Omega$.

We begin the proof of \rf{good lambda 3}. Given $R_W\in W_\lambda$, let $j\in \Z$ be so that $\ell (R_W) = 2^{-j}$ (recall $R_W\in \DD(\R^{n+1})$), fix any $x\in R_W\cap\partial\Omega$ (if such point does not exist then the left-hand side of \rf{good lambda 3} is zero), and let $R \in \DD_{\sigma,j}$ be the cube with $x \in R$. We fix $m=3$ so that $\partial\Omega \cap 3R_W\subset mR$. Indeed, since $\diam\, R \geq 2^{-j}$, for all $y\in R_W$ we have $|y-x|\leq 2\diam\, R_W = 2\cdot 2^{-j}\leq 2\diam\, R$, thus $\dist(y,R)\leq |y-x|\leq 2\diam\, R$, see \rf{dilated boundary cube}.

For every $Q\in U_m(R)$ and every $x\in Q$,
$$
|x-c_R|\leq \diam\, Q + (m-1)\diam\, R + \diam\, R \leq C_\DD (m+1)\diam\, R,
$$
where $C_\DD$ is the constant in \rf{boundary dyadic cube diam control} in the boundary dyadic lattice \cref{lemma:def of boundary dyadic cubes}. Also, if $x\in B_Q\supset B(Q)$ (not necessarily in $Q$), then
$$
|x-c_R|\leq |x-c_Q| + |c_Q-c_R|\leq \diam\, Q + C_\DD (m+1)\diam\, R \leq C_\DD (m+2) \diam\, R.
$$
This implies that $\bigcup_{Q\in U_m (R)} B_Q \supset\bigcup_{Q\in U_m (R)} B(Q)$ and $\bigcup_{Q\in U_m (R)} Q$ are subsets of
\begin{equation}\label{def of ball containing U_m R}
\mathfrak B_R\coloneqq C_\DD (m+2) B_R,
\end{equation}
with $B_R$ as in \rf{ball inside outside the cube}.

Since $\sigma$ is ADR, we will use that
$$
\sigma(3R_W)\approx\sigma(R)\approx_m\sigma(\mathfrak B_R)\approx_m \sigma\left(\bigcup_{Q\in U_m (R)} Q\right).
$$

\begin{rem}\label{rem:candidate of approx plane of cube}
    Under the choice in \rf{choice of plane L in betas large scales}, the estimate \rf{line approx cube} holds with $\varepsilon = 2k\delta_\beta$ and $L_Q = L_{kB_Q}$. Indeed, as $kQ\subset kB_Q$ we have
    $$
    \sup_{x\in k Q} \dist(x,L_Q)
    \leq  \sup_{x\in \partial\Omega\cap kB_Q} \dist(y,L_Q)
    \overset{\text{\rf{choice of plane L in betas large scales}}}{\leq} 2 \delta_\beta k \diam\, Q.
    $$
    Recall in \cref{sec:construction of Lip graph} we needed $\varepsilon=2k\delta_\beta \ll\alpha$, which is granted by \rf{parameter condition 4}.
\end{rem}

We define the following subfamilies of $\Stop(R)$ (see \cref{stopping conditions}) of big and small angle respectively as
$$
\begin{aligned}
\BA(R) &\coloneqq \{Q\in \Stop(R) : \angle(L_R,L_Q) >\alpha\}, \text{ and}\\
\SA(R) &\coloneqq \{Q\in \Stop(R) : \angle(L_R,L_Q) \leq \alpha\},
\end{aligned}
$$
so that $\Stop (R) = \BA(R) \cup \SA(R)$ and the union is disjoint. By the stopping conditions in \cref{stopping conditions} in the construction of the Lipschitz graph, we have
$$
\begin{aligned}
\BA(R) &\supset \{Q\in \Stop(R) : \diam\, Q \geq S\} , \text{ and}\\
\SA(R) &\subset \{Q\in \Stop(R) : \diam\, Q <S\}.
\end{aligned}
$$

\begin{rem}
    If $Q\in \SA(R)$, then $\diam\, Q<S$ but its parent $\widehat Q\in\tree(R)$ satisfies $\diam\, \widehat Q\geq S$ as $\widehat Q\not\in \Stop(R)$, and therefore $\diam\, Q \approx S$.
\end{rem}

Let us define
$$
G_R \coloneqq \bigcup_{Q\in \SA (R)} Q. 
$$
The notation $G_R$ stands for ``good'' in the following sense:

\begin{lemma}\label{lemma:big angle has small measure}
We have $\sigma\left(\bigcup_{Q\in U_m (R)} Q\setminus G_R\right) = \sigma\left(\bigcup_{Q\in\BA(R)} Q\right) \lesssim \delta^{1/3} \sigma(3R_W)$.
	\begin{proof}
        For each $Q\in \BA(R)$, let $L_Q$ be as in \rf{line approx cube} (with $\varepsilon = 2k\delta_\beta$) and $L_{B_Q}$ as in \rf{choice of plane L in betas large scales}. By the same computations in \cref{rem:candidate of approx plane of cube}, both $L_Q$ and $L_{B_Q}$ satisfy the assumptions of \cref{lemma:DS91 lemma 5.13} with $\varepsilon = 2k\delta_\beta$, and hence we conclude that $\angle(L_Q,L_{B_Q})\lesssim k\delta_\beta$. By the same argument, $\angle (L_R,L_{\mathfrak B_R})\lesssim k\delta_\beta$, where $\mathfrak B_R$ is the ball in \rf{def of ball containing U_m R} so that $\bigcup_{Q\in U_m (R)} Q\subset \mathfrak B_R$.
        
        From $\angle(L_Q,L_{B_Q})\lesssim k\delta_\beta$, $\angle (L_R,L_{\mathfrak B_R})\lesssim k\delta_\beta$, and the stopping condition $\angle (L_R,L_Q) >\alpha$, we get $\angle (L_{\mathfrak B_R},L_{B_Q})>\alpha/2$ provided $\delta_\beta \ll \alpha/k$, see \rf{parameter condition 4}. From this and \rf{two normal vectors are almost parallel} we have $\angle (m_{\mathfrak B_R}\nu, m_{B_Q}\nu)>\alpha/3$ provided $\delta^{\frac{1}{2n}}\ll \alpha$, which is granted by \rf{parameter condition 4}. In particular $|m_{\mathfrak B_R} \nu - m_{B_Q}\nu| \gtrsim 1-\cos \alpha$ by the law of cosines.
		
		Therefore, by Chebysheff's inequality we have
		$$
		\begin{aligned}
		\sigma\left(\bigcup_{Q\in\BA(R)} Q\right) 
		&= \sum_{Q\in \BA(R)} \sigma (Q) 
		\lesssim \frac{1}{1-\cos\alpha}\sum_{Q\in \BA(R)} \int_{Q} |m_{B_Q}\nu - m_{\mathfrak B_R}\nu| \, d\sigma\\
		&\leq \frac{1}{1-\cos\alpha} \left(\sum_{Q\in \BA(R)} \int_{Q} |\nu  - m_{B_Q}\nu| \, d\sigma + \sum_{Q\in \BA(R)} \int_Q |\nu - m_{\mathfrak B_R}\nu| \, d\sigma \right) \\
            &\leq \frac{1}{1-\cos\alpha} \left(\sum_{Q\in \BA(R)} \int_{B_Q} |\nu  - m_{B_Q}\nu| \, d\sigma + \int_{\mathfrak B_R} |\nu - m_{\mathfrak B_R}\nu| \, d\sigma \right),
		\end{aligned}
		$$
        where we simply used in the last step that $Q\subset B_Q$ and $\bigcup_{Q\in\BA(R)} Q\subset \bigcup_{Q\in U_m (R)} Q \subset \mathfrak B_R$. The first term is controlled by
        $$
        \begin{aligned}
        \sum_{Q\in \BA(R)} \int_{B_Q} |\nu  - m_{B_Q}\nu| \, d\sigma
        \overset{\text{\rf{bmo norm in large scales}}}&{\leq}  \delta_* \sum_{Q\in \BA(R)} \sigma (B_Q)
        \approx \delta_* \sum_{Q\in \BA(R)} \sigma (Q)\\
        &= \delta_* \sigma\left(\bigcup_{Q\in\BA(R)} Q\right)
        \leq \delta_* \sigma \left(\bigcup_{Q\in U_m (R)} Q\right)
        \lesssim \delta_* \sigma (\mathfrak B_R).
        \end{aligned}
        $$
        The other term is controlled by
        $$
        \int_{\mathfrak B_R} |\nu - m_{\mathfrak B_R}\nu| \, d\sigma
        \overset{\text{\rf{bmo norm in large scales}}}{\leq} \delta_* \sigma(\mathfrak B_R).
        $$
        
        All in all, from the last three inline equations, and using $\delta_*\leq \delta$, $\delta \leq \alpha^{3}$ (see \rf{parameter condition 4}) and $\frac{\delta}{1-\cos(\delta^{1/3})} \leq 3 \delta^{1/3}$ for small enough $\delta$, we have
        $$
        \sigma\left(\bigcup_{Q\in\BA(R)} Q\right) \lesssim \frac{\delta_*}{1-\cos\alpha} \sigma(\mathfrak B_R)
        \lesssim \delta^{1/3} \sigma(\mathfrak B_R).
        $$
        The lemma follows as $\sigma(3R_W)\approx_m\sigma(\mathfrak B_R)$.
	\end{proof}
\end{lemma}

Using the set $G_R$ defined above, we decompose
$$
	\sigma(\{x\in R_W : K_{l,*} (f\characteristic_{3R_W})(x) > 40\lambda, \M_{1+\gamma} f(x) \leq A\lambda\})
    \leq \boxed{1} + \boxed{2} + \boxed{3},
$$
where
$$
\begin{aligned}
    \boxed{1} &\coloneqq \sigma(R_W\setminus G_R),\\
    \boxed{2} &\coloneqq \sigma(\{x\in G_R : K_{l,*}(f\characteristic_{3R_W\cap G_R})(x) > 20\lambda, \M_{1+\gamma} f(x) \leq A\lambda\}),\text{ and}\\
    \boxed{3} &\coloneqq \sigma(\{x\in \partial\Omega : K_{l,*}(f\characteristic_{3R_W\setminus G_R})(x) > 20\lambda\}).
\end{aligned}
$$

Using the control of the measure of $\bigcup_{Q\in U_m (R)} Q\setminus G_R$ in \cref{lemma:big angle has small measure}, we estimate \boxed{1} and \boxed{3}:

\begin{lemma}\label{term bad-bad has small measure}
     We have $\boxed{1} \lesssim \delta^{1/3} \sigma(3R_W)$.
     \begin{proof}
         This is just by the inclusion $\partial\Omega \cap 3R_W\subset mR \subset \bigcup_{Q\in U_m(R)} Q$ and \cref{lemma:big angle has small measure}.
     \end{proof}
\end{lemma}

\begin{lemma}\label{term good-bad has small measure}
	We have $\boxed{3}\lesssim A^{\sqrt{1+\gamma}} \delta^{\frac{1}{3}(1-\frac{1}{\sqrt{1+\gamma}})} \sigma(3R_W)$. In particular $\boxed{3}\leq c_\delta \sigma(3R_W)$ with $c_\delta\to 0$ as $\delta\to 0$, by \rf{parameter condition 3}.
	\begin{proof}
            By Chebysheff's inequality and the fact that $K_{l,*} g \leq K_* g$ for any $g\in L^p(\sigma)$, we have
            $$
            \boxed{3}
            \lesssim \frac{1}{\lambda^{\sqrt{1+\gamma}}} \int_{\partial\Omega} |K_{l,*} (f\characteristic_{3R_W\setminus G_R})|^{\sqrt{1+\gamma}}\, d\sigma
            \leq \frac{1}{\lambda^{\sqrt{1+\gamma}}}\int_{\partial\Omega} |K_{*} (f\characteristic_{3R_W\setminus G_R})|^{\sqrt{1+\gamma}}\, d\sigma.
            $$
            Since $L^p(\sigma)\subset L^{\sqrt{1+\gamma}}_{\loc} (\sigma)$ (as we are assuming that $0<\gamma<p-1$), by \rf{maximal operator of K is bounded} (since $\partial\Omega$ is uniformly rectifiable) and Hölder's inequality respectively, the right-hand side term is bounded by
            $$
            \begin{aligned}
            \frac{1}{\lambda^{\sqrt{1+\gamma}}}\int_{\partial\Omega} |K_{*} (f\characteristic_{3R_W\setminus G_R})|^{\sqrt{1+\gamma}}\, d\sigma
            &\lesssim \frac{1}{\lambda^{\sqrt{1+\gamma}}} \int_{\partial\Omega} |f\characteristic_{3R_W\setminus G_R}|^{\sqrt{1+\gamma}}\, d\sigma\\
            &\lesssim \frac{\sigma(3R_W\setminus G_R)^{1-\frac{1}{\sqrt{1+\gamma}}}}{\lambda^{\sqrt{1+\gamma}}} \left(\int_{3R_W} |f|^{1+\gamma} \, d\sigma\right)^{\frac{1}{\sqrt{1+\gamma}}}\\
            &=\frac{\sigma(3R_W\setminus G_R)^{1-\frac{1}{\sqrt{1+\gamma}}} \sigma(3R_W)^{\frac{1}{\sqrt{1+\gamma}}}}{\lambda^{\sqrt{1+\gamma}}} \left(\avint_{3R_W} |f|^{1+\gamma} \, d\sigma \right)^{\frac{1}{\sqrt{1+\gamma}}}.
            \end{aligned}
            $$
            Using that $\left(\avint_{3R_W} |f|^{1+\gamma} \, d\sigma\right)^{\frac{1}{\sqrt{1+\gamma}}}\lesssim \M_{1+\gamma} f(x_0)^{\sqrt{1+\gamma}}$, $\M_{1+\gamma} f(x_0)\leq A\lambda$ and $3R_W \subset U_m (R)$ in the latter term, we conclude
            $$
            \boxed{3}
            \lesssim A^{\sqrt{1+\gamma}} \left(\frac{\sigma\left(\bigcup_{Q \in U_m (R)} Q \setminus G_R\right)}{\sigma(3R_W)}\right)^{1-\frac{1}{\sqrt{1+\gamma}}} \sigma(3R_W),
            $$
            and the lemma follows by applying \cref{lemma:big angle has small measure} in the last term.
	\end{proof}
\end{lemma}

It remains to study the term $\boxed{2}$. As in the terms \boxed{1} and \boxed{3}, we want to see:
\begin{lemma}\label{term good-good has small measure}
	We have $\boxed{2}\lesssim (\alpha A^2 + \delta^{\frac{1-1/\sqrt{1+\gamma}}{4n}} A^2) \sigma(3R_W)$. In particular $\boxed{2}\leq c_\delta \sigma(3R_W)$ with $c_\delta \to 0$ as $\delta\to 0$, by \rf{parameter condition 3} and \rf{parameter condition 5}.
\end{lemma}

The rest of the section is devoted to the proof of \cref{term good-good has small measure}. By \cref{prop:construction Lip graph} there exists a Lipschitz graph $\Gamma_R$ given by $A:L_R\to L_R^\perp$ with norm $\leq C\alpha$ such that
$$
\dist(x,(\Pi_R(x),A(\Pi_R(x)))) \leq C k \delta_\beta d_R (x) \text{ for all } x\in k_0R\supset mR\supset 3R_W.
$$
In particular, for $x \in Q\in \SA(R)$, if $\widehat Q\supset Q$ is the parent of $Q$ then $\widehat Q\in \tree(R)$, and hence
$$
d_R (x) =\inf_{Q\in \tree(R)} \{\dist(x,Q)+\diam\, Q\}\leq \diam\, \widehat Q\approx \diam\, Q\approx S
$$
and
$$
\dist(x,(\Pi_R(x),A(\Pi_R(x)))) \leq Ck\delta_\beta d_R (x) \lesssim k\delta_\beta S.
$$
So, if $\delta_\beta \ll 1/k$ is small enough we have
$$
\Gamma_R \cap \frac{1}{2} B(Q) \not=\emptyset;
$$
see \rf{ball inside outside the cube} for the definition of $B(Q)$. This in particular implies
$$
\HH^n|_{\Gamma_R} (B(Q))\approx \sigma(Q).
$$

Given $f\characteristic_{3R_W \cap G_R}\in L^p (\sigma)$, let us define an auxiliary function $f_{\Gamma_R}$ in $\Gamma_R$. First, for each $Q\in \SA(R)$, let
\begin{equation}\label{bound of the new function}
f_{\Gamma_R} (Q) \coloneqq \frac{1}{\HH^n (\Gamma_R \cap B(Q))} \int_Q f(y)\characteristic_{3R_W \cap G_R}(y)\, d\sigma(y) \approx \avint_Q f(y)\characteristic_{3R_W \cap G_R} (y)\, d\sigma(y),
\end{equation}
and we define the function $f_{\Gamma_R}$ in $\Gamma_R$ as
\begin{equation}\label{def:def of f_gamma}
	f_{\Gamma_R} (x)
	\coloneqq \sum_{Q\in \SA(R)} f_{\Gamma_R} (Q) \characteristic_{B(Q)} (x), \quad x\in\Gamma_R.
\end{equation}
Note that $f_{\Gamma_R}$ is a piecewise constant function, since $\{B(Q)\}_{Q\in \SA(R)}$ is a pairwise disjoint family. This follows from the fact that $\{Q\}_{Q\in\SA(R)}$ is a pairwise disjoint family and the choice of $c_1$ in \rf{ball inside outside the cube}. By definition of $f_{\Gamma_R}$ and \rf{def of ball containing U_m R}, we have
\begin{equation}\label{support of function in lip graph}
    \supp f_{\Gamma_R} \subset \bigcup_{Q\in \SA(R)} \overline{B(Q)} \subset \mathfrak B_R.
\end{equation}

From the definition of $f_{\Gamma_R}$ from $f\characteristic_{3R_W \cap G_R}$ in \rf{def:def of f_gamma} we have that both $f\characteristic_{3R_W \cap Q} \, d\sigma$ and $f_{\Gamma_R} \characteristic_{B(Q)} \, d\HH^n|_{\Gamma_R}$ have the same mass for each $Q\in \SA(R)$. That is,
\begin{equation}\label{same mass}
	\int_{\Gamma_R} f_{\Gamma_R} \characteristic_{B(Q)} \, d\HH^n
	= \HH^n (\Gamma_R \cap B(Q)) f_{\Gamma_R} (Q)
	= \int_Q f \characteristic_{3R_W\cap Q} \, d\sigma.
\end{equation}

Before seeing some properties of $f_{\Gamma_R}$, let us define the subfamily $I_R$ of $\SA(R)$ as
$$
I_R \coloneqq \{Q\in \SA (R) : \exists x\in Q \text{ with } K_{l,*}(f\characteristic_{3R_W\cap G_R})(x)>6\lambda \text{ and } \M_{1+\gamma} f(x) \leq A\lambda\}.
$$
For each $Q\in I_R$, we fix $z_Q\in Q$ (do not confuse with $c_Q$, the ``center'' of $Q$) such that
\begin{equation}\label{eq:maximal bound for cubes in IR}
\M_{1+\gamma} f(z_Q) \leq A\lambda.
\end{equation}

The following lemma shows that $f_{\Gamma_R}$ inherits a localized maximal bound from $f$, restricted to the cubes $Q\in I_R$.

\begin{lemma}\label{lemma:maximal properties new function}
	For each $Q\in I_R$ we have
	$$
	\left(\avint_{\Gamma_R\cap B(Q)} |f_{\Gamma_R}|^{1+\gamma}\, d\HH^n\right)^{\frac{1}{1+\gamma}}\lesssim A\lambda
	\quad
	\text{and}
	\quad
	\left(\avint_{\Gamma_R\cap \mathfrak B_R} |f_{\Gamma_R}|^{1+\gamma}\, d\HH^n\right)^{\frac{1}{1+\gamma}}\lesssim A\lambda.
	$$
	\begin{proof}
		Let $z_Q$ denote the point in $Q$ such that $\M_{1+\gamma} f(z_Q)\leq A\lambda$, see \rf{eq:maximal bound for cubes in IR}. The first inequality follows from the definition of $f_{\Gamma_R}$, \rf{bound of the new function} and $Q\in I_R$. That is,
		$$
		\left(\avint_{\Gamma_R\cap B(Q)} |f_{\Gamma_R}(y)|^{1+\gamma}\, d\HH^n(y)\right)^{\frac{1}{1+\gamma}}
		=|f_{\Gamma_R} (Q)|
		\lesssim \avint_Q |f(y)| \, d\sigma(y)
		\lesssim \M_{1+\gamma} f(z_Q) \leq A\lambda.
		$$

		Similarly, the second inequality is obtained as
		$$
		\begin{aligned}
		\avint_{\Gamma_R \cap \mathfrak B_R} |f_{\Gamma_R} (y)|^{1+\gamma} \, d\HH^n (y) 
		&\lesssim \frac{1}{\HH^n (\Gamma_R\cap \mathfrak B_R)} \int_{\Gamma_R \cap \mathfrak B_R} \left(\sum_{Q\in \SA(R)} |f_{\Gamma_R} (Q)| \characteristic_{B(Q)}(y) \right)^{1+\gamma} \, d\HH^n (y)\\
		&\approx \frac{1}{\ell(R)^n} \sum_{Q\in \SA(R)} \int_{\Gamma_R \cap B(Q)} |f_{\Gamma_R} (Q)|^{1+\gamma} \, d\HH^n (y) 
		 \\
		\overset{\text{\rf{bound of the new function}}}&{\lesssim} \frac{S^n}{\ell(R)^n} \sum_{Q\in \SA(R)} \left(\avint_Q |f(y)|\, d\sigma(y) \right)^{1+\gamma}\\
		\overset{\text{(Jensen)}}&{\lesssim} \frac{1}{\ell(R)^n} \sum_{Q\in \SA(R)} \int_Q |f(y)|^{1+\gamma} \, d\sigma(y)
            \overset{\text{\rf{def of ball containing U_m R}}}{\lesssim} \avint_{\mathfrak B_R} |f(y)|^{1+\gamma} \, d\sigma(y)\\
		&\lesssim \M_{1+\gamma} f(z_Q)^{1+\gamma} \leq (A\lambda)^{1+\gamma},
		\end{aligned}
		$$
		where we took any $z_Q$ with $Q\in I_R$ in the last inequality.
	\end{proof}
\end{lemma}

To shorten the notation, for $\varepsilon>0$ we denote
$$
\Phi_\varepsilon (z) \coloneqq \frac{z}{|z|^{n+1}} \characteristic_{\{|z|> \varepsilon\}} \text{ for }z\in \R^{n+1}\setminus \{0\}.
$$

We present two lemmas that we will need in the proof of \cref{term good-good has small measure}.

\begin{lemma}\label{difference with approximation is small}
	For all $x\in Q\in I_R$, all $x^\prime \in B(Q)\cap \Gamma_R$ and every $\varepsilon\geq T$,
	$$
	\left|\int_{\partial\Omega} \Phi_\varepsilon (x-y) f(y)\characteristic_{3R_W\cap G_R}(y) \, d\sigma(y) 
	- \int_{\Gamma_R} \Phi_\varepsilon (x^\prime-y) f_{\Gamma_R}(y) \, d\HH^n(y)\right| \leq \frac{1}{2}\lambda.
	$$
	\begin{proof}
		We write 
		$$
		\begin{aligned}
			\boxed{M} =&\, \left|\int \left( \Phi_\varepsilon (x-y) f(y)\characteristic_{3R_W\cap G_R}(y) \, d\sigma(y) - \Phi_\varepsilon(x^\prime-y) f_{\Gamma_R}(y) \, d\HH^n|_{\Gamma_R}(y) \right) \right| \\
			\leq&\, \left|\int_{\partial\Omega} (\Phi_\varepsilon (x-y)-\Phi_\varepsilon (x^\prime-y)) f(y) \characteristic_{3R_W\cap G_R}(y) \, d\sigma(y) \right|\\
			&\, +\left| \int \Phi_\varepsilon (x^\prime-y) f(y)\characteristic_{3R_W\cap G_R}(y) \, d\sigma(y) - \int \Phi_\varepsilon (x^\prime-y)f_{\Gamma_R}(y) \, d\HH^n|_{\Gamma_R}(y) \right|\\
			\eqqcolon&\, \boxed{M1} + \boxed{M2}.
		\end{aligned}
		$$
		
		Let us bound the term $\boxed{M1}$. Since $x\in Q$ and $x^\prime\in B(Q)$ satisfy $|x-x^\prime|\leq 2\diam \, Q < 2S$, and $\Phi_\varepsilon$ is the truncated kernel of big scales $\varepsilon\geq T \gg 100S$, in particular $|x-y|\geq T\gg 100S$, we have $|x-x^\prime|/|x-y|\leq S/T\leq 1/2$. Thus, we can use the cancellation property of the Calderón-Zygmund kernel (or just the mean value theorem) to obtain
		$$
			|\Phi_\varepsilon (x-y)-\Phi_\varepsilon (x^\prime-y)| 
			\lesssim \frac{|x-x^\prime|}{|x-y|^{n+1}} \characteristic_{B(x,\varepsilon/2)^c}(y)
			\lesssim \frac{S}{|x-y|^{n+1}}\characteristic_{B(x,T/2)^c}(y),
		$$
		where we used in the last step that $\varepsilon\geq T$ and that $x\in Q$ and $x^\prime \in B(Q)$ satisfy $|x-x^\prime|\leq 2 S$. From this we get
		$$
		\begin{aligned}
			\boxed{M1} &\leq \int_{\partial\Omega} |\Phi_\varepsilon (x-y)-\Phi_\varepsilon (x^\prime-y)||f(y)| \characteristic_{G_R} (y) \, d\sigma(y)
			\lesssim  S\int_{B(x,T/2)^c} \frac{|f(y)|}{|x-y|^{n+1}} \, d\sigma(y)\\
			&= S \sum_{k=0}^\infty \int_{B(x,2^k T)\setminus B(x, 2^{k-1}T)} \frac{|f(y)|}{|x-y|^{n+1}} \, d\sigma(y)
            \lesssim \frac{S}{T} \sum_{k=0}^\infty \frac{1}{2^k} \avint_{B(x,2^k T)} |f| \, d\sigma(y)
			\overset{\text{\rf{eq:maximal bound for cubes in IR}}}{\lesssim} \frac{S}{T} A\lambda,
		\end{aligned}
		$$
		where we used in the last step that $Q\in I_R$.
		
		It remains to bound the term $\boxed{M2}$. We have
		$$
		\begin{aligned}
			\boxed{M2}=&\left| \int \Phi_\varepsilon (x^\prime-y) f(y)\characteristic_{3R_W\cap G_R}(y) \, d\sigma(y) 
			-\int \Phi_\varepsilon (x^\prime-y) f_{\Gamma_R}(y) \, d\HH^n|_{\Gamma_R} (y) \right|\\
			=&\left|\sum_{\widetilde Q\in \SA(R)} \int_{\widetilde Q} \Phi_\varepsilon (x^\prime-y) f(y)\characteristic_{3R_W}(y) \, d\sigma(y) 
			-\sum_{\widetilde Q\in \SA(R)} \int_{B(\widetilde Q)} \Phi_\varepsilon (x^\prime-y) f_{\Gamma_R}(y) \, d\HH^n|_{\Gamma_R} (y)\right|\\
			\leq&\sum_{\widetilde Q\in \SA(R)} \left|\int \Phi_\varepsilon (x^\prime-y) f(y)\characteristic_{3R_W\cap \widetilde Q} (y) \, d\sigma(y) - \int \Phi_\varepsilon (x^\prime-y) f_{\Gamma_R}(y) \characteristic_{B(\widetilde Q)}(y) \, d\HH^n|_{\Gamma_R} (y) \right|.
		\end{aligned}
		$$
            To shorten the notation, we write $\boxed{M2}\leq \sum_{\widetilde Q\in \SA(R)} \boxed{M2(\widetilde Q)}$ where 
            $$
            \boxed{M2(\widetilde Q)} \coloneqq 
            \left|\int \Phi_\varepsilon (x^\prime-y)\left( f(y)\characteristic_{3R_W\cap \widetilde Q} (y) \, d\sigma(y) - f_{\Gamma_R}(y) \characteristic_{B(\widetilde Q)}(y) \, d\HH^n|_{\Gamma_R} (y) \right) \right|.
            $$
            
            Let us study $\boxed{M2(\widetilde Q)}$ for each $\widetilde Q\in \SA(R)$. Since we have that $f(y)\characteristic_{3R_W\cap \widetilde Q} (y) \, d\sigma(y) - f_{\Gamma_R}(y) \characteristic_{B(\widetilde Q)}(y) \, d\HH^n|_{\Gamma_R} (y)$ has zero mass for each $\widetilde Q\in \SA(R)$, see \rf{same mass}, we can add $\Phi_\varepsilon (x^\prime-c_{\widetilde Q})$ inside the integral sign in the last term, where $c_{\widetilde Q}$ is the ``center'' of $\widetilde Q$,\footnote{Any point in $\widetilde Q$ would do the job.} and therefore
            $$
            \boxed{M2(\widetilde Q)}=
            \left|\int (\Phi_\varepsilon (x^\prime-y)-\Phi_\varepsilon (x^\prime-c_{\widetilde Q}))( f(y)\characteristic_{3R_W \cap \widetilde Q} (y) \, d\sigma(y) - f_{\Gamma_R}(y) \characteristic_{B(\widetilde Q)}(y) \, d\HH^n|_{\Gamma_R} (y) ) \right|.
            $$
            With this we have
            $$
            \begin{aligned}
            \boxed{M2(\widetilde Q)}
            \leq&\, \int |\Phi_\varepsilon (x^\prime-y)-\Phi_\varepsilon (x^\prime-c_{\widetilde Q})| |f(y)|\characteristic_{3R_W \cap \widetilde Q} (y) \, d\sigma(y)\\
            &+ \int |\Phi_\varepsilon (x^\prime-y)-\Phi_\varepsilon (x^\prime-c_{\widetilde Q})||f_{\Gamma_R}(y)| \characteristic_{B(\widetilde Q)}(y) \, d\HH^n|_{\Gamma_R} (y).
            \end{aligned}
            $$
            For $y\in \widetilde Q$ or $y\in B(\widetilde Q)$, note first that $|y-c_{\widetilde Q}|\leq \diam\, \widetilde Q < S$, and second, $|x^\prime-y|> \varepsilon\geq T\gg 100S$ by the truncation of the kernel $\Phi_\varepsilon$. Hence, $|y-c_{\widetilde Q}|/|y-x^\prime|\leq 1/2$. Therefore, by the cancellation of the Calderón-Zygmund kernel $\Phi_\varepsilon$ we have
            $$
            |\Phi_\varepsilon (x^\prime-y)-\Phi_\varepsilon (x^\prime-c_{\widetilde Q})|
            \lesssim \frac{|y-c_{\widetilde Q}|}{|x^\prime-y|^{n+1}}\leq \frac{S}{|x^\prime-y|^{n+1}}.
            $$
            Thus, defining
            $$
            \begin{aligned}
                \boxed{M2(\widetilde Q)_\sigma}
                &\coloneqq \int_{\widetilde Q\cap\{|y-x^\prime|\geq T/2\}} \frac{S}{|x^\prime-y|^{n+1}} |f(y)| \, d\sigma(y),\text{ and}\\
                \boxed{M2(\widetilde Q)_{\HH^n}}
                &\coloneqq \int_{B(\widetilde Q)\cap \Gamma_R \cap\{|y-x^\prime|\geq T/2\}} \frac{S}{|x^\prime-y|^{n+1}} |f_{\Gamma_R}(y)| \, d\HH^n (y).
            \end{aligned}
            $$
            we have
            $$
            \boxed{M2(\widetilde Q)} \lesssim \boxed{M2(\widetilde Q)_\sigma} + \boxed{M2(\widetilde Q)_{\HH^n}}.
            $$

		Summing $\boxed{M2(\widetilde Q)_\sigma}$ over $\widetilde Q \in \SA(R)$, as we did in the bound of $\boxed{M1}$ we have
		$$
		\begin{aligned}
			\sum_{\widetilde Q \in \SA(R)} \boxed{M2(\widetilde Q)_\sigma}
			&\leq S \int_{B(x^\prime, T/2)^c} \frac{|f(y)|}{|x^\prime-y|^{n+1}} \, d\sigma(y)
                \lesssim \frac{S}{T} \sum_{k=0}^\infty \frac{1}{2^k} \avint_{B(x^\prime, 2^k T)} |f| \, d\sigma
                \overset{\text{\rf{eq:maximal bound for cubes in IR}}}{\lesssim} \frac{S}{T} A\lambda,
		\end{aligned}
		$$
		where we used in the last step that $Q\in I_R$.

        We obtain the same bound for $\sum_{\widetilde Q\in \SA(R)} \boxed{M2(\widetilde Q)_{\HH^n}}$. Indeed, by the same computations above we have
        \begin{equation}\label{eq:step 1 of bound in lip graph}
        \sum_{\widetilde Q \in \SA(R)} \boxed{M2(\widetilde Q)_{\HH^n}}
        \lesssim S \sum_{k=0}^\infty \int_{B(x^\prime,T/2)^c} \frac{|f_{\Gamma_R}(y)|}{|x^\prime-y|^{n+1}} \, d\HH^n|_{\Gamma_R}(y) \lesssim \frac{S}{T} \sum_{k=0}^\infty \frac{1}{2^k} \avint_{B(x^\prime,2^k T)} |f_{\Gamma_R}| \, d\HH^n|_{\Gamma_R}.
        \end{equation}
        For each $k\geq 0$, we define
        $$
        \SA_k(R)\coloneqq\{\widetilde Q \in \SA(R) : \widetilde Q \cap B(x^\prime,2^k T)\not=\emptyset\}.
        $$
        Since $|f_{\Gamma_R}(\widetilde Q)|\HH^n|_{\Gamma_R}(B(\widetilde Q))\leq \int_{\widetilde Q} |f|\, d\sigma$ for all $\widetilde Q\in \SA(R)$ (see \rf{bound of the new function}) and $\widetilde Q\subset B(x^\prime,2^{k+1}T)$ if $\widetilde Q\in \SA_k (R)$ (because $\diam\, \widetilde Q<S\ll T$ for all $\widetilde Q\in \SA(R)$), we have
         \begin{equation}\label{eq:step 2 of bound in lip graph}
        \begin{aligned}
        \int_{B(x^\prime,2^k T)} |f_{\Gamma_R}| \, d\HH^n|_{\Gamma_R}
        &\leq \sum_{\widetilde Q\in\SA_k(R)} \int_{B(\widetilde Q)} |f_{\Gamma_R}| \, d\HH^n|_{\Gamma_R}
        = \sum_{\widetilde Q\in\SA_k(R)} |f_{\Gamma_R}(\widetilde Q)| \HH^n|_{\Gamma_R}(B(\widetilde Q))\\
        &\lesssim \sum_{\widetilde Q\in\SA_k(R)}  \int_{\widetilde Q} |f| \, d\sigma
        \leq \int_{B(x^\prime,2^{k+1}T)} |f| \, d\sigma.
        \end{aligned}
        \end{equation}
        By \rf{eq:step 1 of bound in lip graph} and \rf{eq:step 2 of bound in lip graph}, using now that $Q\in I_R$ (i.e., \rf{eq:maximal bound for cubes in IR}), we conclude
        \begin{equation}\label{eq:step 3 of bound in lip graph}
        \sum_{\widetilde Q \in \SA(R)} \boxed{M2(\widetilde Q)_{\HH^n}}
        \lesssim \frac{S}{T} \sum_{k=0}^\infty \frac{1}{2^k} \avint_{B(x^\prime,2^{k+1}T)} |f| \, d\sigma
        \lesssim \frac{S}{T}A\lambda,
        \end{equation}
        as claimed.
    
        All in all, we conclude
        $$
        \boxed{M} \leq \boxed{M1} + \boxed{M2} \leq C \frac{S}{T} A\lambda \leq \frac{1}{2}\lambda,
        $$
        where we use \rf{parameter condition 6} in the last step.
	\end{proof}
\end{lemma}

From now on, we write the double layer potential in $\Gamma_R$ as $K^{\Gamma_R}$, replacing $\partial\Omega$ and its surface measure $\sigma$ by $\Gamma_R$ and $\HH^n|_{\Gamma_R}$ in \rf{eq:truncated epsilon dlp}.

The second lemma for the proof of \cref{term good-good has small measure} is about a continuity-type result for the large scale double layer potential in $\Gamma_R$.

\begin{lemma}\label{one points has big image and all have big image}
	For any $Q\in I_R$ we have $|K_{l,*}^{\Gamma_R} (f_{\Gamma_R}) (x) - K_{l,*}^{\Gamma_R} (f_{\Gamma_R}) (z)| \leq \lambda$ for all $x,z\in B(Q)\cap \Gamma_R$.
	\begin{proof}
		For $\varepsilon\geq T$, since $x,z\in B(Q)\cap \Gamma_R$, $\diam\, Q < S$ and $\varepsilon\geq T\gg S$, we have $|x-y|\approx |z-y|$ if $|y-x|> \varepsilon$. So, using the cancellation of the Calderón-Zygmund kernel, we have
		$$
		\begin{aligned}
			|K_{\varepsilon}^{\Gamma_R} (f_{\Gamma_R})(x)-K_{\varepsilon}^{\Gamma_R} (f_{\Gamma_R})(z)|
			&= \left|\int \langle \Phi_\varepsilon(x-y) - \Phi_\varepsilon(z-y), \nu_{\Gamma_R}(y) \rangle f_{\Gamma_R} \, d\HH^n|_{\Gamma_R} (y)\right|\\
			&\lesssim \int_{|y-x|\geq T/2} \frac{|x-z|}{|x-y|^{n+1}} |f_{\Gamma_R}(y)| \, d\HH^n|_{\Gamma_R} (y)\\
			&\leq S \int_{|y-x|\geq T/2} \frac{1}{|x-y|^{n+1}} |f_{\Gamma_R}(y)| \, d\HH^n|_{\Gamma_R} (y).
		\end{aligned}
		$$
            Note that the last element is in fact $\sum_{\widetilde Q\in \SA(R)} \boxed{M2(\widetilde Q)_{\HH^n}}$ in the proof of \cref{difference with approximation is small}, replacing $x^\prime$ by $x$. Thus, by the same computations in \rf{eq:step 1 of bound in lip graph}, \rf{eq:step 2 of bound in lip graph} and \rf{eq:step 3 of bound in lip graph}, we have $|K_{\varepsilon}^{\Gamma_R} (f_{\Gamma_R})(x)-K_{\varepsilon}^{\Gamma_R} (f_{\Gamma_R})(z)|\lesssim \frac{S}{T}A\lambda$. By the condition \rf{parameter condition 6}, we conclude that this is bounded by $\leq \lambda$.
	\end{proof}
\end{lemma}

For the proof of \cref{term good-good has small measure}, we will compare the double layer potential in $\partial\Omega$ with the double layer potential in the Lipschitz graph $\Gamma_R$, so that we can use \cref{difference with approximation is small,one points has big image and all have big image}, and known properties of $K^{\Gamma_R}$. The following is a simplified version of \cite[Theorem 4.34]{Hofmann2010}, see also \cite[(4.4.9)]{Hofmann2010}.

\begin{theorem}\label{theorem:Lp norm dlp on lip graph}
    Let $\varphi : \R^n \to \R$ be a Lipschitz function, $\Gamma = \{(x,\varphi(x)) : x\in \R^n\}$ its graph, and $1<p<\infty$. Then
    \begin{equation}\label{eq:Lp norm dlp on lip graph}
    \|K_*^{\Gamma}\|_{L^p(\Gamma)\to L^p(\Gamma)} \leq C_{n,p} \|\nabla \varphi\|_{L^\infty (\R^n)}(1+\|\nabla \varphi\|_{L^\infty (\R^n)})^N,
    \end{equation}
    for some $N=N(n)$.
\end{theorem}

Let us see the proof of \cref{term good-good has small measure}. During the proof, $\RR=\RR^{\partial\Omega}$ and $\RR^{\Gamma_R}$ denote the Riesz transform in $\partial\Omega$ and $\Gamma_R$ respectively, see \cref{def:Riesz transform}, and $\RR_*$ and $\RR^{\Gamma_R}_*$ their respective maximal operators as in \rf{eq:CZ maximal operator}.

\begin{proof}[Proof of \cref{term good-good has small measure}]
	For each $Q\in I_R$, let us study first
	$$
	\sigma(\{x\in Q : K_{l,*}(f\characteristic_{3R_W\cap G_R})(x) > 20\lambda, \M_{1+\gamma} f(x) \leq A\lambda\}).
	$$
	For any $x\in Q\in I_R$ and any $z \in B(Q)\cap \Gamma_R$, we have
	$$
		\begin{aligned}
    		K_{l,*} (f\characteristic_{3R_W \cap G_R}) (x)
                \leq &\, K_{l,*}^{\Gamma_R} (f_{\Gamma_R}) (z)
                +\RR_* (f\characteristic_{3R_W \cap G_R} (\nu_{\partial\Omega}(\cdot)-N_{\mathfrak B_R})) (x)\\
                &+\RR_*^{\Gamma_R} (f_{\Gamma_R} (\nu_{\Gamma_R}(\cdot)-N_{\mathfrak B_R})) (z)\\
                &+\sup_{\varepsilon\geq T}\left|\int_{\partial\Omega} \Phi_\varepsilon (x-\cdot) f\characteristic_{3R_W\cap G_R} \, d\sigma 
    			- \int_{\Gamma_R} \Phi_\varepsilon (z-\cdot) f_{\Gamma_R} \, d\HH^n\right|,
		\end{aligned}
	$$
	where $N_{\mathfrak B_R}$ is the orthogonal unit vector to $L_{\mathfrak B_R}$ in \rf{choice of plane L in betas large scales} for the ball $\mathfrak B_R$ defined in \cref{def of ball containing U_m R}. Therefore, if we fix (once and for all, for the cube $Q$) a point $z(Q)\in B(Q)\cap \Gamma_R$ satisfying
	$$
	\RR_*^{\Gamma_R} (f_{\Gamma_R} (\nu_{\Gamma_R}(\cdot)-N_{\mathfrak B_R})) (z(Q))
	\leq \inf_{z\in B(Q)\cap \Gamma_R} \RR_*^{\Gamma_R} (f_{\Gamma_R} (\nu_{\Gamma_R}(\cdot)-N_{\mathfrak B_R})) (z) + \frac{\lambda}{2},
	$$
	then we have that\footnote{The reader may notice that the first set on the right-hand side below is either empty or all of $Q$, since $z(Q)$ does not depend on $x$. In the latter case we will see that it is nonetheless small, as the continuity-type estimate \cref{one points has big image and all have big image} reduces it to the infimum of $K_{l,*}^{\Gamma_R}(f_{\Gamma_R})$ over $B(Q)\cap\Gamma_R$, which is controlled below.}
    \begin{multline*}
    \sigma(\{x\in Q : K_{l,*} (f\characteristic_{3R_W\cap G_R}) (x) >6\lambda\})
    \leq \, \sigma(\{x\in Q : K_{l,*}^{\Gamma_R} (f_{\Gamma_R}) (z(Q)) >2\lambda\})\\
    \begin{aligned}
    &+\sigma(\{x\in Q : \RR_* (f\characteristic_{3R_W\cap G_R} (\nu_{\partial\Omega}(\cdot)-N_{\mathfrak B_R})) (x) >\lambda\})\\
    &+\sigma(\{x\in Q : \inf_{z\in B(Q)\cap \Gamma_R} \RR_*^{\Gamma_R} (f_{\Gamma_R} (\nu_{\Gamma_R}(\cdot)-N_{\mathfrak B_R})) (z) >\lambda\})\\
    &\begin{aligned}
    \,+\, \sigma(\{x\in Q : \sup_{\varepsilon\geq T}\,\biggl|&\int_{\partial\Omega} \Phi_\varepsilon (x-y) f(y)\characteristic_{3R_W \cap G_R}(y) \, d\sigma(y) \\
    &- \int_{\Gamma_R} \Phi_\varepsilon (z(Q)-y) f_{\Gamma_R}(y) \, d\HH^n(y)\biggr| >\lambda\}).
    \end{aligned}
    \end{aligned}
	\end{multline*}
    By \cref{difference with approximation is small}, the last term drops out as the corresponding set is empty, and, since $z(Q)\in B(Q)\cap \Gamma_R$, \cref{one points has big image and all have big image} implies that the first is controlled by $\sigma(\{x\in Q : \inf_{z\in B(Q)\cap \Gamma_R} K_{l,*}^{\Gamma_R} (f_{\Gamma_R}) (z) >\lambda\})$. That is,
	\begin{multline*}
		\sigma(\{x\in Q : K_{l,*} (f\characteristic_{3R_W\cap G_R}) (x) >6\lambda\})
		\leq \sigma(\{x\in Q : \inf_{z\in B(Q)\cap \Gamma_R} K_{l,*}^{\Gamma_R} (f_{\Gamma_R}) (z) >\lambda\})\\
        \begin{aligned}
		&+\sigma(\{x\in Q : \RR_* (f\characteristic_{3R_W\cap G_R} (\nu_{\partial\Omega}(\cdot)-N_{\mathfrak B_R})) (x) >\lambda\})\\
		&+\sigma(\{x\in Q : \inf_{z\in B(Q)\cap \Gamma_R} \RR_*^{\Gamma_R} (f_{\Gamma_R} (\nu_{\Gamma_R}(\cdot)-N_{\mathfrak B_R})) (z) >\lambda\})\\
		\eqqcolon \, & \boxed{2_{Q,i}} + \boxed{2_{Q,ii}} + \boxed{2_{Q,iii}}.
        \end{aligned}
	\end{multline*}
	Using the infimum property, $\boxed{2_{Q,i}}$ is treated as
	$$
	\begin{aligned}
	\boxed{2_{Q,i}}&\leq \frac{\sigma(Q)}{\lambda^{1+\gamma}} \left(\inf_{z\in B(Q)\cap \Gamma_R} K_{l,*}^{\Gamma_R} (f_{\Gamma_R}) (z)\right)^{1+\gamma}
		\approx \frac{\HH^n|_{\Gamma_R}(B(Q))}{\lambda^{1+\gamma}} \left(\inf_{z\in B(Q)\cap \Gamma_R} K_{l,*}^{\Gamma_R} (f_{\Gamma_R}) (z)\right)^{1+\gamma}\\
		 &\leq \frac{1}{\lambda^{1+\gamma}} \int_{B(Q)\cap \Gamma_R} |K_{*}^{\Gamma_R} (f_{\Gamma_R}) (y)|^{1+\gamma} \, d\HH^n (y),
	\end{aligned}
	$$
	in the last step we removed the restriction at large scales. By the same computations above,
	$$
		\boxed{2_{Q,iii}}\lesssim \frac{1}{\lambda^{1+\gamma}} \int_{B(Q)\cap \Gamma_R} |\RR_*^{\Gamma_R} (f_{\Gamma_R} (\nu_{\Gamma_R}(\cdot)-N_{\mathfrak B_R})) (y)|^{1+\gamma} \, d\HH^n (y).
	$$
	Using Chebysheff's inequality, the term $\boxed{2_{Q,ii}}$ is controlled directly by
	$$
		\boxed{2_{Q,ii}}
		\leq \frac{1}{\lambda^{\sqrt{1+\gamma}}} \int_Q |\RR_* (f\characteristic_{3R_W\cap G_R} (\nu_{\partial\Omega}(\cdot)-N_{\mathfrak B_R})) (y)|^{\sqrt{1+\gamma}} \, d\sigma(y).
	$$
	All in all,
	\begin{multline*}
		\sigma(\{x\in Q : K_{l,*} (f\characteristic_{3R_W\cap G_R}) (x) >6\lambda\})
		\lesssim \, \frac{1}{\lambda^{1+\gamma}} \int_{B(Q)\cap \Gamma_R} |K_{*}^{\Gamma_R} (f_{\Gamma_R}) (y)|^{1+\gamma} \, d\HH^n (y) \\
    \begin{aligned}
		&+\frac{1}{\lambda^{\sqrt{1+\gamma}}} \int_Q |\RR_* (f\characteristic_{3R_W\cap G_R} (\nu_{\partial\Omega}(\cdot)-N_{\mathfrak B_R})) (y)|^{\sqrt{1+\gamma}} \, d\sigma(y)\\
		&+\frac{1}{\lambda^{1+\gamma}} \int_{B(Q)\cap\Gamma_R} |\RR_*^{\Gamma_R} (f_{\Gamma_R} (\nu_{\Gamma_R}(\cdot)-N_{\mathfrak B_R})) (y)|^{1+\gamma} \, d\HH^n (y).
	\end{aligned}
    \end{multline*}
	
	Summing over all $Q\in I_R\subset \SA(R)$ and using the above inequality in the last step, we have
	$$
            \begin{aligned}
		\boxed{2}=&\, \sigma(\{x\in G_R : K_{l,*}(f\characteristic_{3R_W\cap G_R})(x) > 20\lambda, \M_{1+\gamma} f(x) \leq A\lambda\})\\
			=& \sum_{Q\in I_R} \sigma(\{x\in Q : K_{l,*}(f\characteristic_{3R_W\cap G_R})(x) > 20\lambda, \M_{1+\gamma} f (x) \leq A\lambda\})\\
			\leq& \sum_{Q\in I_R} \sigma(\{x\in Q : K_{l,*}(f\characteristic_{3R_W\cap G_R})(x) > 6\lambda\})\\
			\lesssim&\frac{1}{\lambda^{1+\gamma}} \int_{\Gamma_R} |K_{*}^{\Gamma_R} (f_{\Gamma_R}) (y)|^{1+\gamma} \, d\HH^n (y) \\
			&+\frac{1}{\lambda^{\sqrt{1+\gamma}}} \int_{\partial\Omega} |\RR_* (f\characteristic_{3R_W\cap G_R} (\nu_{\partial\Omega}(y)-N_{\mathfrak B_R})) (y)|^{\sqrt{1+\gamma}} \, d\sigma(y)\\
			&+\frac{1}{\lambda^{1+\gamma}} \int_{\Gamma_R} |\RR_*^{\Gamma_R} (f_{\Gamma_R} (\nu_{\Gamma_R}(y)-N_{\mathfrak B_R})) (y)|^{1+\gamma} \, d\HH^n (y)\\
			\eqqcolon & \boxed{2_{i}} + \boxed{2_{ii}} + \boxed{2_{iii}}.
		\end{aligned}
	$$
	The smallness of the term $\boxed{2_i}$ will come from the $L^p(\Gamma_R)$ norm of $K^{\Gamma_R}_*$ in the Lipschitz graph with small constant $\lesssim\alpha$, while we will bound the terms $\boxed{2_{ii}}$ and $\boxed{2_{iii}}$ using the smallness assumption on the unit normal vectors, after using that the Riesz transform is bounded in $L^p$ (with constant depending on $p$). We do this in the rest of the proof.

	Using that the Lipschitz norm of $\Gamma_R$ is $\lesssim \alpha$ small, by \cref{theorem:Lp norm dlp on lip graph} we have that the $L^q(\Gamma_R)$ norm (in every $1<q<\infty$) of $K_{*}^{\Gamma_R}$ is $\lesssim_q \alpha$, and hence we get
	$$
	\boxed{2_i} \overset{\text{\rf{eq:Lp norm dlp on lip graph}}}{\lesssim} \frac{\alpha^{1+\gamma}}{\lambda^{1+\gamma}} \int_{\Gamma_R} |f_{\Gamma_R} (y)|^{1+\gamma} \, d\HH^n(y)
	\approx \frac{\alpha^{1+\gamma}}{\lambda^{1+\gamma}} \sigma(\mathfrak B_R) \avint_{\Gamma_R \cap \mathfrak B_R} |f_{\Gamma_R}(y)|^{1+\gamma}\, d\HH^n(y),
	$$
	where we used in the last step that $\supp f_{\Gamma_R}\subset \mathfrak B_R$, see \rf{support of function in lip graph}, and we allow the implicit constant to depend on $\gamma$. By \cref{lemma:maximal properties new function}, the last averaged integral is controlled by $\lesssim (A\lambda)^{1+\gamma}$ and hence
	$$
	\boxed{2_i} \lesssim \alpha^{1+\gamma} A^{1+\gamma} \sigma(\mathfrak B_R) \leq \alpha A^2\sigma(\mathfrak B_R)\approx \alpha A^2\sigma(3R_W),
	$$
	where the second-to-last inequality is simply because $0<\gamma<1$.
	
        Using that the Riesz transform is bounded on UR sets, Hölder's inequality, and that there exists some point $z\in 3R_W\subset\mathfrak B_R$ such that $\M_{1+\gamma} f(z) \leq A\lambda$, we have
        $$
        \begin{aligned}
        \boxed{2_{ii}} &\lesssim \frac{1}{\lambda^{\sqrt{1+\gamma}}}\int_{3R_W} (|f(y)| |\nu_{\partial\Omega}(y)-N_{\mathfrak B_R}|)^{\sqrt{1+\gamma}} \, d\sigma(y) \\
        & \lesssim \frac{1}{\lambda^{\sqrt{1+\gamma}}} \sigma(3R_W) \left(\avint_{\mathfrak B_R} |f(y)| ^{1+\gamma} \, d\sigma(y) \right)^{1/\sqrt{1+\gamma}}
        \left(\avint_{\mathfrak B_R} |\nu_{\partial\Omega}(y)-N_{\mathfrak B_R}| \, d\sigma(y)\right)^{1-1/\sqrt{1+\gamma}} \\
        &\lesssim \frac{1}{\lambda^{\sqrt{1+\gamma}}} \sigma(3R_W) (A\lambda)^{\sqrt{1+\gamma}} \left(\avint_{\mathfrak B_R} |\nu_{\partial\Omega}(y)-N_{\mathfrak B_R}|\, d\sigma(y)\right)^{1-1/\sqrt{1+\gamma}}.
        \end{aligned}
        $$
        In the second step we also used that $|\nu_{\partial\Omega}(\cdot)-N_{\mathfrak B_R}|\leq 2$ $\sigma$-a.e.\ to reduce the exponent of $|\nu_{\partial\Omega}(y)-N_{\mathfrak B_R}|$ to $1$ (as in \rf{eq:reduce high power by bounded oscillation}). Using the smallness assumption of the oscillation of the unit normal vector in \rf{bmo norm in large scales} and that $|m_{\sigma, \mathfrak B_R} \nu_{\partial\Omega}-N_{\mathfrak B_R}|\lesssim \delta^{\frac{1}{4n}}$ by \rf{two normal vectors are almost parallel in distance sense}, we have
        $$
        \avint_{\mathfrak B_R} |\nu_{\partial\Omega}(y)-N_{\mathfrak B_R}|\, d\sigma(y)
        \leq \avint_{\mathfrak B_R} |\nu_{\partial\Omega}(y)-m_{\sigma, \mathfrak B_R} \nu_{\partial\Omega}|\, d\sigma(y)
        +|m_{\sigma, \mathfrak B_R} \nu_{\partial\Omega}-N_{\mathfrak B_R}| \lesssim \delta^{\frac{1}{4n}}.
        $$
        We conclude then that
        $$
        \boxed{2_{ii}} \leq C\delta_*^{\frac{1-1/\sqrt{1+\gamma}}{4n}} A^2 \sigma(3R_W).
        $$
	
	The term $\boxed{2_{iii}}$ is treated similarly, using the $L^\infty$ bound of the oscillation of the unit vector as $\Gamma_R$ is a Lipschitz graph of norm $\lesssim \alpha$. That is, since $\Gamma_R$ is the Lipschitz graph of $A:L_R\to L_R^\perp$ in \cref{prop:construction Lip graph} with norm $\lesssim \alpha$, then $\angle (\nu_{\Gamma_R} (y),N_R)\lesssim \alpha$ for all $y\in \Gamma_R$, where $N_R \perp L_R$, and on the other hand, since $\dist(x,L_R)\leq 2 k\delta_\beta \diam\, R$ for all $x\in kR\supset R$ by \rf{line approx cube} and \cref{rem:candidate of approx plane of cube}, and $\dist(x,L_{\mathfrak B_R})\leq 2\delta_\beta r_{\mathfrak B_R}$ for all $x\in \partial\Omega\cap \mathfrak B_R \supset R$ by \rf{choice of plane L in betas large scales}, by \cref{lemma:DS91 lemma 5.13} we have $\angle (N_{\mathfrak B_R},N_R) = \angle (L_{\mathfrak B_R}, L_R)\lesssim k\delta_\beta\leq \alpha$, provided $\delta_\beta \leq \alpha/k$. So, $\angle (\nu_{\Gamma_R} (y),N_{\mathfrak B_R})\lesssim \alpha$ for all $y\in \Gamma_R$, and by the law of cosines we conclude
    $$
    |\nu_{\Gamma_R} (y)-N_{\mathfrak B_R}|^2 = 2-2|{\cos(\angle (\nu_{\Gamma_R} (y),N_{\mathfrak B_R}))}|
    \leq 2-2\cos(C\alpha)\approx \alpha^2.
    $$
    Using that the Riesz transform is bounded on Lipschitz graphs, $\supp f_{\Gamma_R}\subset \mathfrak B_R$ and $|\nu_{\Gamma_R} (y)-N_{\mathfrak B_R}|\lesssim \alpha$, and \cref{lemma:maximal properties new function}, we have
	$$
	\begin{aligned}
		\boxed{2_{iii}}&\lesssim \frac{1}{\lambda^{1+\gamma}}\int_{\Gamma_R} (|f_{\Gamma_R}(y)||\nu_{\Gamma_R}(y)-N_{\mathfrak B_R}|)^{1+\gamma} \, d\HH^n(y)\\
		&\lesssim \alpha^{1+\gamma} \HH^n|_{\Gamma_R} (\mathfrak B_R) \frac{1}{\lambda^{1+\gamma}}\avint_{\Gamma_R \cap \mathfrak B_R} |f_{\Gamma_R} (y)|^{1+\gamma} \, d\HH^n(y)\\
		&\lesssim \alpha^{1+\gamma} A^{1+\gamma} \sigma(3R_W) \leq \alpha A^2 \sigma(3R_W).
	\end{aligned}
	$$
	
	All in all,
	$$\boxed{2} \leq C\left(\boxed{2_i} + \boxed{2_{ii}} + \boxed{2_{iii}}\right) \leq C(\alpha A^2 + \delta_*^{\frac{1-1/\sqrt{1+\gamma}}{4n}} A^2)\sigma(3R_W),
	$$
	and the lemma follows since $\delta_* \leq \max\{\delta_\beta, \delta_*\} = \delta$.
\end{proof}

With this we have that \cref{term bad-bad has small measure,term good-bad has small measure,term good-good has small measure} are proved. Recall that this implies the good lambda inequality \rf{good lambda 3}, a rewrite of the good lambda inequality \rf{good lambda 2}. As we already saw below \rf{good lambda 2}, it implies the good lambda inequality \rf{good lambda 1}, and in particular \rf{eq:exact Lp norm of large scales dlp} and \rf{eq:Lp norm of large scales dlp}. This concludes the proof of \cref{large scale bounded by bmo norm}.\qed

\section{Uniqueness of the solution of the Neumann problem}\label{sec:uniqueness of neumann problem}

In this section we prove the uniqueness (modulo constant) for the Neumann problem in terms of weak derivatives in unbounded $1$-sided CAD with unbounded boundary for the range of $2n/(n+1)\leq p < n+\varepsilon_1$.

\begin{propo}\label{lemma:uniqueness neumann weak bvp}
    Let $\Omega\subset \R^{n+1}$ be an unbounded $1$-sided CAD with unbounded boundary. Then there exists $\varepsilon_1 = \varepsilon_1 (n,\text{CAD})>0$ such that if $u:\Omega\to\R$ satisfies $\NN(\nabla u)\in L^p(\sigma)$ with $p\in [2n/(n+1),n+\varepsilon_1)$ and 
    $$
    \int_\Omega \nabla u (z) \nabla \varphi(z)\, dm(z) = 0, \text{ for all }\varphi \in C^\infty_c (\R^{n+1}),
    $$
    then $u$ is constant.
\end{propo}

In what follows, the nontangential maximal function is taken with its aperture fixed to $\alpha=1$, as agreed in \rf{eq:Lp of nontangential maximal function does not depend on aperture}.

Before the proof of \cref{lemma:uniqueness neumann weak bvp}, let us state some fundamental properties that will be useful in this section. By \cite[Proposition 3.24]{Hofmann2010}, if $\Omega$ is an ADR domain which is either bounded or has an unbounded boundary and $p\in (0,\infty)$, then
\begin{equation}\label{extra regularity from boundary to interior}
\|u\|_{L^{p(n+1)/n}(\Omega)} \lesssim \|\NN u\|_{L^p(\sigma)},
\end{equation}
where the involved constant depends on the ADR parameter, independent of the function $u:\Omega\to \R$. If $\Omega$ is a $1$-sided CAD, by \cite[(2.609)]{Marin-Martell-Mitrea-x3}, there exists a sufficiently large $C>1$ such that if $\NN^\varepsilon (\nabla u)\in L^p_{\loc}(\sigma)$ for some $p\in (0,\infty]$ and some $\varepsilon>0$, then
\begin{equation}\label{eq:nontan of gradient implies nontan of function}
\NN^{\varepsilon/C} u \in L^p_{\loc} (\sigma).
\end{equation}
In particular, for $p\in (0,\infty)$ and a $1$-sided CAD $\Omega$ that is either bounded or has unbounded boundary, we claim
\begin{equation}\label{eq:u is in W more than p local}
\NN(\nabla u) \in L^p(\sigma) \implies u\in W^{1,p(n+1)/n} (\Omega\cap B)\text{ for any ball } B\text{ centered at }\partial\Omega.
\end{equation}
Indeed, directly from \rf{extra regularity from boundary to interior} we obtain $\nabla u\in L^{p(n+1)/n}(\Omega)$. For a ball $B$ centered at $\partial\Omega$, again by \rf{extra regularity from boundary to interior} we have $\|u\|_{L^{p(n+1)/n}(\Omega\cap B)} \lesssim \|\NN^{r_B} (u\characteristic_B)\|_{L^p(\sigma)}$. We note that $\supp (\NN^{r_B} (u\characteristic_B))\subset 3B$, whence we get $\|\NN^{r_B} (u\characteristic_B)\|_{L^p(\sigma)} \leq \|\NN^{r_B} u\|_{L^p(\sigma|_{3B})}$, which is bounded using that we can take any $\varepsilon>0$ in \rf{eq:nontan of gradient implies nontan of function}.

For the type of domains in \cref{lemma:uniqueness neumann weak bvp}, the following lemma relates the nontangential and the weak derivatives.

\begin{lemma}\label{nontangential derivative and weak derivative}
    Let $\Omega\subset \R^{n+1}$ be an ADR domain with the $2$-sided corkscrew condition, and let $u :\Omega\to \R$ be harmonic. If $\NN (\nabla u) \in L^p (\sigma)$ for some $p\in (1,\infty)$ and the pointwise nontangential limit $(\nabla u)|_{\partial\Omega}$ exists $\sigma$-a.e.\ on $\partial\Omega$, then
    $$
    \int_\Omega \nabla u(z) \nabla \varphi(z) \, dm(z) =\int_{\partial\Omega} \varphi(z) \partial_\nu u (z)\, d\sigma(z), \text{ for all }\varphi\in C_c^\infty (\R^{n+1}).
    $$
    \begin{proof}
        If $\Omega$ is bounded or has unbounded boundary, we define $D=\Omega$, and if $\Omega$ is unbounded and has bounded boundary, then we take a ball $B$ so that $\partial\Omega\subset B/2$ and we define $D=\Omega\cap B$. Note that $\partial D$ is ADR, and we denote $\sigma_D = \HH^n|_{\partial D}$. During the proof $\NN_\Omega$ and $\NN_D$ denote the nontangential maximal operators in $\Omega$ and $D$ respectively. We claim
        \begin{equation}\label{nontang in D is in Lp}
        \NN_D (\nabla u) \in L^p(\sigma_D).
        \end{equation}
        This is by assumption when $D=\Omega$. If $\Omega$ is unbounded with bounded boundary, we have $\NN_D (\nabla u) \leq \NN_\Omega (\nabla u)$ on $\partial\Omega$, and since $u$ is harmonic in $\Omega$ (in particular $u\in C^\infty (\Omega)$) and for $x\in \partial B$ we have that $\Gamma_{D} (x)$ is uniformly far from $\partial\Omega$ (depending on $r_B$), we therefore get $\NN_D (\nabla u) < +\infty$ on $\partial B$. This gives \rf{nontang in D is in Lp}.
        
        Fixing $\varphi\in C^\infty_c (\R^{n+1})$ and assuming also $B\supset \supp\varphi$, this is an almost direct consequence of the divergence theorem in \cite[Theorem 2.8]{Hofmann2010}, applied to the domain $D$ and the vector field $\varphi\nabla u$. Let us recall it in the form we use. Let $U\subset \R^{n+1}$ be an ADR domain with the $2$-sided corkscrew condition, write $\sigma_U=\HH^n|_{\partial U}$, and let $\nu$ be the outward unit normal of $U$. Let $v\in C^0 (U)$ be a vector field with $\divv v\in L^1 (U)$, $\NN v\in L^p (\sigma_U)$ for some $p\in (1,\infty)$, and whose pointwise nontangential trace $v|_{\partial U}$ exists $\sigma_U$-a.e.\ on $\partial U$. If either $U$ is bounded or $v$ has bounded support in $U$ (the latter case being noted in \cite[lines 1-6 in proof of Theorem 3.25]{Hofmann2010}), then Green's formula
\begin{equation}\label{green formula in divergence theorem}
\int_U \divv v \, dm = \int_{\partial U} \langle \nu, v|_{\partial U}\rangle \, d\sigma_U
\end{equation}
holds.

        We apply this to $D$ and $v=\varphi\nabla u$, whose support is bounded since $B\supset \supp\varphi$. Let us check these hypotheses. First, $\partial D$ is ADR and $D$ satisfies the $2$-sided corkscrew condition, inherited from $\Omega$. Next, $\varphi\nabla u\in C^0 (D)$ holds because $u$ is harmonic (and hence smooth) in $\Omega$. For the nontangential maximal function, $\NN_D (\varphi \nabla u)\in L^p(\sigma)$ follows from the bound $\NN_D (\varphi \nabla u)\leq C_\varphi \NN_D(\nabla u)$ and the fact that $\NN_D (\nabla u)\in L^p (\sigma_D)$ by \rf{nontang in D is in Lp}. The pointwise nontangential limit $(\varphi\nabla u)|_{\partial D}$ exists $\sigma_D$-a.e., since $(\nabla u)|_{\partial\Omega}$ exists $\sigma$-a.e.\ by assumption. Finally, $\divv (\varphi\nabla u)\in L^1(D)$ because
        $$
        \int_{D} |{\divv (\varphi\nabla u)}|\, dm
        = \int_{D} |\nabla u \nabla \varphi + \varphi \Delta u|\, dm
        \leq C_\varphi \int_{D\cap B} |\nabla u|\, dm <\infty,
        $$
        where we used $\Delta u =0$ in $D$ and that $\nabla u \in L^{p(n+1)/n}(D)$ by \rf{nontang in D is in Lp} and \rf{extra regularity from boundary to interior}.

        Having checked these conditions, as $\Delta u = 0$ in $\Omega$ and by Green's formula \rf{green formula in divergence theorem} for the domain $D$ and the vector field $\varphi\nabla u$ (with bounded support in $D$), we have
        $$
        \int_{D} \nabla u \nabla \varphi \, dm = \int_{D} \divv (\varphi \nabla u) \, dm = \int_{\partial\Omega} \varphi \partial_\nu u \, d\sigma,
        $$
        as claimed.
    \end{proof}
\end{lemma}

As a consequence of \cref{lemma:uniqueness neumann weak bvp,nontangential derivative and weak derivative}, we get the uniqueness of the Neumann problem with zero nontangential derivative.

\begin{coro}\label{lemma:uniqueness neumann bvp}
    Let $\Omega\subset \R^{n+1}$ and $\varepsilon_1>0$ be as in \cref{lemma:uniqueness neumann weak bvp}. If $p\in [2n/(n+1),n+\varepsilon_1)$, then constant functions in $\Omega$ are the unique solutions of the Neumann boundary value problem
    $$
    \begin{cases*}
    	\Delta u = 0 \text{ in } \Omega,\\
        \NN(\nabla u) \in L^p(\sigma),\\
    	\partial_{\nu} u = 0, \text{ }\sigma\text{-a.e.\ on }\partial\Omega.
    \end{cases*}
    $$
\end{coro}

The rest of the section is devoted to the proof of \cref{lemma:uniqueness neumann weak bvp}. The proof will use the extra regularity in \rf{eq:u is in W more than p local} for functions satisfying $\NN(\nabla u)\in L^p(\sigma)$, the Hölder regularity (up to the boundary) of $W^{1,2}$-solutions of the Neumann problem in terms of weak derivatives (see \cref{thm:holder cont of for zero neumann} below), and the well-known Poincaré inequality for uniform domains (see \cref{lemma:poincare ineq} below).

\begin{theorem}[{\cite[Corollary 4.49]{HofmannSparrius-2025}}]\label{thm:holder cont of for zero neumann}
    Let $\Omega\subset \R^{n+1}$ be as in \cref{lemma:uniqueness neumann weak bvp} and let $B\coloneqq B(\xi,R)$ be a ball of radius $R>0$ centered at $\xi\in\partial\Omega$. There exist $K = K(\text{CAD},n)\geq 1$, $C=C(\text{CAD},n)<\infty$ and $\alpha=\alpha(\text{CAD},n)>0$ such that if $u\in W^{1,2}(\Omega\cap B_{KR}(\xi))$ satisfies
    $$
    \int_\Omega \nabla u (z) \nabla \varphi (z) \, dm(z) = 0, \text{ for all }\varphi\in W_c^{1,2}(B_{KR}(\xi)),
    $$
    then for all $0<r<R/2$ there holds
    $$
    |u(x)-u(y)|\leq C \left(\frac{|x-y|}{r}\right)^\alpha \left(\avint_{\Omega\cap B_{2r}(\xi)} u(z)^2 \, dm(z)\right)^{1/2} \text{ for all }x,y\in \Omega\cap B_r(\xi).
    $$
\end{theorem}

\begin{lemma}[Poincaré inequality]\label{lemma:poincare ineq}
    Let $\Omega\subset \R^{n+1}$ be a uniform domain and let $B$ be a ball of radius $R>0$ centered at $\partial\Omega$. There exists $K>2$, depending only on the uniformity constants of $\Omega$, such that if $u\in W^{1,p}(KB\cap\Omega)$, for $1<p<n+1$, then with $p+\varepsilon_p = p(n+1)/(n+1-p)$ there holds
    \begin{equation}\label{eq:poincare ineq}
        \left(\avint_{\Omega\cap B} |u(z)-m_{\Omega\cap B} u|^{p+\varepsilon_p}\, dm(z) \right)^{1/(p+\varepsilon_p)} \lesssim R \left(\avint_{\Omega\cap KB} |\nabla u(z)|^p\, dm(z)\right)^{1/p}.
    \end{equation}
    \begin{proof}
        See \cite[Theorem 4.1]{Mourgoglou-Tolsa-2024-solvabilityneumannproblemelliptic-arxiv} for instance. A quick inspection of \cite[p.~25, l.~5]{Mourgoglou-Tolsa-2024-solvabilityneumannproblemelliptic-arxiv} reveals that $p+\varepsilon_p = p(n+1)/(n+1-p)$.
    \end{proof}
\end{lemma}

We now turn to the proof of \cref{lemma:uniqueness neumann weak bvp}.

\begin{proof}[Proof of \cref{lemma:uniqueness neumann weak bvp}]
    Let us fix a constant $K>2$ for which both \cref{thm:holder cont of for zero neumann} and the Poincaré inequality in \cref{lemma:poincare ineq} hold, and we allow all subsequent involved constants to depend on $K$. Fixed $\xi\in\partial\Omega$, we write $B_R\coloneqq B_R (\xi)$ for $R>0$. By \rf{eq:u is in W more than p local} we have $u\in W^{1,p(n+1)/n}(\Omega\cap KB_{10R})\subset W^{1,p}(\Omega\cap KB_{10R})$, in particular $u\in W^{1,2}(\Omega\cap KB_{10R})$ since we are assuming $p\geq 2n/(n+1)$. By assumption (and a density argument since $u\in W^{1,2}(\Omega\cap KB_{10R})$), it follows that we can apply \cref{thm:holder cont of for zero neumann} and therefore, for any $x,y\in B(\xi,R/2)$ we have
    \begin{equation}\label{eq:proof uniqueness neumann step 1}
    |u(x)-u(y)| \lesssim \left(\frac{|x-y|}{R}\right)^\alpha \left(\avint_{\Omega\cap B_{R}} |u-m_{\Omega\cap B_R} u|^2 \, dm\right)^{1/2},
    \end{equation}
    with $\alpha = \alpha (n,\text{CAD})>0$.

    If $p\geq 2$, then by the Poincaré inequality \rf{eq:poincare ineq} and Hölder's inequality we have
    $$
    \left(\avint_{\Omega\cap B_{R}} |u-m_{\Omega\cap B_R} u|^2 \, dm\right)^{1/2}
    \lesssim R \left(\avint_{\Omega\cap B_{10R K}} |\nabla u|^2\, dm \right)^{1/2}
    \leq R\left(\avint_{\Omega\cap B_{10R K}} |\nabla u|^p\, dm \right)^{1/p}.
    $$
    If $2n/(n+1)\leq p<2$, then by Hölder's inequality and the Poincaré inequality \rf{eq:poincare ineq} (with $p+\varepsilon_p \coloneqq (n+1)p/(n+1-p)$ which is $\geq 2$ since $p\geq 2n/(n+1)\geq 2(n+1)/(n+3)$),
    $$
    \begin{aligned}
    \left(\avint_{\Omega\cap B_{R}} |u-m_{\Omega\cap B_R(\xi)} u|^2 \, dm\right)^{1/2}
    &\leq \left(\avint_{\Omega\cap B_{R}} |u-m_{\Omega\cap B_R} u|^{p+\varepsilon_p} \, dm\right)^{1/(p+\varepsilon_p)}\\
    &\lesssim R \left(\avint_{\Omega\cap B_{10R K}} |\nabla u|^p\, dm \right)^{1/p}.
    \end{aligned}
    $$
    That is, in any case we obtain the same bound. From this, Hölder's inequality and \rf{extra regularity from boundary to interior} we have
    \begin{equation}\label{eq:proof uniqueness neumann step 2}
    \begin{aligned}
        \left(\avint_{\Omega\cap B_{R}} |u-m_{\Omega\cap B_R} u|^2 \, dm\right)^{1/2}
        &\lesssim R \left(\avint_{\Omega\cap B_{10R K}} |\nabla u|^p\, dm \right)^{1/p}\\
        &\lesssim \frac{R}{R^{(n+1)\frac{n}{p(n+1)}}} \|\nabla u\|_{L^{p(n+1)/n} (\Omega)}
        \lesssim R^{1-n/p} \|\NN(\nabla u)\|_{L^p (\sigma)}.
    \end{aligned}
    \end{equation}

    From \rf{eq:proof uniqueness neumann step 1} and \rf{eq:proof uniqueness neumann step 2} we get
    $$
    |u(x)-u(y)| \lesssim |x-y|^\alpha R^{1-n/p-\alpha}\|\NN(\nabla u)\|_{L^p (\sigma)}.
    $$
    Having fixed $x,y\in \Omega$ and $\xi\in \partial\Omega$, this bound holds as long as $R>0$ is big enough so that $x,y\in B(\xi,R/2)$. Taking $R\to \infty$, we obtain that $u$ is constant in $\Omega$ provided $p<n/(1-\alpha)$, i.e., $p<n+\varepsilon_1$ with $\varepsilon_1 = \varepsilon_1 (\alpha,n)=\varepsilon_1 (\text{CAD},n)>0$.
\end{proof}

\section{Injectivity of \texorpdfstring{$\pm\frac{1}{2}Id + K^*$ in $L^p (\sigma)$ for $p\in [2n/(n+1),n+\varepsilon)$}{±(1/2)Id + K* in Lp(sigma) for p in [2n/(n+1), n+epsilon)}}\label{sec:+-id+k* are injective}

Recall that for \cref{thm:lp dirichlet problem,thm:lp neumann problem}, we are interested in the invertibility of the operators $\frac{1}{2}Id + K$ in $L^{p^\prime}(\sigma)$ and $-\frac{1}{2}Id + K^*$ in $L^p(\sigma)$. By \cref{reduction:all equivalent}, it suffices to show that $\frac{1}{2}Id + K^*$ and $-\frac{1}{2}Id + K^*$ are injective in $L^p(\sigma)$.

More precisely, in this section we see the injectivity of $\pm\frac{1}{2}Id + K^*$ in $L^p (\sigma)$ with $2n/(n+1)\leq p<n+\varepsilon_3$, for $2$-sided CAD's with unbounded boundary. This is stated in \cref{coro:+id + K* is injective}, and follows from the injectivity of $-\frac{1}{2}Id + K^*$ for unbounded $1$-sided CAD with unbounded boundary.

\begin{propo}\label{prop -id+adjoint k is injective}
    Let $\Omega\subset \R^{n+1}$ and $\varepsilon_1>0$ be as in \cref{lemma:uniqueness neumann weak bvp}. Then there exists $0<\varepsilon_2 = \varepsilon_2(n,\text{CAD})\leq \varepsilon_1$ such that $-\frac{1}{2}Id + K^*$ is injective in $L^p (\sigma)$ for all $p\in [2n/(n+1),n+\varepsilon_2)$.
\end{propo}

\Cref{prop -id+adjoint k is injective} is proved in \cref{sec:proof of +-id+k* are injective}. Applying this result to $\overline\Omega^c$ we deduce the following corollary.

\begin{coro}\label{coro:+id + K* is injective}
    Assume that $\R^{n+1}\setminus\overline\Omega$ is an unbounded $1$-sided CAD with unbounded boundary, and let $\varepsilon_2>0$ be as in \cref{prop -id+adjoint k is injective} for the domain $\R^{n+1}\setminus\overline\Omega$. Then $\frac{1}{2}Id + K^*$ is injective in $L^p(\sigma)$ for all $p\in [2n/(n+1),n+\varepsilon_2)$.
    \begin{proof}
        A quick computation reveals that $K_\Omega^* = - K_{\overline{\Omega}^c}^*$. Hence, $\frac{1}{2}Id + K_\Omega^* = \frac{1}{2}Id - K_{\overline{\Omega}^c}^*$, which is injective in $L^p (\sigma)$ for all $p\in [2n/(n+1),n+\varepsilon_2)$ by \cref{prop -id+adjoint k is injective}.
    \end{proof}
\end{coro}

\begin{coro}
    If $\Omega\subset \R^{n+1}$ is an unbounded $2$-sided CAD with unbounded boundary, then there exists $\varepsilon_3 = \varepsilon_3 (n,\text{CAD})>0$ such that both $\pm\frac{1}{2}Id + K^*$ are injective in $L^p (\sigma)$ for all $p\in [2n/(n+1),n+\varepsilon_3)$.
    \begin{proof}
        Take $\varepsilon_3 \coloneqq \min\{\varepsilon_2,\widetilde\varepsilon_2\}$, with $\varepsilon_2$ and $\widetilde\varepsilon_2$ as in \cref{prop -id+adjoint k is injective,coro:+id + K* is injective} for the domains $\Omega$ and $\R^{n+1}\setminus \overline\Omega$ respectively.
    \end{proof}
\end{coro}

By \cref{reduction:all equivalent}, this immediately implies that for $p\in [2n/(n+1),n+\varepsilon_3)$, the operators $\pm\frac{1}{2}Id + K$ and $\pm\frac{1}{2}Id + K^*$ are invertible in $L^{p^\prime}(\sigma)$ and $L^p(\sigma)$ respectively, under small enough flatness assumption of the $\delta$-$(s,S;R)$ domain $\Omega$, see \cref{def:final assumptions of domain}.

\begin{rem}
    For bounded $2$-sided corkscrew ADR domains, the injectivity of $\frac{1}{2}Id + K^*$ in $L^p(\sigma)$ for $p\in[2n/(n+1),\infty)$ is proved in \cite[Theorem 1.7.2 (4)]{MitreaX3-GeometricHarmonicAnalysis-Vol4}. Since the domain is bounded, the injectivity for the full range $[2n/(n+1),\infty)$ follows from the particular case $p=2n/(n+1)$ and Hölder's inequality. The case $p=2$ was already proved in \cite[Proposition 5.11]{Hofmann2010}. Both proofs are based on the divergence theorem. In fact, for the so-called regular SKT domains (see \cite[Definition 4.8]{Hofmann2010}), the injectivity in $L^2(\sigma)$ of $\frac{1}{2}Id + K^*$ is enough to show that all four operators $(\pm)\frac{1}{2}Id + K^{(*)}$ are invertible in $L^p(\sigma)$, for all $1<p<\infty$, see \cite[Proposition 5.12]{Hofmann2010}.
\end{rem}

\subsection{The single layer potential}\label{sec:single layer potential}

Let $\Omega\subset \R^{n+1}$ be an ADR domain with the $2$-sided corkscrew condition (in particular $\partial\Omega$ is UR by \cref{thm:2-sided corkscrew + ADR implies UR boundary}) and $f\in L^1 \left(\frac{d\sigma(x)}{1+|x|^n}\right)$. The modified single layer potential operator associated with $\Omega$ is
\begin{equation}\label{single layer potential}
\SSSSmod f (x) = \SSSS_{{\rm mod},\Omega} f (x) \coloneqq \frac{1}{w_n (1-n)} \int_{\partial\Omega} \left(\frac{1}{|x-y|^{n-1}}-\frac{\characteristic_{B_1(0)^c}(y)}{|y|^{n-1}}\right) f(y)\, d\sigma(y), \!\!\! \quad x\in \R^{n+1}\setminus \partial\Omega,
\end{equation}
and its boundary version is defined as
$$
\Smod f (x) = S_{{\rm mod},\Omega} f (x) \coloneqq \frac{1}{w_n (1-n)} \int_{\partial\Omega} \left(\frac{1}{|x-y|^{n-1}}-\frac{\characteristic_{B_1(0)^c}(y)}{|y|^{n-1}}\right) f(y) \, d\sigma(y), \quad x\in\partial\Omega.
$$

\begin{rem}
    The classical single layer potential $\SSSS$ as well as its boundary version $S$ are defined as $\SSSSmod$ and $\Smod$ but replacing their kernel by $\frac{1}{|x-y|^{n-1}}$. However, this only makes sense for domains with bounded boundary, and in this case, the difference between the classical and the modified single layer potential is constant.
\end{rem}

The single layer potential satisfies
$$
\nabla \SSSSmod f (x) = \frac{1}{w_n} \int_{\partial\Omega} \frac{x-y}{|x-y|^{n+1}} f(y)\, d\sigma(y), \quad x\in \R^{n+1}\setminus \partial\Omega,
$$
see \cite[(3.39)]{Marin-Martell-Mitrea-x3}. Note that this is (modulo a multiplicative constant) the Riesz transform $\RR_\sigma f$. Moreover, it satisfies $\Delta (\SSSSmod f)=0$ in $\R^{n+1}\setminus\partial\Omega$, see \cite[(3.40)]{Marin-Martell-Mitrea-x3}, and $\SSSSmod f$ is continuous through the boundary $\partial\Omega$ in the nontangential sense
\begin{equation}\label{eq:single layer is continuous up to the boundary}
(\SSSSmod f)|_{\partial\Omega}^{\rm nt}(z) = \Smod f(x),\text{ for }\sigma\text{-a.e.\ }x\in\partial\Omega,
\end{equation}
see \cite[(3.47)]{Marin-Martell-Mitrea-x3}.

Recall the interior and exterior normal nontangential derivatives in \rf{def:interior nontang derivative} and \rf{def:exterior nontang derivative}.
There is a jump formula for the interior derivative, see \cite[(3.67)]{Marin-Martell-Mitrea-x3}, which is
\begin{equation}\label{jump formula gradient single layer interior}
\partial_{\nu_\Omega}^{\rm int} \SSSSmod f(x)
= \left(-\frac{1}{2}Id +K^* \right) f(x), \text{ for }\sigma\text{-a.e.\ } x\in\partial\Omega.
\end{equation}
A quick computation reveals that $K_{\overline\Omega^c}^* = - K_\Omega^*$, and using this in the last equality, for $\sigma$-a.e.\ $x\in \partial\Omega$ we have the jump formula for the exterior derivative:
\begin{equation}\label{jump formula gradient single layer exterior}
\begin{aligned}
\partial_{\nu_\Omega}^{\rm ext} \SSSSmod f(x)
&= \lim_{\substack{z\to x\\z\in\Gamma^{\overline\Omega^c} (x)}} \langle \nu (x), \nabla \SSSSmod f(z) \rangle 
= -\lim_{\substack{z\to x\\z\in\Gamma^{\overline\Omega^c} (x)}} \langle \nu_{\overline\Omega^c} (x), \nabla \SSSS_{{\rm mod},\overline\Omega^c} f(z) \rangle \\
\overset{\text{\rf{jump formula gradient single layer interior}}}&{=} -\left(-\frac{1}{2}Id +K_{\overline\Omega^c}^* \right) f(x)
= \left(\frac{1}{2}Id + K^*\right) f(x).
\end{aligned}
\end{equation}

Given any $\varepsilon>0$, by \cite[(3.41)]{Marin-Martell-Mitrea-x3} we have
\begin{equation}\label{eq:nontang eps of Smod f is in intersection of Lp}
\NN^\varepsilon (\SSSSmod f) \in \bigcap_{0<p<\frac{n}{n-1}} L^p_{\loc}(\sigma).
\end{equation}
Arguing as in the proof of \rf{eq:u is in W more than p local} (using \rf{eq:nontang eps of Smod f is in intersection of Lp} with $p=1$ instead of \rf{eq:nontan of gradient implies nontan of function}) we obtain $\SSSSmod f \in L^{(n+1)/n} (\Omega\cap B)$ for any ball $B$ centered at $\partial\Omega$. Note that the same holds for $\overline{\Omega}^c$ (as it satisfies the same assumptions as $\Omega$), and therefore we have
\begin{equation}\label{single layer is in L 1+}
\SSSSmod f \in L^{(n+1)/n}_{\loc} (\R^{n+1}).
\end{equation}

For $f\in L^p(\sigma)$ with $1<p<\infty$, the single layer potential satisfies
\begin{equation}\label{lp norm of nontangential of gradient single layer potential}
\|\NN (\nabla \SSSSmod f)\|_{L^p(\sigma)} \lesssim_{p,\text{UR},\alpha} \|f\|_{L^p(\sigma)},
\end{equation}
see \cite[(3.127)]{Marin-Martell-Mitrea-x3}. By \rf{lp norm of nontangential of gradient single layer potential} (and arguing as in \rf{nontang in D is in Lp} if $\Omega$ is unbounded with bounded boundary) and \rf{extra regularity from boundary to interior}, the single layer potential satisfies
\begin{equation}\label{single layer potential is loc sobolev with extra exponent just omega}
    \nabla (\SSSSmod f) \in L^{p(n+1)/n} (\Omega\cap B) \text{ for any ball } B \text{ centered at }\partial\Omega.
\end{equation}
In fact, if $\Omega$ is bounded or has unbounded boundary, a direct application of \rf{lp norm of nontangential of gradient single layer potential} and \rf{extra regularity from boundary to interior} gives $\nabla (\SSSSmod f) \in L^{p(n+1)/n} (\Omega)$. As before, the same holds for $\overline{\Omega}^c$. Arguing as in the proof of \cref{nontangential derivative and weak derivative}\footnote{By the divergence theorem in \cite[Theorem 2.8]{Hofmann2010} using now \rf{eq:single layer is continuous up to the boundary}, \rf{eq:nontang eps of Smod f is in intersection of Lp}, \rf{single layer is in L 1+} and \rf{single layer potential is loc sobolev with extra exponent just omega}.}, weak derivatives of the single layer potential in $\R^{n+1}$ exist and its gradient is $\characteristic_{\R^{n+1}\setminus \partial\Omega}\nabla \SSSSmod f$, whence we get from \rf{single layer potential is loc sobolev with extra exponent just omega} (with both domains $\Omega$ and $\overline{\Omega}^c$) that in the weak sense there holds
\begin{equation}\label{single layer potential is loc sobolev with extra exponent}
    \nabla (\SSSSmod f) \in L^{p(n+1)/n}_{\loc} (\R^{n+1}),
\end{equation}
and in fact $\nabla (\SSSSmod f) \in L^{p(n+1)/n} (\R^{n+1})$ if $\partial\Omega$ is unbounded, i.e., both $\Omega$ and $\overline{\Omega}^c$ are unbounded.

\subsection{Proof of \cref{prop -id+adjoint k is injective}}\label{sec:proof of +-id+k* are injective}

The proof of \cref{prop -id+adjoint k is injective} is presented in several steps. In \cref{step:0} we define the notation used throughout the proof, in \cref{step:1} we show that the single layer potential is constant in $\Omega$ and $\sigma$-a.e.\ on $\partial\Omega$, in \cref{step:2} we find a submean value property for the modulus of the single layer potential minus the constant found in \cref{step:1}, and from this we deduce in \cref{step:3} that the single layer potential is also constant in the complement of $\Omega$, with the same value. From this, in \cref{step:4} we conclude the proof of the injectivity.

\begin{subtask}[wide, labelindent=0pt,start=0]
\setlength{\itemsep}{5pt}
\item\label{step:0} Let $\varepsilon_2\in (0,\varepsilon_1]$, with $\varepsilon_1>0$ as in \cref{lemma:uniqueness neumann bvp}. During the proof we write $\Omega^+ \coloneqq \Omega $ and $\Omega^- \coloneqq \R^{n+1}\setminus \overline\Omega$. Let $f\in L^p (\sigma)$ satisfy $(-\frac{1}{2}Id + K^*) f = 0$ in $L^p(\sigma)$, equivalently,
\begin{equation}\label{-0.5 id + k* is zero ae}
\left(-\frac{1}{2}Id + K^*\right) f = 0\text{ }\sigma\text{-a.e.\ on }\partial\Omega.
\end{equation}
We want to show that $f=0$ $\sigma$-a.e.\ on $\partial\Omega$. To this end, we consider the modified single layer potential
$$
u\coloneqq \SSSSmod f,
$$
which is harmonic in $\R^{n+1}\setminus \partial\Omega$.

\item\label{step:1} In this step we prove
\begin{equation}\label{claim step 1}
    u\equiv c_0\text{ constant in }\Omega\text{ and }\sigma\text{-a.e.\ on }\partial\Omega.
\end{equation}
By \rf{jump formula gradient single layer interior} and \rf{-0.5 id + k* is zero ae} we have $\partial_{\nu_\Omega} u = (-\frac{1}{2}Id + K^*) f = 0$ $\sigma$-a.e.\ on $\partial\Omega$. Recall that $u:\Omega\to \R$ is harmonic and $\NN (\nabla u)\in L^p(\sigma)$ by \rf{lp norm of nontangential of gradient single layer potential}. All in all, the assumptions of \cref{lemma:uniqueness neumann bvp} are satisfied, and therefore we conclude that $u$ is constant in $\Omega$. By \cref{eq:single layer is continuous up to the boundary}, we conclude the proof of \rf{claim step 1}.

\item\label{step:2} In this step we prove the following submean value property
\begin{equation}\label{submean value property for |u-c|}
    |u(z)-c_0| \leq \avint_{B_r(z)} |u(x)-c_0|\, dm(x) \text{ for all }z\in \R^{n+1}\setminus\partial\Omega \text{ and }r>0,
\end{equation}
where $c_0$ is the constant in \rf{claim step 1}. Note that if $u$ were a continuous function in $\R^{n+1}$, this would immediately follow from the fact that $h$ is harmonic in $\overline{\Omega}^c$ and $h=0$ in $\overline\Omega$. However, due to the lack of continuity, a more careful argument is required.

We define $h(\cdot) \coloneqq u(\cdot)-c_0$ to shorten the notation. To prove \rf{submean value property for |u-c|}, we first see that it satisfies
\begin{equation}\label{|h| is weak subharmonic}
\int_{\R^{n+1}} \nabla |h|\nabla \varphi \, dm\leq 0 \text{ for all }\varphi \in C^\infty_c (\R^{n+1}) \text{ with } \varphi \geq 0,
\end{equation}
called weak subharmonic condition.

Note that, by \rf{single layer potential is loc sobolev with extra exponent} and since $|\nabla |h|| = |\nabla u|$ $m$-a.e.\ in $\R^{n+1}$, we have
\begin{equation}\label{grad |h| in L p+}
\nabla |h| \in L^{p(n+1)/n}_{\loc}(\R^{n+1}).
\end{equation}
Let us fix $B\supset \supp \varphi$ (not necessarily centered at $\partial\Omega$), and write $\Omega^-_B \coloneqq \Omega^-\cap B$ from now on. By $\supp \varphi\subset B$ and $h=u-c_0\equiv 0$ in $\Omega$ and $\sigma$-a.e.\ in $\partial\Omega$, we have
\begin{equation}\label{subharmonicity:from full space to just omega}
\int_{\R^{n+1}} \nabla |h|\nabla \varphi \, dm
= \int_{\Omega^-_B} \nabla |h|\nabla \varphi \, dm.
\end{equation}

As in \cite[(2.3.37)]{Hofmann2010} for the domain $\Omega^-_B$, for $\delta >0$ we define
$$
\chi_\delta (x) \coloneqq
\begin{cases}
    \begin{aligned}
    &1 &&\text{if } x\in \Omega^-_B \setminus U_\delta (\partial\Omega^-_B),\\
    &2\delta^{-1} \dist(x,\partial U_{\delta/2} (\partial\Omega^-_B)) &&\text{if } x\in \Omega^-_B \cap (U_\delta (\partial\Omega^-_B) \setminus U_{\delta/2} (\partial\Omega^-_B)),\\
    &0 &&\text{if } x\in U_{\delta/2} (\partial\Omega^-_B) \cup (\R^{n+1}\setminus \Omega^-_B).
    \end{aligned}
\end{cases}
$$
Moreover, we fix $\psi \in C^\infty_c (B_1(0))$ satisfying $\psi \geq 0$ and $\int_{B_1 (0)} \psi \, dm = 1$, and define
\begin{equation}\label{approx of identity to mollify}
    \psi_\varepsilon(\cdot) \coloneqq \frac{1}{\varepsilon^{n+1}} \psi \left(\frac{\cdot}{\varepsilon}\right), \text{ for }\varepsilon>0.  
\end{equation}

Since by Hölder's inequality and \rf{grad |h| in L p+} we have
$$
\int_{\Omega^-_B} \big| (1-\chi_\delta) \nabla |h|\nabla \varphi \big| \, dm
\lesssim \|\nabla |h|\|_{L^p(m|_{B})} m(U_\delta (\partial\Omega^-_B))^{1-\frac{1}{p}}
\lesssim \|\nabla |h|\|_{L^p(m|_{B})} \delta^{1-\frac{1}{p}} \overset{\delta\to 0}{\longrightarrow} 0,
$$
we get
\begin{equation}\label{subharmonicity:integral as a lim of delta}
\int_{\Omega^-_B} \nabla |h|\nabla \varphi \, dm
= \lim_{\delta \to 0} \int_{\Omega^-_B} \chi_\delta \nabla |h|\nabla \varphi \, dm.
\end{equation}
Also, since $\nabla |h|\in L^1 (m|_{B})$ by \rf{grad |h| in L p+}, for every $\delta>0$ we have
\begin{equation}\label{subharmonicity:integral in delta using convolution with epsilon}
\int_{\Omega^-_B} \chi_\delta \nabla |h|\nabla \varphi \, dm
= \lim_{\varepsilon\to 0} \int_{\Omega^-_B} \chi_\delta \nabla (|h|*\psi_\varepsilon)\nabla \varphi \, dm
\end{equation}
For every $\delta>0$ and every $0<\varepsilon<\delta/2$, by $\nabla (\chi_\delta \varphi)=\varphi \nabla \chi_\delta + \chi_\delta \nabla \varphi$ and $\divv (\chi_\delta \varphi \nabla(|h|*\psi_\varepsilon)) = \nabla (\chi_\delta \varphi) \nabla(|h|*\psi_\varepsilon) + \chi_\delta \varphi \Delta(|h|*\psi_\varepsilon)$, we have
$$
\begin{aligned}
\int_{\Omega^-_B} \chi_\delta \nabla (|h|*\psi_\varepsilon)\nabla \varphi \, dm
=& \int_{\Omega^-_B} \divv(\chi_\delta \varphi \nabla (|h|*\psi_\varepsilon)) \, dm\\
&- \int_{\Omega^-_B} \chi_\delta \varphi \Delta (|h|*\psi_\varepsilon) \, dm
- \int_{\Omega^-_B} \varphi \nabla (|h|*\psi_\varepsilon)\nabla \chi_\delta \, dm
\end{aligned}
$$
The first term on the right-hand side is zero by the classical divergence theorem. Also, we have that $-\Delta(|h|*\psi_\varepsilon)\leq 0$ in $\Omega^-\setminus U_\delta (\partial\Omega)$ for $\varepsilon<\delta/2$, because $|h|*\psi_\varepsilon$ is smooth and $(|h|*\psi_\varepsilon)(z)\leq m_{B(z,r)} (|h|*\psi_\varepsilon)$ for all $z\in \Omega^- \setminus U_\delta (\partial\Omega)$ and $0<r<\delta/2$, which follows from the fact that $|h|(z)\leq m_{B(z,r)} |h|$ for all $z\in \Omega^-$ and $0<r<\dist(z,\partial\Omega)$ (since $\Delta h=0$ in $\Omega^-$). Hence, we get
\begin{equation}\label{subharmonicity:inequality integral in delta and convolution epsilon}
\int_{\Omega^-_B} \chi_\delta \nabla (|h|*\psi_\varepsilon)\nabla \varphi \, dm
\leq - \int_{\Omega^-_B} \varphi \nabla (|h|*\psi_\varepsilon)\nabla \chi_\delta \, dm.
\end{equation}

From \rf{subharmonicity:integral as a lim of delta}, \rf{subharmonicity:integral in delta using convolution with epsilon}, \rf{subharmonicity:inequality integral in delta and convolution epsilon}, and that the latter integral in \rf{subharmonicity:inequality integral in delta and convolution epsilon} converges to $-\int_{\Omega^-_B} \varphi \nabla |h|\nabla \chi_\delta \, dm$ as $\varepsilon\to 0$ because $\nabla |h|\in L^1(m|_{B})$, we get
\begin{equation}\label{subharmonicity:inequality integral just in delta}
\int_{\Omega^-_B} \nabla |h|\nabla \varphi \, dm
\leq \lim_{\delta \to 0} -\int_{\Omega^-_B} \nabla \chi_\delta (\varphi \nabla |h|) \, dm.
\end{equation}
Rewriting the latter integral using the notation $\nu_\delta (x) \coloneqq - \nabla (\dist(x,\partial U_{\delta/2} (\partial\Omega^-_B)))$ and $\widetilde U_\delta (\partial\Omega^-_B) \coloneqq \Omega^-_B\cap (U_\delta(\partial\Omega^-_B)\setminus U_{\delta/2} (\partial\Omega^-_B))$, we now claim (and prove below) that
\begin{equation}\label{subharmonicity:claim limit in delta converges to boundary integral}
\lim_{\delta \to 0} \frac{2}{\delta}\int_{\widetilde U_\delta(\partial\Omega^-_B)} \langle\nu_\delta ,\varphi \nabla |h| \rangle \, dm 
=
\int_{\partial\Omega^-_B} \langle \nu_{\Omega^-_B}, (\varphi \nabla|h|)|_{\partial\Omega^-_B}\rangle \, d\HH^n.
\end{equation}

Before its proof, let us conclude the proof of \rf{|h| is weak subharmonic}. Using that $\varphi\equiv 0$ on $\partial B$ and the fact that $\partial_{\nu_{\Omega^-}} |h|\leq 0$ because $|h|\geq 0$ converges to zero nontangentially $\sigma$-a.e.\ at $\partial\Omega$ (by \rf{claim step 1} in \cref{step:1}), we have
$$
\int_{\partial\Omega^-_B} \langle \nu_{\Omega^-_B}, (\varphi \nabla|h|)|_{\partial\Omega^-_B}\rangle \, d\HH^n
= \int_{B\cap \partial\Omega} \varphi \partial_{\nu_{\Omega^-}} |h| \, d\sigma
\leq 0.
$$
This, together with \rf{subharmonicity:from full space to just omega}, \rf{subharmonicity:inequality integral just in delta} and \rf{subharmonicity:claim limit in delta converges to boundary integral}, concludes the proof of \rf{|h| is weak subharmonic}, in the absence of the justification of \rf{subharmonicity:claim limit in delta converges to boundary integral}.

Let us see the claim in \rf{subharmonicity:claim limit in delta converges to boundary integral}. We first need to check that
\begin{equation}\label{nontang of nabla |h| is in L1 and nontan limit exists}
\NN(\nabla |h|)\in L^p(\sigma)\text{ and }(\nabla |h|)|_{\partial\Omega} \text{ exists }\sigma\text{-a.e.\ on }\partial\Omega.
\end{equation}
The first condition is satisfied by $|\nabla |h||=|\nabla u|$ $m$-a.e.\ in $\Omega^-$ and \rf{lp norm of nontangential of gradient single layer potential}. Regarding the second condition, we want to see that for $\sigma$-a.e.\ $x\in \partial\Omega$ there exists a vector $v_x$ such that, for all $\varepsilon>0$ there is $\delta>0$ satisfying that
$$
z\in \Gamma^{\Omega^-}(x)\cap B_\delta(x)
\implies 
|\nabla |h| - v_x|<\varepsilon.
$$
By the jump formula \rf{jump formula gradient single layer exterior} and the assumption $(-\frac{1}{2}Id+K^*)f=0$, see \rf{-0.5 id + k* is zero ae}, we have $\partial_{\nu_{\Omega^-}} h = f$ for $\sigma$-a.e., in particular the pointwise nontangential limit $(\nabla h)|_{\partial\Omega}$ exists $\sigma$-a.e., which we denote by $w_x$ for such $x\in \partial\Omega$. If $w_x =0$, then $v_x=w_x$ does the job as $|\nabla |h|| = |\nabla h|$. Since $h\equiv 0$ $\sigma$-a.e., if $w_x\not =0$, then for $\delta>0$ small enough we have that either $h>0$ or $h<0$ in $\Gamma^{\Omega^-}(x)\cap B_\delta(x)$. If $\pm h>0$, then we take $v_x = \pm w_x$, and hence we have
$$
|\nabla |h|- \pm w_x|
=|\pm \nabla u \mp w_x|
=|\nabla u - w_x|<\varepsilon,
$$
where we used that $\nabla |h| = \frac{h}{|h|}\nabla u = \pm \nabla u$ in $\Gamma^{\Omega^-}(x)\cap B_\delta(x)$ (since we are in the case $\pm h >0$ there) and the fact that the pointwise nontangential limit exists $\sigma$-a.e.\ and is $(\nabla u)|_{\partial\Omega}=w_x$.

As argued in \cite[p.~2596, lines~6-8]{Hofmann2010}, the equality \rf{subharmonicity:claim limit in delta converges to boundary integral} directly holds if one replaces $\varphi \nabla |h|$ by $v\in C^{0,1} (\overline{\Omega^-_B})$. We adapt the end of Step I in the proof of \cite[Theorem 1.3.1]{MitreaX3-GeometricHarmonicAnalysis-Vol1} to our situation. By the density of $C^\infty_c (\R^{n+1})$ functions in $L^1 (\sigma|_{B\cap \partial\Omega})$, for any $\eta>0$ there exists $w \in C^\infty_c (\R^{n+1})$ such that
$$
\|(\nabla |h|)|_{\partial\Omega} - w\|_{L^1 (\sigma|_{B\cap \partial\Omega})}<\eta.
$$
First,
$$
\int_{\partial\Omega^-_B} |\langle\nu_{\Omega^-_B}, \varphi w-(\varphi \nabla |h|)|_{\partial\Omega}\rangle| \, d\HH^n
\lesssim \|(\nabla |h|)|_{\partial\Omega}-w\|_{L^1 (\sigma|_{B\cap\partial\Omega})} < \eta,
$$
second, as we already noted, since $\varphi w$ is in particular Lipschitz in $\overline{\Omega^-_B}$, the equality in \rf{subharmonicity:claim limit in delta converges to boundary integral} holds in this case and so we have
$$
\lim_{\delta \to 0} \frac{2}{\delta}\int_{\widetilde U_\delta(\partial\Omega^-_B)} \langle\nu_\delta , \varphi w \rangle \, dm 
=
\int_{\partial\Omega^-_B} \langle \nu_{\Omega^-_B}, \varphi w \rangle \, d\HH^n,
$$
and third, by \cref{lemma:2.3.25 in prop 2.12 in HMT}, for small enough $\delta>0$ so that $\varphi =0$ in $B\setminus (1-\delta) B$ we have
$$
\frac{2}{\delta} \int_{\widetilde U_\delta(\partial\Omega^-_B)} |\langle \nu_\delta, \varphi\nabla |h| - \varphi w\rangle| \, dm
\lesssim \|\NN^\delta (\nabla |h| - w)\|_{L^1 (\sigma|_{B\cap \partial\Omega})},
$$
recall the definition of $\NN^\delta$ in \rf{delta-nontangential maximal function}. Note that the second condition in \rf{nontang of nabla |h| is in L1 and nontan limit exists} implies that $\NN^\delta (\nabla |h| - w)(x) \to (\nabla |h|)|_{\partial\Omega}(x)-w(x)$ as $\delta\to 0$ for $\sigma$-a.e.\ $x\in\partial\Omega$. From this and $\NN(\nabla |h|)\in L^1(\sigma|_{B\cap\partial\Omega})$, see \rf{nontang of nabla |h| is in L1 and nontan limit exists}, by the dominated convergence theorem we have
$$
\lim_{\delta \to 0} \|\NN^\delta (\nabla |h| - w)\|_{L^1 (\sigma|_{B\cap \partial\Omega})} = \|(\nabla |h|)|_{\partial\Omega}-w\|_{L^1 (\sigma|_{B\cap \partial\Omega})} < \eta.
$$
All in all, we have
$$
\left|
\lim_{\delta \to 0} \frac{2}{\delta}\int_{\widetilde U_\delta(\partial\Omega^-_B)} \langle\nu_\delta ,\varphi \nabla |h| \rangle \, dm 
-
\int_{\partial\Omega^-_B} \langle \nu_{\Omega^-_B}, (\varphi \nabla|h|)|_{\partial\Omega}\rangle \, d\HH^n
\right|\\
\lesssim \eta,
$$
with uniform constant, and as $\eta>0$ is arbitrary, we conclude the claim in \rf{subharmonicity:claim limit in delta converges to boundary integral}, and therefore the proof of \rf{|h| is weak subharmonic} is now complete.

Let us see now how \rf{|h| is weak subharmonic} implies the submean value property \rf{submean value property for |u-c|}. With $\psi_\varepsilon$ as in \rf{approx of identity to mollify}, the function $|h|*\psi_\varepsilon$ is smooth and for $z\in \R^{n+1}$ there holds
\begin{equation}\label{|h|*psi_epsilon is subharmonic}
\begin{aligned}
-\Delta(|h|*\psi_\varepsilon)(z)
&=\sum_{i=1}^{n+1} -(\partial_{x_i}|h|*\partial_{x_i} \psi_\varepsilon) (z)
= -\int \nabla |h|(y) \nabla \psi_\varepsilon (z-y)\, dm(y)\\
&= \int \nabla |h|(y) \nabla (\psi_\varepsilon (z-y))\, dm(y)
\overset{\text{\rf{|h| is weak subharmonic}}}{\leq} 0,
\end{aligned}
\end{equation}
whence we get that $|h|*\psi_\varepsilon$ is subharmonic in $\R^{n+1}$ and in particular we get the submean value property
$$
(|h|*\psi_\varepsilon) (z) \leq \avint_{B_r(z)} (|h|*\psi_\varepsilon) (x)\, dm(x) \text{ for all }z\in \R^{n+1} \text{ and }r>0.
$$
For $z\in \R^{n+1}\setminus\partial\Omega$, the left-hand side converges to $|h(z)|$ because $h\in C^0 (\R^{n+1}\setminus\partial\Omega)$. On the other hand, since $\||h|*\psi_\varepsilon-|h|\|_{L^1(m|_{B_r(z)})}\to 0$ as $\varepsilon\to 0$ by \rf{single layer is in L 1+}, the right-hand side term converges to $\avint_{B_r(z)} |h(x)|\, dm(x)$ as $\varepsilon\to 0$. This concludes the proof of the submean value property \rf{submean value property for |u-c|}.\qed

\item\label{step:3} In this step we prove
\begin{equation}\label{claim step 3}
u\equiv c_0\text{ constant in }\Omega^-.
\end{equation}

For a fixed $z\in\Omega^-$, taking $r>2\dist(z,\partial\Omega)$ and $\xi\in\partial\Omega$ so that $|z-\xi|=\dist(z,\partial\Omega)$, we have $B(z,r)\subset B(\xi,2r)$ and therefore
$$
|u(z)-c_0| \overset{\text{\rf{submean value property for |u-c|}}}{\leq} \avint_{B(z,r)} |u(x)-c_0|\, dm(x)
\lesssim \avint_{B(\xi,2r)} |u(x)-c_0|\, dm(x).
$$
Note that $c_0$ is the mean of $u$ over any ball inside $\Omega$. So, taking $\widetilde B$ the interior corkscrew ball of $B(\xi,2r)$, we have that $\widetilde B\subset B(\xi,2r)\cap\Omega$ with $r_{\widetilde B}\approx r$, and therefore we get that the last term above is
$$
\avint_{B(\xi,2r)} |u(x)-c_0|\, dm(x)
= \avint_{B(\xi,2r)} |u(x)-m_{\widetilde B} u|\, dm(x).
$$
Adding $\pm m_{B(\xi,2r)} u$, by the classical Poincaré inequality since $u\in W^{1,1}_{\loc}(\R^{n+1})$ (see \rf{single layer is in L 1+} and \rf{grad |h| in L p+}), Hölder's inequality and the estimates \rf{extra regularity from boundary to interior} and \rf{lp norm of nontangential of gradient single layer potential}, the latter term is controlled by
$$
\begin{aligned}
\avint_{B(\xi,2r)} |u(x)-m_{\widetilde B} u|\, dm(x)
&\lesssim \avint_{B(\xi,2r)} |u(x)-m_{B(\xi,2r)} u|\, dm(x)
\lesssim r\avint_{B(\xi,2r)} |\nabla u(x)| \, dm(x).\\
&\lesssim r^{1-n/p} \|\nabla u\|_{L^{p(n+1)/n}(\Omega)}
\overset{\text{\rf{extra regularity from boundary to interior}},\text{\rf{lp norm of nontangential of gradient single layer potential}}}{\lesssim} r^{1-n/p} \|f\|_{L^p(\sigma)}
\end{aligned}
$$
All in all, for any $z\in \Omega^-$ and any $r>2\dist(z,\partial\Omega)$ we get
\begin{equation}\label{bound of subharmonic function of a fixed point and any radius}
|u(z)-c_0|\lesssim r^{1-n/p} \|f\|_{L^p(\sigma)}.
\end{equation}

In the case $2n/(n+1)\leq p<n$, letting $r\to \infty$ we deduce that $u(z)=c_0$, as claimed. For $n\leq p< n+\varepsilon_2$, we will use the following lemma to see that \rf{bound of subharmonic function of a fixed point and any radius} implies a slightly better bound, as in \rf{eq:proof uniqueness neumann step 1}.

\begin{lemma}\label{Holder cont at the boundary for subharmonic functions}
    Let $D \subset\R^{n+1}$ be a domain satisfying the exterior corkscrew condition with constant $M>1$. There are $\alpha=\alpha(M,n)>0$ and $C=C(M,n)\geq 1$ such that for every nonnegative smooth subharmonic function $u$ in $\R^{n+1}$, all $\xi\in \partial\Omega$ and all $r>0$, if $u\equiv 0$ in $\R^{n+1}\setminus U_{r/(4M)^{k_0}} (D)$ for some integer $k_0\geq 1$, then
    \begin{equation}\label{eq:Holder cont at the boundary for subharmonic functions}
    u(x)\leq C \left(\frac{|x-\xi|}{r}\right)^\alpha \sup_{B_r(\xi)} u, \text{ for all }x\in B(\xi,r)\setminus \overline{B(\xi,r/(4M)^{k_0+1})}.
    \end{equation}
\end{lemma}

This is well known for harmonic functions vanishing continuously on $\partial\Omega$, even for more general domains. However, to keep track of the parameters, the proof is written in full detail at the end of this section. Assume the lemma to be true for the moment.

Let us fix $z\in \Omega^-$, and take $r>2\dist(z,\partial\Omega)$ and $\xi\in\partial\Omega$ so that $\dist(z,\partial\Omega) = |z-\xi|$. With $\psi_\varepsilon$ as in \rf{approx of identity to mollify} for $\varepsilon>0$, as we already saw in \rf{|h|*psi_epsilon is subharmonic}, $|h|*\psi_\varepsilon$ is nonnegative, smooth and subharmonic in $\R^{n+1}$, and identically zero in $\R^{n+1}\setminus U_\varepsilon (\Omega^-)$. So, by \cref{Holder cont at the boundary for subharmonic functions} there is $\alpha_2 >0$ such that if $0<\varepsilon\ll r$, then
$$
(|h|*\psi_\varepsilon) (z)\overset{\text{\rf{eq:Holder cont at the boundary for subharmonic functions}}}{\lesssim} \left(\frac{|z-\xi|}{r}\right)^{\alpha_2} \sup_{B(\xi,r)} (|h|*\psi_\varepsilon).
$$
Also, for $x\in B(\xi,r)$ we have
$$
(|h|*\psi_\varepsilon) (x) = \int_{B(x,\varepsilon)} |h|(y) \psi_\varepsilon (x-y)\, dm(y)
\overset{\text{\rf{bound of subharmonic function of a fixed point and any radius}}}{\lesssim} r^{1-n/p} \|f\|_{L^p(\sigma)},
$$
as long as $\varepsilon<r$. All in all, and as $(|h|*\psi_\varepsilon) (z)\to |h(z)| = |u(z)-c_0|$ as $\varepsilon\to 0$ because $h\in C^0 (\Omega^-)$, we conclude
$$
|u(z)-c_0|\lesssim |z-\xi|^{\alpha_2} r^{1-\frac{n}{p}-\alpha_2} \|f\|_{L^p(\sigma)}.
$$
This goes to zero as $r\to \infty$ if $1<p<n/(1-\alpha_2)$. The parameter $\varepsilon_2>0$ is taken so that $n+\varepsilon_2 =\min\{n+\varepsilon_1,n/(1-\alpha_2)\}$. We conclude that $u\equiv c_0$ in $\Omega^-$.\qed

\item\label{step:4} Let us finish the proof by seeing that $f\equiv 0$ in $L^p (\sigma)$. By \rf{jump formula gradient single layer exterior} and since $u\equiv c_0$ in $\Omega^-$ (see \rf{claim step 3} in \cref{step:3}), we have $(\frac{1}{2}Id + K^*)f = \partial_{\nu_{\overline{\Omega}^c}} u =0$. On the other hand, recall that $(-\frac{1}{2}Id + K^*)f=0$ by assumption, see \rf{-0.5 id + k* is zero ae}. We conclude that
$$
f = \left(\frac{1}{2}Id + K^*\right)f - \left(-\frac{1}{2}Id + K^*\right)f = 0\quad\sigma\text{-a.e.\ on }\partial\Omega,
$$
as claimed.\qed
\end{subtask}

We conclude this section with the proof of \cref{Holder cont at the boundary for subharmonic functions}. We adapt the standard argument to show the Hölder continuity of harmonic functions that vanish continuously on the boundary of a domain satisfying the exterior corkscrew condition.

\begin{proof}[Proof of \cref{Holder cont at the boundary for subharmonic functions}]
    We first claim that there is $c_1=c_1(M,n) \in (0,1)$ such that for a fixed $\xi\in\partial\Omega$ and $s>0$, there holds
    \begin{equation}\label{fast decayment harmonic function}
    \sup_{B(\xi,s/(4M))} u
    \leq (1-c_1)\sup_{B(\xi,s)} u, \text{ as long as } u\equiv 0 \text{ in }\R^{n+1}\setminus U_{s/(2M)}(D).
    \end{equation}
    Indeed, by the exterior corkscrew condition there is $B(A_s(\xi),s/M)\subset B(\xi,s)\setminus \overline D$, and therefore there exists $\rho \in (s/(2M),s)$ such that
    \begin{equation}\label{lower bound n hausdorff measure}
    \HH^n (\partial B(\xi,\rho)\setminus D)\geq \HH^n (\partial B(\xi,\rho)\cap B(A_s(\xi),s/(2M)))\gtrsim_{M,n} \rho^n.
    \end{equation}
    We are using $2M>2$ to ensure that $u\equiv 0$ in $B(A_s(\xi),s/(2M))\subset \R^{n+1}\setminus U_{s/(2M)}(D)$. Let $H_u$ be the harmonic extension in $B(\xi,\rho)$ of $u|_{\partial B(\xi,\rho)}$. For $x\in B(\xi,\rho)$, using the Poisson kernel and that $u\equiv 0$ in $B(A_s(\xi),s/(2M))\subset \R^{n+1}\setminus U_{s/(2M)}(D)$, we write $H_u$ as
    $$
    \begin{aligned}
    H_u(x) &= \frac{\rho^2-|x-\xi|^2}{w_n \rho} \int_{\partial B(\xi,\rho)} \frac{u(z)}{|x-z|^{n+1}} \, d\HH^n (z)\\
    &=\frac{\rho^2-|x-\xi|^2}{w_n \rho} \int_{\partial B(\xi,\rho)\setminus B(A_s(\xi),s/(2M))} \frac{u(z)}{|x-z|^{n+1}} \, d\HH^n (z),
    \end{aligned}
    $$
    where $w_n = \HH^n(\partial B(0,1))$. We bound $H_u (x)$ for $x\in B(\xi,\rho/2)$. Bounding $u$ by its supremum in $B(\xi,s)\supset B(\xi,\rho)$, we have
    $$
    \begin{aligned}
    H_u (x)\leq&\, \frac{\rho^2-|x-\xi|^2}{w_n \rho} \int_{\partial B(\xi,\rho)\setminus B(A_s(\xi),s/(2M))} \frac{\sup_{B(\xi,s)}u}{|x-z|^{n+1}} \, d\HH^n (z)\\
    =&\, \frac{\rho^2-|x-\xi|^2}{w_n \rho} \int_{\partial B(\xi,\rho)} \frac{\sup_{B(\xi,s)}u}{|x-z|^{n+1}} \, d\HH^n (z)\\
    &-\frac{\rho^2-|x-\xi|^2}{w_n \rho} \int_{\partial B(\xi,\rho) \cap B(A_s(\xi),s/(2M))} \frac{\sup_{B(\xi,s)}u}{|x-z|^{n+1}} \, d\HH^n (z)\\
    =& \left(1-\frac{\rho^2-|x-\xi|^2}{w_n \rho} \int_{\partial B(\xi,\rho) \cap B(A_s(\xi),s/(2M))} \frac{1}{|x-z|^{n+1}} \, d\HH^n (z)\right)\sup_{B(\xi,s)} u, 
    \end{aligned}
    $$
    where we used in the last equality that the harmonic extension (using the Poisson kernel) of a constant function is the constant itself. Now, using that $|x-\xi|\leq \rho/2$ and $|x-z|\leq 3\rho/2$ if $x\in B(\xi,\rho/2)$ and $z\in \partial B(\xi,\rho)$, we get
    \begin{multline*}
    \frac{\rho^2-|x-\xi|^2}{w_n \rho} \int_{\partial B(\xi,\rho) \cap B(A_s(\xi),s/(2M))} \frac{1}{|x-z|^{n+1}} \, d\HH^n (z)\\
    \geq \frac{2^{n-1}}{3^n w_n \rho^n} \HH^n (\partial B(\xi,\rho) \cap B(A_s(\xi),s/(2M)))
    \overset{\text{\rf{lower bound n hausdorff measure}}}{\geq} c_1 \in (0,1), 
    \end{multline*}
    for some $c_1 = c_1 (M,n)\in (0,1)$. All in all, by the inclusion, the maximum principle, and the estimates above, we get
    $$
    \sup_{B(\xi,s/(4M))} u
    \leq \sup_{B(\xi,\rho/2)} u
    \leq \sup_{B(\xi,\rho/2)} H_u
    \leq (1-c_1)\sup_{B(\xi,s)} u,
    $$
    provided that $u\equiv 0$ in $\R^{n+1}\setminus U_{s/(2M)}(D)$, as claimed in \rf{fast decayment harmonic function}.
    
    As we are assuming $u\equiv 0$ in $\R^{n+1}\setminus U_{r/(4M)^{k_0}}(D)$, iterating \rf{fast decayment harmonic function} we get
    $$
    \sup_{B(\xi,r/(4M)^k)} u \leq (1-c_1)^k \sup_{B(\xi,r)} u \text{ for all }0\leq k\leq k_0.
    $$
    This readily proves the lemma.
\end{proof}

\section{Uniqueness of the solution of the Dirichlet problem}\label{sec:uniqueness of dirichlet problem}

Given a Wiener regular\footnote{This ensures that the Green function is well-defined. For instance, domains with the exterior corkscrew condition or whose boundary is ADR are Wiener regular. In fact, satisfying one of these conditions at small scales is sufficient.} domain $\Omega$ (not necessarily bounded) and $x\in \Omega$, we denote by $g_x=g_x^\Omega$ the harmonic Green function of $\Omega$ with pole at $x$, and by $\omega^x = \omega_\Omega^x$ the harmonic measure of $\Omega$ with pole at $x$. If the domain $\Omega$ is ADR and $\omega^x$ is absolutely continuous with respect to $\sigma$ ($\omega^x \ll \sigma$ for short) for some (and hence for all) $x\in \Omega$, then the Radon-Nikodym derivative $d\omega^x/d\sigma$ (also known as the Poisson kernel) exists, with
$$
\omega^x (E) = \int_E \frac{d\omega^x}{d\sigma} \, d\sigma, \text{ for any Borel set }E\subset \R^{n+1},
$$
see \cite[Theorem 2.17]{Mattila1995} for instance.

This section is dedicated to the uniqueness of the solution of the $L^p$ Dirichlet problem in \rf{dirichlet bvp} under the assumption that the harmonic measure satisfies the $p^\prime$-reverse Hölder inequality.

\begin{propo}\label{reverse holder implies uniqueness dirichlet}
    Let $\Omega\subset\R^{n+1}$ be a CAD, $x\in\Omega$, $p\in (1,\infty)$, and let $p^\prime$ so that $1/p+1/p^\prime=1$. Assume $\omega_\Omega^x\ll\sigma$ and that there exists $\gamma>0$ such that
    \begin{equation}\label{reverse holder ineq}
        \left(\avint_{B(\xi,r)} \left|\frac{d\omega_\Omega^x}{d\sigma}\right|^{p^\prime}\, d\sigma \right)^{1/{p^\prime}} \lesssim \frac{\omega^x (B(\xi,r))}{\sigma(B(\xi,r))}.
    \end{equation}
    for every $\xi\in\partial\Omega$ and $0<r\leq \gamma$. Then the zero function is the unique solution of the homogeneous Dirichlet problem
    \begin{equation}\label{homogeneous dirichlet bvp}
    \begin{cases*}
    	\Delta u = 0 \text{ in } \Omega\\
        \NN u \in L^p(\sigma)\\
    	u|_{\partial\Omega}^{\rm nt} = 0 \text{ }\sigma\text{-a.e.\ on }\partial\Omega.
    \end{cases*}
    \end{equation}
\end{propo}

This is well known when the domain is assumed to be bounded, see for instance \cite[Theorem 5.7.7]{MitreaX3-GeometricHarmonicAnalysis-Vol3}. For completeness, we provide a detailed proof in the general case, where $\Omega$ is not necessarily bounded. Its proof is a direct combination of \cref{decayment of green function,uniqueness dirichlet just from decayment of green} below. A quick inspection of the proof of \cref{decayment of green function} below reveals that the NTA condition is only used for balls with sufficiently small radii.

We state the extrapolation of the solvability of $(D_p)$ with uniqueness of \rf{dirichlet bvp}.

\begin{coro}\label{extrapolation dirichlet problem}
    Let $\Omega\subset\R^{n+1}$ be a CAD. If $(D_{p_0})$ is solvable for some $1<p_0<\infty$, then there exists $\varepsilon>0$ such that $(D_p)$ is solvable for all $p\in (p_0-\varepsilon,\infty)$ and the solution of $(D_p)$ is the unique solution of \rf{dirichlet bvp}.
    \begin{proof}
        Using that the harmonic measure is doubling in NTA domains (see \cite[Lemma 4.9]{Jerison1982}\footnote{\label{original JK for bounded domains}This is originally for bounded domains, but the same continues to hold for unbounded domains.}), this follows by the well-known equivalence between the solvability of the Dirichlet problem and the reverse Hölder for the harmonic measure (see \cite[Proposition 2.14]{MourgoglouPoggiTolsa2025}\footnote{The solvability of the Dirichlet problem in \cite[Definition 1.1]{MourgoglouPoggiTolsa2025} is stated for $C_c (\partial\Omega)$ functions, instead of $L^p(\sigma)$. However, it is well known that this is equivalent to the $L^p$ solvability definition of $(D_p)$ in \rf{dirichlet bvp} and \rf{lp norm nontang of dirichlet solution}, by a density argument.} for instance), Gehring's lemma (see \cite[Theorem 6.38]{Giaquinta2012} for instance), and \cref{reverse holder implies uniqueness dirichlet}.
    \end{proof}
\end{coro}

We now turn to the first step in the proof of \cref{reverse holder implies uniqueness dirichlet}.

\begin{lemma}\label{decayment of green function}
    Let $\Omega\subset\R^{n+1}$ be a CAD domain, $x\in\Omega$, $p\in (1,\infty)$, and let $p^\prime$ so that $1/p+1/p^\prime=1$. Assume that there exists $\gamma>0$ such that \rf{reverse holder ineq} holds for every $\xi\in\partial\Omega$ and $0<r\leq \gamma$. Then for any $\delta < \min\{\diam(\partial\Omega),\dist(x,\partial\Omega)\}/100$ there holds
    $$
    \begin{aligned}
    \|\NN^\delta g_x^\Omega\|_{L^{p^\prime}(\sigma)} &\lesssim \gamma^{-n/p} \delta,\\
    \|\NN^\delta (\nabla g_x^\Omega)\|_{L^{p^\prime}(\sigma)} &\lesssim \gamma^{-n/p},
    \end{aligned}
    $$
    where the involved constant depends on $p$, the CAD character of $\Omega$, the constant in \rf{reverse holder ineq} (and the aperture on the definition of $\NN^\delta$ in \rf{delta-nontangential maximal function}).
    \begin{proof}
        We define the Hardy-Littlewood maximal function
        $$
        M_\sigma \omega^x (\xi) \coloneqq \sup_{r>0} \frac{\omega^x (B_r (\xi))}{\sigma (B_r (\xi))}
        = \sup_{r>0} \avint_{B_r(\xi)} \frac{d\omega^x}{d\sigma}\, d\sigma, \quad \xi\in\partial\Omega.
        $$
        
        Let us fix $\xi\in\partial\Omega$. For any $y\in \Gamma(\xi)\cap \overline{B_{2\delta}(\xi)}$, in particular, $|y-\xi|< (1+\alpha)\dist(y,\partial\Omega)$, we have that $y$ is essentially a corkscrew point at scale $2|y-\xi|$ in the sense that both $y$ and $A_{2|y-\xi|} (\xi)$ (the corkscrew point of $\xi$ at scale $2|y-\xi|$) belong to $B_{2|y-\xi|}(\xi)$ and are uniformly far from the boundary with distance $\gtrsim 2|y-\xi|$. Thus, by Harnack inequality and the relation in \cite[Lemma 4.8]{Jerison1982}\textsuperscript{\ref{original JK for bounded domains}} between the Green function and the harmonic measure in NTA domains, we have
        $$
        g_x (y)\approx g_x (A_{2|y-\xi|}(\xi)) \lesssim \frac{\omega^x (B(\xi,2|y-\xi|))}{|y-\xi|^{n-1}}.
        $$
        Therefore, from this and the definition of $M_\sigma \omega^x$, we have
        $$
        \NN^\delta g_x (\xi) = \sup_{y\in \Gamma(\xi)\cap \overline{B_{2\delta} (\xi)}} g_x (y)
        \lesssim \sup_{y\in \Gamma(\xi)\cap \overline{B_{2\delta} (\xi)}} |y-\xi| \frac{\omega^x (B(\xi,2|y-\xi|))}{|y-\xi|^n}
        \lesssim \delta M_\sigma \omega^x (\xi).
        $$
        So, from this and the $L^{p^\prime}$-bound of the Hardy-Littlewood maximal function $M_\sigma \omega^x$ we have
        $$
        \int_{\partial\Omega} |\NN^\delta g_x|^{p^\prime} \, d\sigma \lesssim \delta^{p^\prime} \int_{\partial\Omega} |M_\sigma \omega^x|^{p^\prime} \, d\sigma
        \lesssim_{p^\prime} \delta^{p^\prime} \int_{\partial\Omega} \left|\frac{d\omega^x}{d\sigma}\right|^{p^\prime} \, d\sigma.
        $$
        Given by the $5R$-covering theorem, let $\{B_k\}_{k\geq 1}$ be a subfamily of $\{B(\xi,\gamma)\}_{\xi\in\partial\Omega}$ such that $\partial\Omega \subset \bigcup_{k\geq 1} B_k$ and $\{B_k/5\}_{k\geq 1}$ is pairwise disjoint. So, \rf{reverse holder ineq} applies for each $B_k$ since $r(B_k)=\gamma$. Also, since $\{B_k/5\}_{k\geq 1}$ is a pairwise disjoint family and each ball has the same radius, we also have that the family $\{B_k\}_{k\geq 1}$ has finite overlap. Thus, we have
        $$
        \begin{aligned}
        \int_{\partial\Omega} \left|\frac{d\omega^x}{d\sigma}\right|^{p^\prime} \, d\sigma
        &\leq \sum_{k\geq 1} \int_{B_k} \left|\frac{d\omega^x}{d\sigma}\right|^{p^\prime} \, d\sigma
        \approx \gamma^n \sum_{k\geq 1} \avint_{B_k} \left|\frac{d\omega^x}{d\sigma}\right|^{p^\prime} \, d\sigma
        \lesssim \gamma^n \sum_{k\geq 1} \left(\frac{\omega^x (B_k)}{\sigma(B_k)}\right)^{p^\prime}\\
        &\approx \gamma^{n(1-p^\prime)} \sum_{k\geq 1} \omega^x (B_k)^{p^\prime}
        \leq \gamma^{n(1-p^\prime)} \sum_{k\geq 1} \omega^x (B_k)
        \lesssim \gamma^{n(1-p^\prime)}.
        \end{aligned}
        $$
        All in all, $\|\NN^\delta g_x\|_{L^{p^\prime}(\sigma)} \lesssim \gamma^{-n/p} \delta$.

        The $L^{p^\prime}(\sigma)$-norm of the nontangential maximal operator of $\nabla g_x$ follows from the same computations above, using that $|\nabla g_x (\cdot)| \lesssim g_x (\cdot)/\dist(\cdot,\partial\Omega)$ far from the pole. Indeed, we now have
        $$
        \NN^\delta (\nabla g_x) (\xi) = \sup_{y\in \Gamma(\xi)\cap \overline{B_{2\delta} (\xi)}} |\nabla g_x (y)|
        \lesssim \sup_{y\in \Gamma(\xi)\cap \overline{B_{2\delta} (\xi)}} \frac{g_x (y)}{\dist(y,\partial\Omega)}
        \lesssim \sup_{y\in \Gamma(\xi)\cap \overline{B_{2\delta} (\xi)}} \frac{g_x (y)}{|y-\xi|},
        $$
        and the other computations that are needed have already been done.
    \end{proof}
\end{lemma}

We now proceed to the second step in the proof of \cref{reverse holder implies uniqueness dirichlet}.

\begin{lemma}\label{uniqueness dirichlet just from decayment of green}
    Let $\Omega\subset \R^{n+1}$ be an ADR domain and let $u$ be a solution of \rf{homogeneous dirichlet bvp} with $1<p<\infty$. Take $p^\prime$ so that $1/p+1/p^{\prime}=1$, and fix $z\in\Omega$ and $g_z=g_z^\Omega$. If $\|\NN^\delta g_z\|_{L^{p^\prime}(\sigma)}\lesssim \delta$ for small enough $\delta>0$, then $u\equiv 0$ in $\Omega$.
        \begin{proof}
            First note that $\|\NN^\delta g_x\|_{L^{p^\prime}(\sigma)}\lesssim \delta$ holds for every $x\in\Omega$ with constant depending on $x\in \Omega$ and provided $\delta$ is small enough also depending on $x\in\Omega$, because by the symmetry of the harmonic Green function we have
            $$
            g_x (y) = g_y (x) \approx_{x,z} g_y (z) = g_z (y).
            $$
        
            Let us fix $x\in \Omega$ and $R\geq 2\dist(x,\partial\Omega)$, and define $\Omega_R \coloneqq \Omega \cap B_R(x)$. We remark that the auxiliary domain $\Omega_R$ and its Green function are only necessary if the domain is unbounded, either with bounded or unbounded boundary. If the domain $\Omega$ is bounded, one can directly remove the dependence on $R$ by taking $R=2\diam(\Omega)$ on the computations below.
            
            We will make use of three Green's functions: $g^\Omega_x$, $g^{\Omega_R}_x$ and $g^{B_R(x)}_x$ are the Green functions with pole at $x$ of $\Omega$, $\Omega_R$ and $B_R(x)$ respectively. Moreover, it holds that $g^{\Omega_R}_x \leq g^\Omega_x$ and $g^{\Omega_R}_x \leq g^{B_R(x)}_x$ in $\Omega_R\setminus \{x\}$.
        
            Let us fix a dimensional constant $C_1\geq 1$. For $\delta< \dist(x,\partial\Omega_R)/(4C_1)$, as in \cite[(7.1.10) and (7.1.11)]{Hofmann2010}, let $\psi_\delta^\Omega\in C^\infty (\Omega)$ such that $0\leq \psi_\delta^\Omega \leq 1$, $\psi_\delta^\Omega =1$ in $\Omega\setminus U_{C_1\delta} (\partial\Omega)$, $\psi_\delta^\Omega = 0$ in $U_{C_1\delta/2} (\partial\Omega) \cap \Omega$, and $|\partial^\alpha \psi_\delta^\Omega| \lesssim_\alpha (C_1\delta)^{-|\alpha|}$ for all multi-indices $\alpha$. Similarly for $B_R(x)$ (with $\delta$ instead of $C_1\delta$), define $\psi_\delta^{B_R(x)}\in C^\infty_c (B_R(x))$ such that $0\leq \psi_\delta^{B_R(x)} \leq 1$, $\psi_\delta^{B_R(x)} =1$ in $B_R(x)\setminus U_{\delta} (\partial B_R(x))$, $\psi_\delta^{B_R(x)} = 0$ in $U_{\delta/2} (\partial B_R(x)) \cap B_R(x)$, and $|\partial^\alpha \psi_\delta^{B_R(x)}| \lesssim_\alpha \delta^{-|\alpha|}$ for all multi-indices $\alpha$. Finally, we define
            $$
            \psi_\delta \coloneqq \psi_\delta^\Omega \psi_\delta^{B_R(x)} \in C^\infty_c (\Omega_R),
            $$
            which in particular satisfies $0\leq \psi_\delta \leq 1$, $|\partial^\alpha \psi_\delta| \lesssim_{\alpha,C_1} \delta^{-|\alpha|}$ for all multi-indices $\alpha$, and moreover
            $$
            \dist(y,\partial\Omega_R) \geq \delta/2\text{ for all }y\in\supp\nabla \psi_\delta.
            $$
            
            As in \cite[(7.1.12)]{Hofmann2010}, we have
            $$
            u(x) = u(x)\psi_\delta(x) 
            = \int_{\Omega_R} \langle\nabla g^{\Omega_R}_x,\nabla (u\psi_\delta)\rangle \, dm
            = 2\int_{\Omega_R} \langle\nabla g^{\Omega_R}_x,\nabla \psi_\delta\rangle u \, dm
            +\int_{\Omega_R} g^{\Omega_R}_x u \Delta\psi_\delta\, dm.
            $$
            Using that $|\nabla g^{\Omega_R}_x(\cdot)|\lesssim g^{\Omega_R}_x(\cdot)/\dist(\cdot,\partial\Omega_R) \lesssim g^{\Omega_R}_x(\cdot)/\delta$ in $\supp \nabla \psi_\delta$, and
            $$
            \supp\nabla\psi_\delta \subset (U_{C_1\delta} (\partial\Omega)\cap B_R(x)) \cup ((B_R(x)\setminus \overline{B_{R-\delta} (x)})\setminus U_{C_1\delta} (\partial\Omega)),
            $$
            we have
            \begin{equation}\label{eq:u bounded by I delta and II delta}
            \begin{aligned}
            |u(x)|&\lesssim \frac{1}{\delta^2} \int_{\supp\nabla\psi_\delta} g_x^{\Omega_R} |u|\, dm|_{\Omega}\\
            &\leq \frac{1}{\delta^2} \int_{U_{C_1\delta} (\partial\Omega)\cap B_R(x)} g_x^{\Omega_R} |u|\, dm|_{\Omega}
            +\frac{1}{\delta^2} \int_{(B_R(x)\setminus \overline{B_{R-\delta} (x)})\setminus U_{C_1\delta} (\partial\Omega)} g_x^{\Omega_R} |u|\, dm|_{\Omega}\\
            &\eqqcolon \rom{1}_\delta + \rom{2}_\delta.
            \end{aligned}
            \end{equation}
        
            Let us study the term $\rom{1}_\delta$. Using that $g^{\Omega_R} \leq g^\Omega$ in $\Omega_R\setminus \{x\}$, \cref{lemma:2.3.25 in prop 2.12 in HMT} and Hölder's inequality respectively, we have
            \begin{equation}\label{uniqueness:bound of I delta}
            \rom{1}_\delta \leq \frac{1}{\delta^2} \int_{U_{C_1\delta} (\partial\Omega)\cap B_R(x)} g^{\Omega} |u|\, dm|_{\Omega}
            \lesssim \frac{1}{\delta} \int_{\partial\Omega} \NN^{C_1\delta} (g^\Omega|u|)\, d\sigma
            \leq \frac{1}{\delta} \|\NN^{C_1\delta} g^\Omega\|_{L^{p^\prime}(\sigma)} \|\NN^{C_1\delta} u\|_{L^p(\sigma)}.
            \end{equation}    
        
            We now turn to the term $\rom{2}_\delta$, which is directly zero if the domain is bounded and $R\geq 2\diam(\Omega)$. Before, recall that
            $$
            g^{B_R(x)}_x (y) = \frac{1}{\kappa_n (n-1)}\left(\frac{1}{|y-x|^{n-1}} - \frac{1}{R^{n-1}}\right),
            $$
            where $\kappa_n$ is the surface area of the unit sphere $\mathbb S^n\subset\R^{n+1}$. So, for $y\in B_R(x)\setminus \overline{B_{R-\delta}(x)}$ we get
            $$
            g^{B_R(x)}_x (y) \lesssim \frac{1}{|R-\delta|^{n-1}} - \frac{1}{R^{n-1}} \lesssim \frac{\delta}{R^n}.
            $$
            Hence, using that $g^{\Omega_R}\leq g^{B_R}$ in $\Omega_R\setminus \{x\}$ and this, we have
            $$
                \rom{2}_\delta
                \lesssim 
                \frac{1}{R^n} \frac{1}{\delta} \int_{(B_R(x)\setminus \overline{B_{R-\delta} (x)})\setminus U_{C_1\delta} (\partial\Omega)} |u|\, dm|_{\Omega}.
            $$
            For shortness we write $F_\delta \coloneqq (B_R(x)\setminus \overline{B_{R-\delta} (x)})\setminus U_{C_1\delta} (\partial\Omega)$. Let $W_\Omega$ denote the family of Whitney cubes of $\Omega$, as in \cref{Whitney decomposition lemma}, and
            $$
            W_\Omega^\delta\coloneqq \{Q\in W_\Omega : Q\cap F_\delta\not=\emptyset\}.
            $$
            Then we have
            $$
            \rom{2}_\delta
            \lesssim \frac{1}{R^n} \frac{1}{\delta} \int_{F_\delta} |u|\, dm|_{\Omega}
            = \frac{1}{R^n} \frac{1}{\delta} \sum_{Q\in W_\Omega^\delta} \int_{Q\cap F_\delta} |u|\, dm.
            $$
            Since $m(Q\cap F_\delta)\lesssim \frac{\delta}{\ell(Q)} m(Q)=\delta \ell(Q)^n$ for each $Q\in W_\Omega^\delta$, we have
            $$
            \rom{2}_\delta
            \lesssim \frac{1}{R^n} \frac{1}{\delta} \sum_{Q\in W_\Omega^\delta} \int_{Q\cap F_\delta} |u|\, dm
            \lesssim \frac{1}{R^n} \sum_{Q\in W_\Omega^\delta} \ell(Q)^n \sup_Q |u|.
            $$
            Note that having fixed $C_1\geq 1$ big enough depending on the dimension, then every $Q\in W_\Omega^\delta$ satisfies $\ell(Q)\geq 100\delta$, and therefore $3Q\cap \partial B_R(x)\not = \emptyset$ with $\HH^n(3Q\cap \partial B_R(x))\approx \ell(Q)^n$. So, we get
            \begin{equation}\label{eq:bound of II delta in terms of the sum}
            \rom{2}_\delta
            \lesssim \frac{1}{R^n} \sum_{Q\in W_\Omega^\delta} \ell(Q)^n \sup_Q |u|
            \lesssim \frac{1}{R^n} \sum_{Q\in W_\Omega} \HH^n(3Q\cap \partial B_R(x)) \sup_Q |u|,
            \end{equation}
            uniformly in $\delta>0$, where we also used in the last step that $W_\Omega^\delta \subset W_\Omega$.
            
            Let us bound the sum on the right-hand side. Defining
            $$
            f(y)\coloneqq \sum_{Q\in W_\Omega} \characteristic_{3Q\cap\partial B_R(x)} (y) \sup_Q |u|,\quad y\in\Omega,
            $$
            we have
            \begin{equation}\label{eq:sum in whitney cubes bounded by integral function f}
            \begin{aligned}
            \sum_{Q\in W_\Omega} \HH^n(3Q\cap \partial B_R(x)) \sup_Q |u|
            &=\sum_{Q\in W_\Omega} \int_{3Q\cap\partial B_R(x)} \sup_Q |u| \, d\HH^n\\
            &\leq \sum_{Q\in W_\Omega} \int_{3Q\cap\partial B_R(x)} f \, d\HH^n
            \lesssim \int_{\partial B_R (x)} f\, d\HH^n,
            \end{aligned}
            \end{equation}
            where we used the finite overlap of $\{3Q\}_{Q\in W_\Omega}$ in the last step. Applying \cref{lemma:interior integral adr to boundary} to $f$ (with $\xi_x \in\partial\Omega$ such that $\dist(x,\partial\Omega)=|x-\xi_x|$, $C_2\geq 1$ is a fixed big enough dimensional constant so that $\supp f\subset B(\xi,C_2 R)$, and with aperture $\beta_1$), the latter term is bounded by
            \begin{equation}\label{eq:integral function f bounded by integral maximal function f}
            \int_{\partial B_R (x)} f\, d\HH^n
            \lesssim 
            \int_{2B(\xi_x,C_2 R)\cap \partial\Omega} \NN_{\beta_1} f \, d\sigma.
            \end{equation}
            Given $\xi\in\partial\Omega$, let
            $$
            W_\Omega(\xi) \coloneqq \{Q\in W_\Omega : 3Q \cap \Gamma_{\beta_1}(\xi) \not = \emptyset\}.
            $$
            Using this definition, for any $y\in \Gamma_{\beta_1} (\xi)$, the sum in the definition of $f(y)$ runs over the cubes in $W_\Omega (\xi)$ (instead of $W_\Omega$), which in particular implies
            \begin{equation}\label{eq:maximal function of f in term of cubes touching cone}
            \NN_{\beta_1} f(\xi) 
            = \sup_{y\in \Gamma_{\beta_1} (\xi)} f(y)
            = \sup_{y\in \Gamma_{\beta_1} (\xi)} \sum_{Q\in W_\Omega(\xi)} \characteristic_{3Q\cap\partial B_R(x)} (y) \sup_Q |u|.
            \end{equation}
            By the relation $\ell(Q)\approx \dist(Q,\partial\Omega)$ for any $Q\in W_\Omega$, it is clear that there exists an aperture $\beta_2 \geq \beta_1$ (depending on $\beta_1$ and the dimension) such that
            $$
            \bigcup_{Q\in W_\Omega(\xi)} 3Q \subset \Gamma_{\beta_2} (\xi).
            $$
            From \rf{eq:maximal function of f in term of cubes touching cone}, this and the finite overlap of $\{3Q\}_{Q\in W_\Omega}$, we get
            $$
            \NN_{\beta_1} f(\xi) \lesssim \NN_{\beta_2} u,
            $$
            Therefore, from \rf{eq:sum in whitney cubes bounded by integral function f}, \rf{eq:integral function f bounded by integral maximal function f} and this, we conclude
            \begin{equation}\label{eq:final bound of the sum in whitney cubes of II delta}
            \sum_{Q\in W_\Omega} \HH^n(3Q\cap \partial B_R(x)) \sup_Q |u|
            \lesssim \int_{2B(\xi_x,C_2 R)\cap \partial\Omega} \NN_{\beta_2} u \, d\sigma.
            \end{equation}
            
            All in all, we have
            $$
            \rom{2}_\delta
            \overset{\text{\rf{eq:bound of II delta in terms of the sum}}}{\lesssim} \frac{1}{R^n} \sum_{Q\in W_\Omega} \HH^n(3Q\cap \partial B_R(x)) \sup_Q |u|
            \overset{\text{\rf{eq:final bound of the sum in whitney cubes of II delta}}}{\lesssim} \frac{1}{R^n} \int_{2B(\xi_x,C_2 R)\cap \partial\Omega} \NN_{\beta_2} u \, d\sigma,
            $$
            and by Hölder's inequality (and \rf{eq:Lp of nontangential maximal function does not depend on aperture}) we get
            \begin{equation}\label{uniqueness:bound of II delta}
            \rom{2}_\delta \lesssim R^{-n/p} \|\NN u\|_{L^p(\sigma)},
            \end{equation}
            and we conclude the control of the term $\rom{2}_\delta$.
        
            Collecting estimates \rf{uniqueness:bound of I delta} and \rf{uniqueness:bound of II delta} in \rf{eq:u bounded by I delta and II delta}, we obtain, for all $R\geq 2\dist(x,\partial\Omega)$ and all $\delta < \dist(x,\partial\Omega_R)/(4C_1)$,
            $$
            |u(x)|\lesssim \frac{1}{\delta} \|\NN^{C_1\delta} g^\Omega\|_{L^{p^\prime}(\sigma)} \|\NN^{C_1\delta} u\|_{L^p(\sigma)} + R^{-n/p} \|\NN u\|_{L^p(\sigma)}.
            $$
            First, for $\delta>0$ small enough (depending on $x$), by \cref{decayment of green function} we have
            $$
            \frac{1}{\delta} \|\NN^{C_1\delta} g^\Omega\|_{L^{p^\prime}(\sigma)} \|\NN^{C_1\delta} u\|_{L^p(\sigma)}
            \lesssim \|\NN^{C_1\delta} u\|_{L^p(\sigma)}.
            $$
            Since $u|_{\partial\Omega}^{\rm nt}=0$ $\sigma$-a.e.\ on $\partial\Omega$ implies $\NN^{C_1\delta} u(x)\to 0$ as $\delta \to 0$ for $\sigma$-a.e.\ $x\in\partial\Omega$, and $\NN u \in L^p(\sigma)$, by the dominated convergence theorem we have $\|\NN^{C_1\delta} u\|_{L^p(\sigma)}\to 0$ as $\delta \to 0$. Second, as we are assuming $\|\NN u\|_{L^p(\sigma)}<\infty$ we have $R^{-n/p} \|\NN u\|_{L^p(\sigma)}\to 0$ as $R\to \infty$. That is, $u(x)=0$. Since the point $x\in \Omega$ is arbitrary, we conclude that $u$ is identically zero in $\Omega$.
        \end{proof}
\end{lemma}

\section{The Dirichlet and Neumann problems}\label{sec:proof of dirichlet and neumann problems}

In this section we solve the Dirichlet and Neumann problems, stated in \cref{thm:lp dirichlet problem,thm:lp neumann problem} respectively.

\begin{proof}[Proof of \cref{thm:lp dirichlet problem}]
    Let $\varepsilon_D >0$ be so that $n/(n-1)-\varepsilon_D$ is the Hölder conjugate exponent of $n+\varepsilon_2$, with $\varepsilon_2>0$ as in \cref{coro:+id + K* is injective}. Fix $n/(n-1)-\varepsilon_D <p_0\leq 2n/(n+1)$, and let $\delta_0 = \delta_0 (p_0,n,\text{CAD})>0$ be given by \cref{reduction:all equivalent}.
    
    By \cref{reduction:all equivalent,coro:+id + K* is injective} we have that $\left(\frac{1}{2}Id + K\right)^{-1} : L^{p_0}(\sigma)\to L^{p_0}(\sigma)$ is a bounded linear operator, as claimed in the theorem. Hence, $\left(\frac{1}{2}Id + K\right)^{-1} f \in L^{p_0}(\sigma)$ with $\|\left(\frac{1}{2}Id + K\right)^{-1} f\|_{L^{p_0}(\sigma)}\lesssim \|f\|_{L^{p_0}(\sigma)}$, where the involved constant does not depend on $f\in L^{p_0}(\sigma)$. It is clear then that the function $u = \DD\left(\left(\frac{1}{2}Id + K\right)^{-1} f\right)$ in \rf{eq:solution of lp dirichlet problem} is harmonic in $\Omega$. Moreover, by the jump formula \rf{jump formula double layer interior}, we have
    $$
    u|_{\partial\Omega}^{\rm nt}(x) = f(x), \text{ for }\sigma\text{-a.e.\ }x\in\partial\Omega,
    $$
    whence we conclude that it solves the Dirichlet problem in \rf{dirichlet bvp} with $p_0$ instead of $p$. Finally, by \rf{lp norm of nontangential of double layer potential} and as $\left(\frac{1}{2}Id + K\right)^{-1}$ is bounded in $L^{p_0}(\sigma)$, we conclude
    $$
    \|\NN u\|_{L^{p_0}(\sigma)}
    \lesssim \left\|\left(\frac{1}{2}Id + K\right)^{-1} f\right\|_{L^{p_0}(\sigma)}
    \lesssim \|f\|_{L^{p_0}(\sigma)},
    $$
    as claimed in \rf{lp norm nontang of dirichlet solution}. That is, $(D_{p_0})$ is solvable with solution $u$.

    By \cref{extrapolation dirichlet problem}, there is $\varepsilon>0$ such that $(D_p)$ is solvable for all $p\in (p_0-\varepsilon,\infty)$, and the solution of $(D_p)$ is the unique solution of \rf{dirichlet bvp}.  
\end{proof}

\begin{proof}[Proof of \cref{thm:lp neumann problem}]
    Let $\varepsilon_N=\varepsilon_2$ with $\varepsilon_2>0$ as in \cref{prop -id+adjoint k is injective}. We fix $2n/(n+1)\leq p < n+\varepsilon_N$, and let $\delta_0 = \delta_0 (p,n,\text{CAD})$ be given by \cref{reduction:all equivalent}.
    
    By \cref{reduction:all equivalent,prop -id+adjoint k is injective} we have that $\left(-\frac{1}{2}Id + K^*\right)^{-1} : L^p(\sigma)\to L^p(\sigma)$ is a bounded linear operator. Hence, $\left(-\frac{1}{2}Id + K^*\right)^{-1} f \in L^p(\sigma)$ with $\|\left(-\frac{1}{2}Id + K^*\right)^{-1} f\|_{L^p(\sigma)}\lesssim \|f\|_{L^p(\sigma)}$, where the involved constant does not depend on $f\in L^p(\sigma)$. It is clear then that the function $u = \SSSSmod\left(\left(-\frac{1}{2}Id + K^*\right)^{-1} f\right)$ in \rf{eq:solution of lp neumann problem} is harmonic in $\Omega$. Moreover, by the jump formula \rf{jump formula gradient single layer interior} we have
    $$
    \partial_{\nu_\Omega} u (x) = f(x), \text{ for }\sigma\text{-a.e.\ }x\in\partial\Omega, 
    $$
    whence we conclude that it solves the Neumann problem in \rf{neumann bvp}. Finally, by \rf{lp norm of nontangential of gradient single layer potential} and as $\left(-\frac{1}{2}Id + K^*\right)^{-1}$ is bounded in $L^p(\sigma)$, we conclude
    $$
    \|\NN (\nabla u)\|_{L^p(\sigma)}
    \lesssim \left\|\left(-\frac{1}{2}Id + K^*\right)^{-1} f\right\|_{L^p(\sigma)}
    \lesssim \|f\|_{L^p(\sigma)},
    $$
    as claimed in \rf{lp norm nontang of neumann solution}. That is, $(N_p)$ is solvable with solution $u$, and by \cref{lemma:uniqueness neumann bvp}, it is the unique (modulo constant) solution of \rf{neumann bvp}.
\end{proof}

\begin{rem}
    It is straightforward to verify that any solution $u$ of the Dirichlet problem \rf{dirichlet bvp} (respectively, Neumann problem \rf{neumann bvp}) satisfies $\|\NN u\|_{L^p (\sigma)} \geq \|f\|_{L^p (\sigma)}$ (respectively, $\|\NN (\nabla u)\|_{L^p (\sigma)} \geq \|f\|_{L^p (\sigma)}$).
\end{rem}

\vv

\begin{acknowledgements}
Supported by Generalitat de Catalunya’s Agency for Management of University and Research Grants (AGAUR) (2021 FI\_B 00637 and 2023 FI-3 00151) and by the Xavier Tolsa Domènech Donation Project from the 2019 Rei Jaume I Award. Also partially supported by 2021-SGR-00071 (AGAUR, Catalonia) and by the Spanish State Research Agency (AEI) projects PID2020-114167GB-I00 and PID2021-125021NAI00. The author thanks CERCA Programme/Generalitat de Catalunya for institutional support.
\end{acknowledgements}


\vv

\vv

\bibliographystyle{alpha}
\bibliography{references.bib} 

\end{document}